\documentclass{article}
\usepackage{graphicx} 

\usepackage{algpseudocode}
\usepackage{algorithm}
\usepackage{amsmath}
\usepackage{amsthm}
\usepackage{amsfonts}
\usepackage{graphicx}
\usepackage{geometry}
\usepackage{bm}
\usepackage{todonotes}
\usepackage{mathtools}
\usepackage{epstopdf}
\usepackage{float}
\usepackage{placeins}
\usepackage{ulem}
\usepackage{amssymb}
\usepackage{latexsym}
\usepackage{url}
\usepackage[percent]{overpic}
\usepackage{xcolor}
\definecolor{newcolor}{rgb}{.8,.349,.1}
\usepackage{mathrsfs} 
\usepackage{hyperref}
\usepackage{cleveref}
\usepackage{subcaption}
\usepackage{amsthm}
\usepackage{algpseudocode}
\usepackage{algorithm}
\usepackage{amsfonts}
\usepackage{graphicx}
\usepackage{geometry}
\usepackage{bm}
\usepackage{mathtools}
\usepackage{comment}

\usepackage{xcolor}
\usepackage{amsmath}
\usepackage{amssymb}
\usepackage{amsthm}
\usepackage{amsfonts}
\usepackage{graphicx}
\usepackage{geometry}
\usepackage{bm}
\usepackage[T1]{fontenc} 							
\usepackage[utf8]{inputenc} 							

\usepackage{mathtools}
\newtheorem{pro}{Problem}

\usepackage[percent]{overpic}
\graphicspath{ {./Figures/} }
\newcommand{\pad}[2]{\frac{\partial #1}{\partial #2}}

\newcommand{\td}[2] {\frac{ {\rm d} #1}{ {\rm d} #2}}
\newcommand{\BB}{\text{\Large$\mathbb{B}$} }

\newcommand{\vphi}{\varphi}

\renewcommand{\vec}[1]{\mathbf{#1}}

\title{High order ghost-FEM for incompressible Navier-Stokes equations on moving domains}
\author{Hridya Dilip, Clarissa Astuto, Armando Coco and Giovanni Russo}
\date{July 2025}

\begin{document}

\maketitle

\begin{abstract}
    We develop a new numerical technique for approximating solutions of the Navier–Stokes equations on moving domains. The method aims at simulating an incompressible fluid past an object whose motion is assigned {\it a priori} using a level-set function. The proposed approach relies on a space discretization based on the ghost finite element method (ghost-FEM), which allows computations on unfitted meshes and avoids costly remeshing as the domain evolves in time. Time integration is performed using an IMplicit-EXplicit (IMEX) scheme to address the nonlinearity of the convective term, ensuring  high-order accuracy for incompressible flows. The error introduced by the geometrical approximation is handled using the Shifted Boundary Method, which allows higher order approximations of boundary conditions on unfitted meshes. Dirichlet boundary conditions are imposed weakly by means of Nitsche’s method. The associated stabilization parameter is chosen  by solving a generalized eigenvalue problem, ensuring stability and accuracy of the numerical scheme. We present a series of numerical experiments designed to validate the accuracy  of the proposed method, as well as comparisons with established benchmark problems involving moving boundaries. 
\end{abstract}
\section{Introduction}
A wide variety of fluid dynamics phenomena can be modeled using the unsteady incompressible Navier–Stokes equations. These equations capture the behavior of flows in environments such as the atmosphere, oceans, rivers, and lakes, as well as in industrial systems such as turbines, aircraft engines, and vehicles \cite{vidal2003thermodynamics}.  Another well-studied problem is the flow of surfactants around a moving trap \cite{tomlinson2023unsteady,tomlinson2023laminar,hansbo2016cut,LEE200643,CiCP-31-707}, that motivated the development of the method proposed in this paper.
Aqueous surfactants play a key role in various biological, biochemical, and industrial processes. They influence foam stability \cite{brown1990foam}, wettability, and coating flows \cite{valentini1991role}, and are commonly applied as sprays to enhance the effectiveness of foliar-applied agrochemicals by improving pesticide penetration into plant foliage \cite{knoche1991performance}. 
The influence of surfactants on the deformation of drops and bubbles \cite{Raudino20168574,CiCP-31-707}, and on the surrounding flow field, remains an area of active research, not only for its relevance to polymer and emulsion technologies, but also for the insights it offers into the complex physics of viscous interfacial flows \cite{velankar2001effect}. Surfactant molecules act as a buffer between the fluid phases across an interface, modifying intermolecular forces and reducing surface tension in proportion to the local surfactant concentration \cite{adamson1967physical}. For this reason, computing two-phase flows with surfactant has been challenging due to the associated interface dynamics, surface convection–diffusion of surfactants, fluid dynamics of the underlying Stokes or Navier–Stokes flows, and the coupling among them via the surface tension \cite{10.1093/imamat/66.1.55,XU20125897,hansbo2016cut,ASTUTO2023111880}.


From a numerical perspective, the Navier–Stokes equations have been extensively studied \cite{LANGTANGEN20021125}, using finite differences \cite{van1986second,COCO2020109623,ASTUTO2023111880,astuto2023multiscale}, finite element \cite{hughes1986new,fortin1981old} and finite volume \cite{leveque2002finite,Osher} methods, and in \cite{GANESAN20123685} a finite element scheme based on coupled arbitrary Lagrangian–Eulerian (for the velocity) and Lagrangian approach (for the concentration of surfactants) is developed for the computation of interface flows with soluble surfactants. Many numerical schemes use a level-set function to compute two-phase flows with  surfactants \cite{XU20125897,hansbo2016cut}, and they are dissolved in the bulk fluid but also exist in adsorbed form on the interface separating two immiscible fluids. 

Among the time discretizations, we have $\theta-$method, such as Crank-Nicolson, or the projection method, proposed by Chorin in \cite{CHORIN1997118,BROWN2001464}. It consists of a fractional step, where the velocity field is advanced in time by discretizing the momentum equation and ignoring the pressure term, obtaining a fictitious velocity that is not divergence-free. This velocity is then corrected by the pressure term by enforcing the divergence-free condition.
Other approaches in time are the IMEX schemes \cite{boscarino2024implicit}, semi-implicit strategies designed to deal with nonlinear terms and high stiffness. They have been employed in a wide variety of applications, including the design of Asymptotic Preserving methods for the
inviscid \cite{boscheri2020second} and viscous \cite{boscheri2021efficient,decaria2021embedded} compressible flows. IMEX schemes offer a suitable and flexible approach when aiming for higher-order time integration \cite{astuto2023self}, and here we illustrate their performance in the context of Navier-Stokes equations and moving domains. 

For the spatial discretization we use the ghost-FEM, a recent unfitted boundary method first introduced in \cite{astuto2024nodal,astuto2024comparison}. {The same space discretization has been extended to the numerical solution of biological network formation in a leaf-shaped domain in  \cite{astuto2024self} and to Multiscale Poisson-Nernst-Planck system \cite{astuto2025asymptotic,astuto2025standard}. In \cite{dilip2025multigrid}, the authors present multigrid methods for solving elliptic partial differential equations on arbitrary domains using the same ghost-FEM discretization.} In this method, a rectangular region discretized by a regular Cartesian mesh is arbitrarily intersected by the domain $\Omega$. 
The active mesh $\Omega_h$ is the subset of the mesh that intersects $\Omega$. The finite element method is based on a variational formulation over $\Omega_h$, with stabilization terms to handle instabilities coming from small cut elements. In this study, the incompressible Navier–Stokes equations are solved in $\mathbb Q_2/\mathbb Q_1$ integrated FEM \cite{ern2021finite}. More importantly, the pressure variable is treated in a fully implicit manner, computing the solution of the pressure/velocity variables in a monolithic way. 

Despite the abundance of unfitted finite element schemes available in the literature, such as CutFEM in \cite{burman2015cutfem} and AgFEM in \cite{BADIA2018533,BADIA202360}, the ghost-FEM is of interest for the following reason: when dealing with arbitrary domains, it is common to face small cut cells problem. While other schemes rely on stabilization techniques, the ghost-FEM employs a very simple algorithm based on snapping-back-to-grid approach to overcome the ill-conditioning issue, and it is straightforward to implement numerically. To the best of our knowledge, the use of IMEX schemes in conjunction with immersed boundary methods remains largely unexplored. Therefore, this work represents a significant contribution to this active area of research.

The paper is organized as follows. In Section~\ref{sec:model}, we present the mathematical model and its variational formulation. Section~\ref{sec:numerics} describes the numerical tools employed in our simulations, starting from the level-set approach used to implicitly define the computational domain, the shifted boundary method adopted to handle geometric errors when prescribing boundary conditions on the approximated domain, the choice of the optimal penalization parameter for Dirichlet conditions, and the Aslam technique to extrapolate the solution outside the fluid domain. In Section~\ref{sec:results}, we first verify the expected order of accuracy in space and time of the proposed numerical scheme, and then present results for both steady and moving domains, including comparisons with established benchmark problems. We point out that the motion of the domain is assigned {\it a priori} (one-way coupling). 
Finally, we draw some conclusions.

\section{Model and variational formulation}
\label{sec:model}
\begin{figure}
\begin{subfigure}[t]{0.49\textwidth}
    \centering
        \begin{overpic}[width=0.8\textwidth]{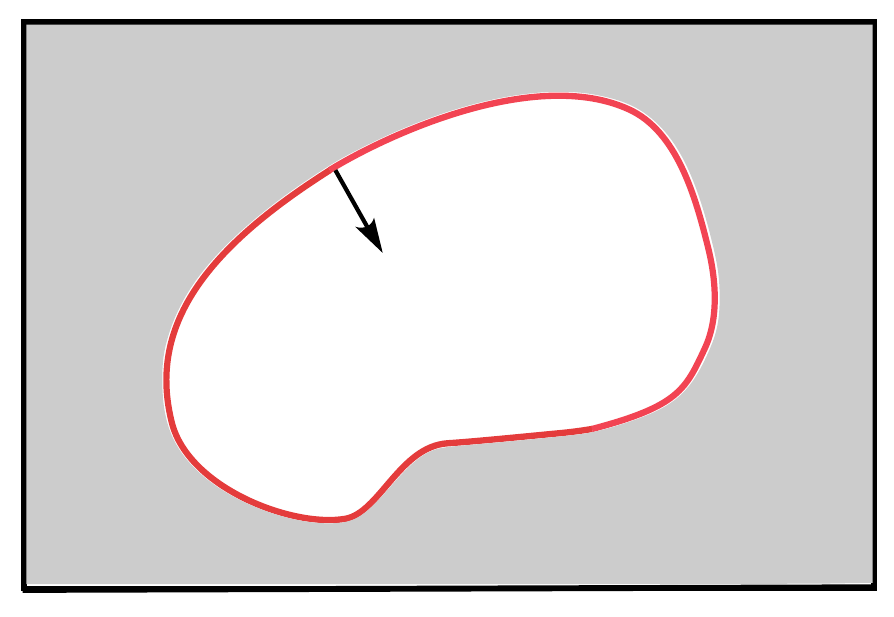}
\put(68,43){\LARGE $\Gamma_c$}
\put(11,43){\LARGE $\Omega$}
\put(99,24){\LARGE $\Gamma_s$}
\put(26,22){\LARGE $\Omega_c$}
\put(41,45){\Large ${\bf n}$}
\end{overpic}
        \caption{\textit{Domain representation: $\Omega$ is the domain, given as the difference between the rectangle $R$ and $\Omega_c$; $\Gamma_c$ is the internal boundary, $\Gamma_s $ the external squared boundary and ${\bf n}$ the outer normal vector to $\Gamma$.}}
        \label{fig:domain_representation}
    \end{subfigure}    
    \begin{subfigure}[t]{0.49\textwidth}
    \centering    
    \begin{overpic}[width=0.6\textwidth]{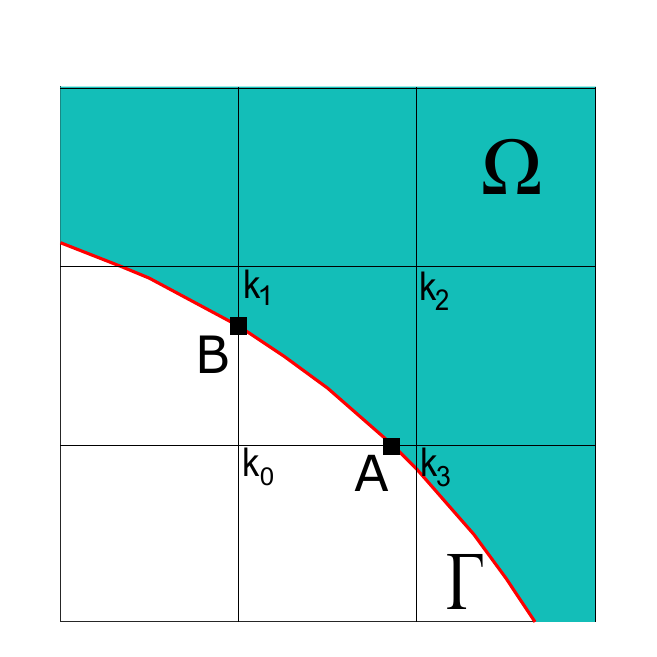}
    \end{overpic}
\caption{\textit{Intersection of the boundary $\Gamma$ with the grid: see Algorithm~\ref{alg:ab} for the computation of the points A and B.}}
\label{fig:ab}
    \end{subfigure} 
\caption{\textit{Geometric features of the computational domain.}}   
\label{fig:domain}
\end{figure}
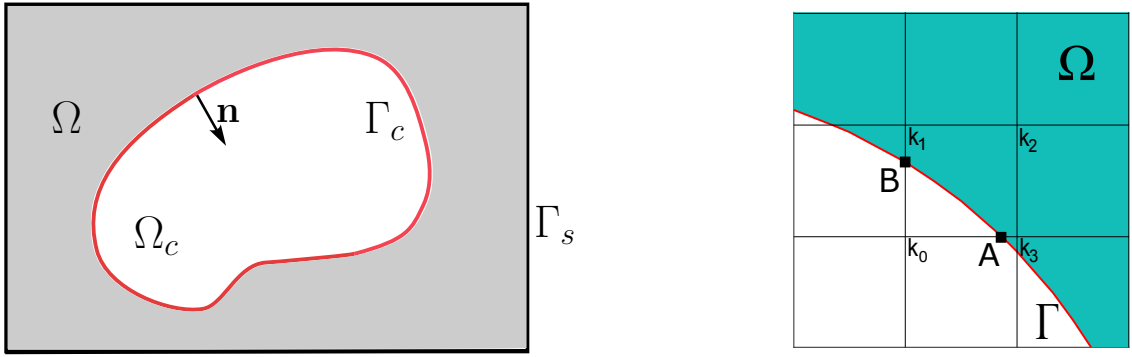

The incompressible Navier-Stokes equations are
given by
\begin{subequations}
\label{eq:NS}
\begin{align}\label{eq:NS1}
\pad{\textbf{u}}{t} + \textbf{u}\cdot \nabla \textbf{u} + \nabla p & = \frac{1}{Re}\Delta \textbf{u} + {\bf f}  \qquad {\rm in } \, \, \Omega(t)\\   \label{eq:NS2}
\nabla \cdot \textbf{u} & = 0 \\
{\bf u } & = {\bf g} \qquad {\rm on } \, \, \Gamma,
\end{align}
\end{subequations}
where ${\bf u} = {\bf u}({\bf x},t) $ is the fluid velocity, $p = p({\bf x},t)$ is the pressure, $Re$ is the Reynolds number and ${\bf f} = {\bf f}({\bf x},t)$ the external forces. 
The domain $\Omega(t)$ is a bounded domain in $\mathbb R^2$, and we assume that it is the difference between a fixed rectangular region $R$ of dimension $L_x \times L_y$ and a cavity $\Omega_c(t) \subseteq R$, that represents a solid object (see Fig.~\ref{fig:domain_representation}).
From now on, we omit the time dependence of the domain and of the boundary when it is clear from the context.
We denote by $\Gamma_c = \partial \Omega_c$, $\Gamma_s = \partial R$, and $\Gamma = \partial \Omega$. Then, $\Gamma = \Gamma_c \cup \Gamma_s$.
The system is closed by imposing Dirichlet boundary conditions on $\Gamma$. The governing equations are based on the fundamental physical principles of conservation and can be derived from the conservation of mass and momentum. Eq.~\eqref{eq:NS1} represents the conservation of momentum, while Eq.~\eqref{eq:NS2} is the incompressibility constraint forcing the conservation of mass. The main challenge in discretizing the system lies in the  structure of the resulting algebraic problem. The saddle-point nature of the formulation, where the pressure acts as a Lagrange multiplier to enforce the incompressibility constraint, induces a coupling between velocity and pressure even at the discrete level \cite{brezzi2012mixed,boffi2013mixed,quarteroni2000factorization}. In Section~\ref{section:saddle} we show the technique that we employ to deal with this feature. 

Here we consider the variational formulation of Eqs.~\eqref{eq:NS}. We introduce the two spaces ${\bf V}$ and $Q$ defined as
\begin{equation}
\label{eq:FEMspaces}
     \mathbf{V} = [H^1(\Omega)]^2, \qquad Q = \biggl\{ q\in H^1(\Omega) : \int_\Omega q\,{ \rm d}\Omega = 0\biggr\}.
\end{equation}
Multiplying \eqref{eq:NS1} by a test function ${\bf w} \in {\bf V}$,  \eqref{eq:NS2} by $q \in Q$, integrating over $\Omega$ and applying the divergence theorem, we obtain
\begin{subequations}
\begin{align}
&\int_\Omega \pad {\bf u} t {\bf w}\,{ \rm d}\Omega + \int_\Omega \left( {\bf u} \cdot \nabla \right){\bf u} \,{\bf w} \, { \rm d}\Omega 
 +  \int_\Omega \nabla {\bf u} \cdot \nabla {\bf w} \, { \rm d}\Omega  -  \int_{\Gamma } \left(\nabla {\bf u} \cdot {\bf n} \right) {\bf w} \, { \rm d}s -  \int_{\Gamma }  {\bf u} \left( \nabla {\bf w} \cdot {\bf n} \right) \, { \rm d}s \\ & + \lambda\int_{\Gamma} {\bf u} \, {\bf w} \, { \rm d}s - \int_\Omega  p \, \nabla \cdot {\bf w} \, { \rm d}\Omega  + \int_{\Gamma } \left(p \, {\bf n} \right) {\bf w} \,{\rm d} s  = \int_\Omega f\, {\bf w}\,{ \rm d}\Omega +\lambda \int_{\Gamma} {\bf g}\, {\bf w} \, { \rm d}s -  \int_{\Gamma }  {\bf g} \left( \nabla {\bf w} \cdot {\bf n} \right) \, { \rm d}s \\ 
&\int_\Omega \nabla \cdot {\bf u} \, q \,{ \rm d}\Omega -  \int_{\Gamma } \left(q \, {\bf n} \right) {\bf u} \,{\rm d} s   = -  \int_{\Gamma } \left(q \, {\bf n} \right) {\bf g} \,{\rm d} s.
\end{align}
\end{subequations}
It is a well-established fact that as the penalization $\lambda\to\infty$ will force zero Dirichlet conditions
\cite{Courant,HilbertCourant}.

We denote by $(\cdot,\cdot)_{L^2(\Omega)}$ and $(\cdot,\cdot)_{L^2(\Gamma)}$ the scalar products in $L^2(\Omega)$ and $L^2(\Gamma)$, respectively. With this notation, the variational formulation of the problem reads as follows.
\begin{pro}\label{pro:variational}
Given ${\bf f}, {\bf g} \in L^2(\Omega)$, find ${\bf u}(t) \in {\bf V}$ and $p(t) \in Q$, for almost every $t\in(0,T)$, such that
\begin{subequations}
\label{eq_var_formulation}
\begin{align*} 
&\left( \pad {\bf u} t,{\bf w} \right)_{L^2(\Omega)} + \left( \left( {\bf u} \cdot \nabla \right){\bf u} , {\bf w} \right) _{L^2(\Omega)}
 + \left( \nabla {\bf u},  \nabla {\bf w} \right)_{L^2(\Omega)} - \left(\nabla {\bf u} \cdot {\bf n} , {\bf w} \right)_{L^2(\Gamma)} -  \left({\bf u} , \nabla  {\bf w}\cdot {\bf n} \right)_{L^2(\Gamma)} \\ & + \lambda ({\bf u},{\bf w})_{L^2(\Gamma)}  - \left(  p, \nabla \cdot {\bf w} \right)_{L^2(\Omega)} + \left(p {\bf n} , {\bf w} \right)_{L^2(\Gamma)} =  \left( {\bf f}, {\bf w}\right)_{L^2(\Gamma)} + \lambda ({\bf g},{\bf w})_{L^2(\Omega)}  -  \left({\bf g} , \nabla  {\bf w}\cdot {\bf n} \right)_{L^2(\Gamma)} \qquad &\forall {\bf w} \in {\bf V}  \\ 
&\left( \nabla \cdot {\bf u} , q\right)_{L^2(\Omega)}  - \left(q {\bf n} , {\bf u} \right)_{L^2(\Gamma)} = - \left(q {\bf n} , {\bf g} \right)_{L^2(\Gamma)} \qquad &\forall q \in Q.
\end{align*}
\end{subequations}
\end{pro}

\subsection{Level-set representation of the domain}
The domain $\Omega$ is implicitly defined by a level set function $\phi(x,y)$ such that~\cite{Osher2002,Russo2000,book:72748, Sussman1994}:
\begin{align}
	\Omega = \{(x,y): \phi(x,y) < 0\}, \qquad
	\Gamma_c = \{(x,y): \phi(x,y) = 0\}.
\end{align}
The unit normal vector on $\Gamma_c$ pointing outward from $\Omega$ can be computed as 
\begin{equation}\label{eq:normal} {\bf n} = \frac{\nabla \phi }{|\nabla \phi|}.
\end{equation}

\section{Numerical approximation}
\label{sec:numerics}
In this section, we describe the space and time discretization for the problem~\eqref{pro:variational}. For the spatial discretization, we follow the strategy in \cite{astuto2024nodal,astuto2024comparison}, a recently developed ghost-FEM method. Time discretization is performed by an IMplicit EXplicit numerical scheme
\cite[Section 1.3]{boscarino2024implicit}, \cite{jin1999efficient}.

\subsection{Spatial Discretization}
The rectangular region $R=(x_\text{min},x_\text{min}+L_x) \times (y_\text{min},y_\text{min}+L_y)$ is discretized by a uniform mesh of squares ({\it cells}) of size $h=L_x/N_x = L_y/N_y$, where $N_x,N_y \in \mathbb N$ are the number of cells in $x$ and $y$ direction, respectively.  
The set of cells is denoted by $\mathcal C$, with $\# \mathcal C = N_x N_y$. A cell $K \in \mathcal C$ is classified as {\it active cell} if $K \cap  \Omega \neq \emptyset$, or {\it inactive cell} otherwise.
The union of all active cells is denoted by $\Omega_\text{act}$.
An active cell $K$ can be {\it internal cell} if $K \subseteq \Omega$ or {\it cut cell} otherwise (namely $K \cap \Gamma \neq \emptyset)$. 

\begin{figure}[h]
    \centering
    \begin{subfigure}{0.45\textwidth}
\begin{overpic}[width=\textwidth]{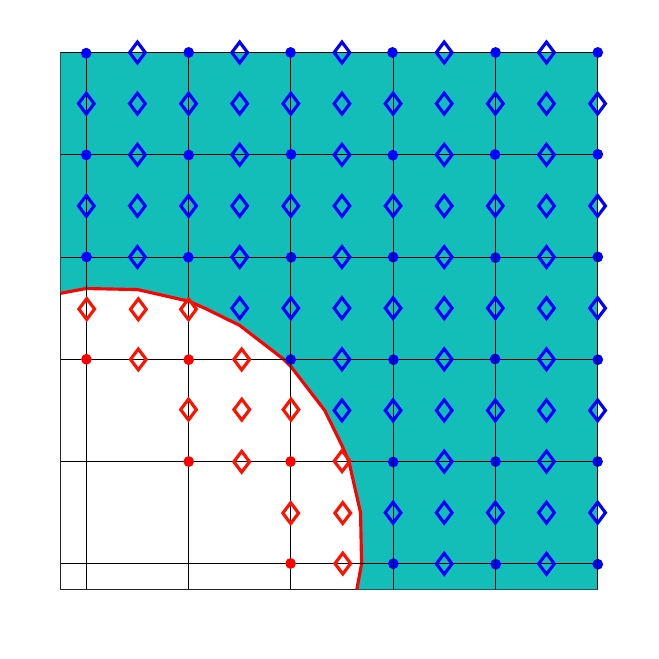}
\put(39.5,39.5){\Large P}
  \linethickness{1pt}
  \put(19,26){\framebox(42,40){}}
\end{overpic}
    \caption{\textit{Before snapping}}
    \label{fig:points_before}
    \end{subfigure}
    \hfill 
    \begin{subfigure}{0.45\textwidth}
\begin{overpic}[width=\textwidth]{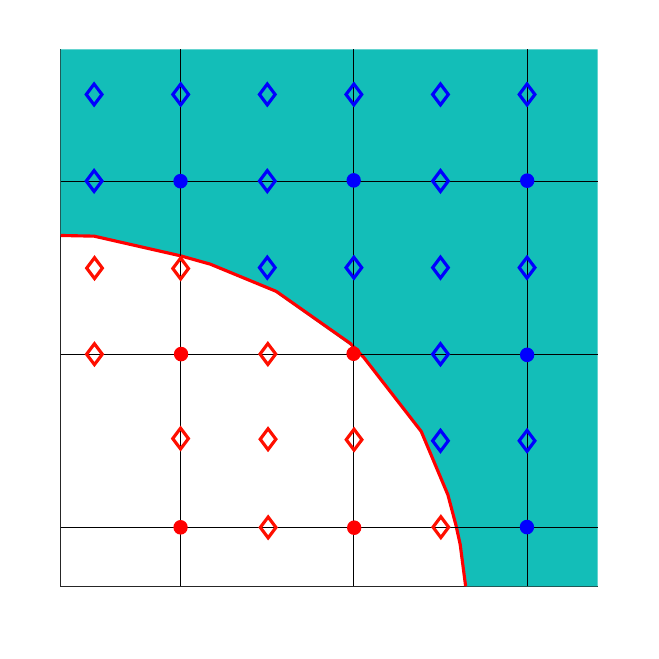}
\put(47,40){\Large G}
\end{overpic}
\caption{\textit{After snapping (zoom-in)}}
\label{fig:points_after}
    \end{subfigure}
    \caption{\textit{Classification of points and results of the snapping-back-to-grid algorithm. The internal point set $\mathcal I_v$ is represented by blue dots and diamonds, while the set of $\mathcal I_p$ is denoted by blue dots only. Red markers denote the respective sets of ghost points $\mathcal G_v$ and $\mathcal G_p$. Panel (a) shows the grid before snapping, while panel (b) shows the grid after snapping. Panel (b) shows a zoom-in of the box highlighted in panel (a), where a point originally classified as internal, $P$, has changed its classification to ghost, G.}}
    \label{fig:points}
\end{figure}
The discrete spaces ${\bf V}_h$ and $Q_h$ are constructed from continuous, piecewise polynomial functions defined on $\Omega_\text{act}$.
The functions ${\bf w}$ and $q$, which are defined in the function spaces $\bf V$ and $Q$, respectively, are approximated by the functions ${\bf w}_h$ and $q_h$ that belong to finite-dimensional subspaces ${\bf V}_h$ and $Q_h$, defined as follows:

\begin{subequations}
\label{eq:basis}
\begin{align}\label{eq:Vh}
{\bf V}_h & = \{ {\bf v}_h \in {\bf V} : {\bf v}_h|_K \in [\mathbb Q_2(K)]^2, \quad \forall K\in \Omega_\text{act}\}, \\ \label{eq:Qh}
Q_h & = \{ q_h \in Q : {q_h}|_K \in \mathbb Q_1(K), \quad \forall K\in \Omega_\text{act}\},
\end{align}
\end{subequations}
where $\mathbb Q_1(K)$ and $\mathbb Q_2(K)$ denote, respectively, the spaces of bilinear and biquadratic Lagrangian shape functions on the element $K$.
This choice of finite element spaces for the velocity/pressure fields is inf-sup stable~\cite{ern2021finite}, leading to third order accuracy for the velocity field and second order for the pressure.

In each cut cell, the portion of boundary is approximated by a straight line (see Algorithm~\ref{alg:ab} for more details). The union of these segments is denoted by $\Gamma_h$ and the corresponding approximated domain by $\Omega_h$.
Consequently, the integrals in \eqref{pro:variational} defined over $\Omega$ and its boundary $\Gamma$ are now evaluated over $\Omega_h$ and $\Gamma_{h}$, respectively, leading to the following discrete problem.

\begin{pro}\label{pro:variational_2}
Given ${\bf f}_h, {\bf g}_h \in L^2(\Omega)$, find ${\bf u}_h \in {\bf V}_h$ and $p_h \in Q_h$ such that, for almost every $t\in(0,T)$, it holds
\begin{subequations}
\label{eq:problem_discr}
\begin{align}
\label{eq:problem_discr_a}
\left( \pad {{\bf u}_h}{t},{\bf w}_h \right)_{L^2(\Omega_h)} &+ a_h({\bf u}_h,{\bf w}_h) + b_h({ p}_h,{\bf w}_h)+ c_h({\bf u}_h,{\bf u}_h,{\bf w}_h)   = \left( {\bf f}_h, {\bf w}_h\right)_{L^2(\Omega_h)} + \\ \notag &+ \lambda ({\bf g}_h,{\bf w}_h)_{L^2(\Gamma_{h})} -  \left({\bf g}_h , \nabla  {\bf w}_h\cdot {\bf n} \right)_{L^2(\Gamma_{h})},  \qquad  \forall {\bf w}_h \in {\bf V}_h\\ 
 b_h(q_h,{\bf u}_h) & = (q_h{\bf n}, {\bf g}_h )_{L^2(\Gamma_h)},  \qquad  \forall q_h \in Q_h
\end{align}
\end{subequations}
where 
\begin{subequations}
    \begin{align}
    a_h({\bf u}_h,{\bf w}_h) & = \left( \nabla {\bf u}_h,  \nabla {\bf w}_h \right)_{L^2(\Omega_h)} - \left(\nabla {\bf u}_h \cdot {\bf n} , {\bf w}_h \right)_{L^2(\Gamma_{h})} -  \left({\bf u}_h , \nabla  {\bf w}_h\cdot {\bf n} \right)_{L^2(\Gamma_{h})} + \lambda ({\bf u}_h,{\bf w}_h)_{L^2(\Gamma_{h})}  \\
    b_h({ p}_h,{\bf w}_h) & =  -\left(  p_h, \nabla \cdot{\bf w}_h \right)_{L^2(\Omega_h)} + \left(p_h {\bf n} , {\bf w}_h \right)_{L^2(\Gamma_{h})} \\
        c_h({\bf v}_h,{\bf u}_h,{\bf w}_h) & = \left( \left({{\bf v}_h} \cdot \nabla \right){{\bf u}_h} , {\bf w}_h \right)_{L^2(\Omega_h)}.
    \end{align}
\end{subequations}
\end{pro}

We define bases for the spaces ${\bf V}_h$ and $Q_h$. In particular, we first define two grids of nodes: $\mathcal N_v$ and $\mathcal N_p$, as follows.
The grid $\mathcal N_v$ is constituted by the points $(x_\text{min} + ih/2,y_\text{min} + jh/2)$, for $i=0,\ldots,2N_x$ and $j=0,\ldots,2N_y$.
The grid $\mathcal N_p$ is constituted by the points $(x_\text{min} + ih,y_\text{min} + jh)$, for $i=0,\ldots,N_x$ and $j=0,\ldots,N_y$. Note that $\mathcal N_p \subset \mathcal N_v$ and that the closure of each cell $K$, denoted as $\bar{K}$, contains nine nodes of $\mathcal N_v$ and four nodes of $\mathcal N_p$.
For both grids, we define the respective sets of active nodes $\mathcal A_v$ and $\mathcal A_p$ by $\mathcal A_{[v,p]} = \mathcal N_{[v,p]} \cap \Omega_\text{act}$.
We partition active node sets into two subsets, respectively:
$\mathcal A_{[v,p]} = \mathcal I_{[v,p]} \cup \mathcal G_{[v,p]}$, where $\mathcal I_{[v,p]}$ are the sets of internal nodes, and $\mathcal G_{[v,p]}$ are the sets of ghost nodes. We define them as follows: $\mathcal I_{[v,p]} = \mathcal A_{[v,p]} \cap \Omega$, $\mathcal G_{[v,p]}=\mathcal A_{[v,p]} \backslash \mathcal I_{[v,p]}$.
In Fig.~\ref{fig:points_before}, internal nodes of $\mathcal I_v$ are depicted by blue markers, while the subset of active nodes of $\mathcal I_p$ are the vertices of the cells and denoted by blue dots. Ghost points are denoted by red markers.

Let $D_{[v,p]} = \left| \mathcal A_{[v,p]}\right|$ be the cardinality of the respective active node sets.
We observe that $\dim{\bf V}_h = 2D_{v}$ and $\dim Q_h = D_{p}$.  
We associate to each active node $\vec{x}_i \in \mathcal A_{[v,p]}$ a global shape function $\vphi^{[v,p]}_i$ uniquely defined by
\[
\varphi^{[v,p]}_i \in \mathbb Q_{[2,1]},
\qquad
\text{supp} \; \varphi^{[v,p]}_i = \bigcup_{\bar{K} \ni \mathbf x_i} K,
\qquad
\varphi^{[v,p]}_i(\mathbf x_j)=\delta_{ij},
\qquad
\forall \mathbf x_j\in\mathcal A_{[v,p]}.
\]
The basis for ${\bf V}_h$ is defined through shape functions $\left\{ \vphi^{v}_1 {\bf e}_1, \ldots, \vphi^{v}_{D_v} {\bf e}_1, \vphi^{v}_1 {\bf e}_2, \ldots, \vphi^{v}_{D_v} {\bf e}_2 \right\}$, where ${\bf e}_1$ and ${\bf e}_2$ are the vectors the canonical basis. Analogously, we define a basis for $Q_h$ as $  \left\{ {\vphi}_i^p, i=1,\ldots,D_{p} \right\}$.

The discrete spaces are therefore
\[
{\bf V}_h =
\left\{
{\bf v}_h = {v}^1_h {\bf e}_1+{v}^2_h {\bf e}_2, \quad 
{v}^{[1,2]}_h = \sum_{i=1}^{D^v} {v}^{[1,2]}_i\, \varphi^v_i
\right\},
\qquad
Q_h =
\left\{
q_h =
\sum_{i=1}^{D^p} q_i\, \varphi^p_i
\right\}.
\]

We reformulate system \eqref{eq:problem_discr} in a compact form, which will be used in the time discretization. To this end,
we reuse (when it is clear from the context) the notation ${\bf u}_h$ and $p_h$ to denote the vectors of nodal values, that correspond to the components with respect to the bases introduced above. Introducing computational matrices representing the spatial derivatives, we obtain the following expression
\begin{subequations}
\label{eq:FEM_semidiscr}
\begin{align}
\mathbb M\pad{ {\bf u}_h}{ t} & = \mathbb H[{\bf u}_h] \, {\bf u}_h + \mathbb G p_h + \mathbb M \left({\bf f}_h + \lambda {\bf g}_h\right) + {\bf \tilde{g}}_h\\ 
\mathbb G^\top\, {\bf u}_h & = {\hat g}_h,
\end{align}
\end{subequations}
where ${\mathbb M} = 
\begin{pmatrix}
M & 0 \\
0 & M
\end{pmatrix}
$, with
$M_{ij}=(\varphi^v_i,\varphi^v_j)_{L^2(\Omega_h)}$ (the mass matrix),
$\mathbb H= 
\begin{pmatrix}
H^1[{\bf u}_h] & 0 \\
0 & H^2[{\bf u}_h]
\end{pmatrix}$,
with
$(H^{[1,2]}[{\bf u}_h])_{ij}=a_h(\vphi^v_i{\bf e}_{[1,2]},\vphi^v_j {\bf e}_{[1,2]})+c_h({\bf u}_h,\vphi^v_i{\bf e}_{[1,2]},\vphi^v_j {\bf e}_{[1,2]})$,
$\mathbb G = 
\begin{pmatrix}
G^1 \\
G^2
\end{pmatrix}$,
with
$G^{[1,2]}_{ij} = b_h(\vphi^p_j,\vphi^v_i {\bf e}_{[1,2]})$, 
${\bf \tilde{g}}_h
= \begin{pmatrix}
\tilde{g}^1 \\
\tilde{g}^2
\end{pmatrix}$,
with 
$\tilde{g}_i^{[1,2]} = \left({\bf g}_h , \nabla  \vphi^v_i {\bf e}_{[1,2]} \cdot {\bf n} \right)_{L^2(\Gamma_{h})}$,
$(\hat{g}_h)_{i} = (\vphi^p_i {\bf n}, {\bf g}_h )_{L^2(\Gamma_h)}$.
We remind that these operators also include the boundary terms.

Eq.~\eqref{eq:FEM_semidiscr} can be summarized as follows
\begin{equation}\label{eq:Theta_eps}
	\BB \td{ \textbf{q}_h }{t} = \Theta [ \textbf{q}_h ]\textbf{q}_h + {\bf \bar{g}}_h,
\end{equation}
where $\textbf{q}_h = ({ \bf u}_h, p_h)^\top$, and:
\begin{equation}
\label{eq:sys_mon}
 \BB  = \begin{pmatrix}
 \mathbb M & 0 \, \\
  0 & 0 \\
\end{pmatrix},\qquad 
  \Theta[{\bf q}_h] = \begin{pmatrix}
\mathbb H[{\bf u}_h] & \mathbb G  \\
\mathbb G^\top & 0 
\end{pmatrix}, 
\qquad
{\bf \bar{g}}_h=
\BB \left({\bf f}_h + \lambda {\bf g}_h\right) +
\begin{pmatrix}
{\bf \tilde{g}}_h \\
{\hat{g}}_h
\end{pmatrix}.
\end{equation}


\begin{algorithm}
\caption{Computation of the intersection of the boundary $\Gamma$ with the grid (see Fig.~\ref{fig:ab})}\label{alg:ab}
\begin{algorithmic}
\State $k_4 = k_0$
\For{i = 0:3}
\If{$\phi(k_i)\phi(k_{i+1})<0$} 
   \State $\theta = \phi(k_i)/(\phi(k_i)-\phi(k_{i+1})) $ 
   \State $P = \theta k_{i+1}+(1-\theta)k_i$
   \If{$\phi(k_i)<0 $}
       \State ${\bf A}:=P$
       \Else
       \State ${\bf B}:=P$
   \EndIf
\EndIf 
\EndFor
\end{algorithmic}
\end{algorithm}


\subsubsection{Geometrical error}
In ghost-FEM, boundary conditions are not applied on the true boundary $\Gamma$ but on an approximate one $\Gamma_h$, and it is essential to preserve the accuracy of their evaluation. To prevent any loss of accuracy that could degrade the overall convergence rates, we correct the boundary conditions via Taylor expansions, following the strategy proposed in the Shifted Boundary Method \cite{Atallah2022, song2018shifted}. We briefly describe the mentioned strategy. 

We first select a mapping between the real boundary $\Gamma$ and the approximate one $\Gamma_h$, as follows
\[ \mathcal M_h : \Gamma_h \to \Gamma\]
which maps any point ${\bf x} \in\Gamma_h$ to its orthogonal projection ${\bf z} = \mathcal{M}_h({\bf x}) $ on the physical boundary $\Gamma$. The projection point is computed via the Algorithm~\ref{alg:closest_point} (see  \cite[Sec.~3.1.2]{Coco2024} for more details). 

The mapping $\mathcal M_h$ can be characterized through a distance vector function ${\bf d}$ defined by
\[ {\bf d}({\bf x}) := {\bf z}-{\bf x} = (\mathcal M_h - I)({\bf x}). \]
Let us assume that $u$ is sufficiently smooth in the strip between $\Gamma_h$ and $\Gamma$ so that it admits a $m$-th order Taylor expansion, and let us denote by $\mathcal D^i_{\bf d} u = \sum_{\alpha \in \mathbb N^2, \, |\alpha| = i} \frac{i!}{\alpha!} \frac{\partial^i u}{\partial z^\alpha}{\bf d}^\alpha$ the $i$-th order directional derivative along ${\bf d}$. Then, we can write
\[ u({\bf z}) = u( {\bf x} + {\bf d}({\bf x})) = u({\bf x}) + \sum_{i=1}^m \frac{\mathcal D^i_{\bf d} u({\bf x})}{i!} + o(||{\bf d}||^m).\]
The Dirichlet boundary condition $u({\bf z})=g({\bf z})$ on $\Gamma$ is replaced by a boundary condition on ${\Gamma}_h$ via
\[
u({\bf{x}}) = g(M_h({\bf x})), \qquad {\bf x} \in \Gamma_h.
\]
In this work, our aim is third-order accuracy, and therefore it suffices to set $m=1$.

\begin{algorithm}
\caption{Computation of the projection point ${\bf z} = \mathcal{M}_h({\bf x}) $}\label{alg:closest_point}
\begin{algorithmic}
\State toll = $10^{-2} h^2$
\State $\epsilon=h/10$
\State $k = 0$
\State $\Delta = 2 \; \text{toll}$
\State ${{\bf x}}^{(0)} = {\bf x}$
\While{ $\Delta$ > toll}
   \State ${\bf a}^{(k)} = -\epsilon \phi({{\bf x}}^{(k)}) \frac{\nabla \phi({{\bf x}}^{(k)}) }{|\nabla \phi({{\bf x}}^{(k)})|^2} $
   \State ${{\bf x}}^{(k+1)} = {{\bf x}}^{(k)} + {\bf a}^{(k)}$ 
   \State $\Delta = ||{{\bf x}}^{(k+1)} - {{\bf x}}^{(k)} ||$
   \State $k = k+1$
\EndWhile
\State ${\bf z}={{\bf x}}^{(k)}$
\end{algorithmic}
\end{algorithm}


\subsubsection{Choice of the penalization parameter and eigenvalue problem}
We devote this section to discussing the choice of the penalization parameter $\lambda$ in Eq.~\eqref{eq:problem_discr_a}. 
In \cite{astuto2024nodal}, it was shown that setting $\lambda = c h^\alpha$ with $c \in \mathbb R$ sufficiently large, and $\alpha = 2$ 
ensures well-posedness for the Poisson's problem. In the present work, we extend this analysis to time-dependent 
problems and higher-order discretizations. 
Let $c_0 >0$ be a constant such that the bilinear form is coercive. Then, for any $c>c_0$, the formulation is well-posed. However, combined with IMEX, the choice of this constant is crucial: in some cases, when the penalization parameter is too small, it can lead to a loss of stability in the discretization,  while large values result in severe ill-conditioning of the linear system. The selection of the stabilization parameter is therefore critical and several approaches for estimating this parameter are based on solving general eigenvalue problems \cite{griebel2003particle,dolbow2009efficient,garhuom2022eigenvalue}. 

The bilinear form in Eq.~\eqref{eq:problem_discr} is coercive if there exists a positive constant $\tilde{C}>0$ such that  $a_h({\bf u}_h,{\bf u}_h) \geq \tilde{C}\|{\bf u}_h\|_{{\bf V}_h}$ for any ${\bf u}_h \in {\bf V}_h$. Here, the norm on the FE space ${\bf V}_h$ is defined as 
\begin{equation}\label{eq:norm_vh}
    \|{\bf u}_h\|_{{\bf V}_h}^2 = \| \nabla {\bf u}_h \|^2_{L^2(\Omega_h)} + \lambda \|{\bf u}_h\|^2_{L^2(\Gamma_h)}.
\end{equation} 

To establish a lower bound on the bilinear form $a_h({\bf u}_h, {\bf u}_h)$, we use the strategy of~\cite{Saberi2023}, as described below.
We apply Young’s inequality, which states that for any two real numbers $a$ and $b$ and a positive parameter $\varepsilon$ we have
$ab \leq a^2/(2\varepsilon) + \varepsilon b^2/2$.
Applying it to the boundary integral term, we obtain: 
\[
\int_{\Gamma_h} 2 {\bf u}_h({\bf n}\cdot\nabla {\bf u}_h) {\rm d} s \leq 
\frac{1}{\varepsilon}\int_{\Gamma_h} {\bf u}_h {\bf u}_h {\rm d} s
+ \varepsilon \int_{\Gamma_h} ({\bf n}\cdot\nabla {\bf u}_h) ({\bf n}\cdot\nabla {\bf u}_h) {\rm d} s.
\]
Therefore:
\begin{align*}
    a_h({\bf u}_h,{\bf u}_h) &= \int_{\Omega_h}  \nabla  {\bf u}_h\cdot \nabla  {\bf u}_h {\rm d} \Omega + \int_{\Gamma_h} \lambda {\bf u}_h {\bf u}_h - 2 {\bf u}_h({\bf n}\cdot\nabla {\bf u}_h) {\rm d} s \\ 
    &\ge \int_{\Omega_h}  \nabla  {\bf u}_h\cdot \nabla  {\bf u}_h {\rm d} \Omega  -  \varepsilon \int_{\Gamma_h}({\bf n}\cdot  \nabla  {\bf u}_h) ({\bf n}\cdot  \nabla  {\bf u}_h)  {\rm d} s + \left(\lambda -\frac{1}{\varepsilon}\right)\int_{\Gamma}  {\bf u}_h {\bf u}_h  {\rm d} s\\
    &\ge \left( 1- C\varepsilon \right)\int_{\Omega_h}  \nabla  {\bf u}_h\cdot \nabla  {\bf u}_h \, {\rm d} \Omega + \left(\lambda -\frac{1}{\varepsilon}\right)\int_{\Gamma_h}  {\bf u}_h {\bf u}_h  {\rm d} s, 
\end{align*}
where in the last expression we assume that there exists a constant $C>0$ such that 
\begin{equation}\label{eq:lowC}
C \int_{\Omega_h}  \nabla  {\bf u}_h\cdot \nabla  {\bf u}_h {\rm d} \Omega \ge \int_{\Gamma_h}({\bf n}\cdot  \nabla  {\bf u}_h) ({\bf n}\cdot  \nabla  {\bf u}_h)  {\rm d} s.
\end{equation}
The bilinear form is coercive if $1-C\varepsilon > 0 $ and $\lambda - \frac{1}{\varepsilon} > 0$. Hence, a necessary condition for the stabilization parameter is $\lambda > C$.

\par The first step is to find an estimate for the smallest value of $C$ such that \eqref{eq:lowC} is satisfied.
This is obtained by solving the generalized global eigenvalue problem 
\begin{equation}\label{eq:geneigprob}
{\bf K} v = \Lambda {\bf M} v,
\end{equation}
where ${\bf K}_{ij} = \int_{\Gamma_h} ({\bf n}\cdot  \nabla  \vphi_i^v)({\bf n}\cdot  \nabla  \vphi_j^v)~\mathrm{d}s$ and ${\bf M}_{ij} = \int_{\Omega_h}   \nabla  \vphi^v_i \cdot   \nabla  \vphi^v_j~\mathrm{d}\Omega$ for $i,j \in \{1,2,\dots,D^v\}$.
The lower bound of $C$ is then estimated by taking $C = \max \Lambda$. The stabilization parameter is then chosen as $\lambda = \gamma C$, for some constant $\gamma > 1$ (we used $\gamma = 1.1$ in our tests).



\subsubsection{Stabilization technique for the saddle-point problem}
\label{section:saddle}
The discretization of the pressure term may introduce instabilities, as the pressure appears only in the form of its gradient, resulting in a non-unique solution of the problem~\ref{eq:sys_mon}. In finite difference schemes, this issue can be addressed by employing a mixed interpolation strategy based on the MAC grid \cite{LANGTANGEN20021125,COCO2020109623}. Analogously, when using finite element methods, a mixed formulation is adopted \cite{boffi2013mixed}. 

In Section~\ref{sec:IMEX}, we will 
employ a semi-implicit time discretization, then a linear system will be solved at each time step. The matrix of the linear system will take the form 
\[ \BB - k\Delta t \, \Theta[{\bf q}_h],\]
for some $ k\in \mathbb R$ that depends on the choice of the IMEX scheme (a detailed explanation is given in Section~\ref{sec:IMEX}.)
The singularity of the matrix $\BB - k\Delta t \, \Theta[{\bf q}_h]$ is addressed by introducing a stabilization diagonal submatrix to $\Theta[{\bf q}_h]$ as follows
\begin{equation} \Theta[{\bf q}_h] = \begin{pmatrix}
\mathbb H[{\bf u}_h] & \mathbb G  \\
\mathbb G^\top & \varepsilon \mathbb I
\end{pmatrix}. 
\end{equation} 
We set $\varepsilon = 10^{-10}$ in our numerical simulations.

\subsection{Temporal discretization}
\label{sec:IMEX}
In this section we describe the time discretization strategy to solve Eq.~\eqref{eq:Theta_eps}. We adopt an IMEX scheme, which offers a flexible framework to obtain arbitrary order of accuracy.

We consider a final time $T$ and define the time step size as $\Delta t = T/N_{\rm ts},\, N_{\rm ts} \in \mathbb N,$ denoting the time steps by $t^n = n\Delta t$ and $\textbf{q}^n \approx \textbf{q}(t^n), \, n = 0,\cdots,N_{\rm ts}$.

Observe that vector $\bar{\bf g}_h$ of Eq.~\eqref{eq:Theta_eps} may depend on time (for example, when source terms and boundary conditions are time dependent), then we write $\bar{\bf g}_h(t)$.

We first set $\textbf{q}^1_E = \textbf{q}^n$, then the stage fluxes $\textbf{q}_E^{i}$ and $\textbf{q}_I^{i}$ are calculated as
\begin{subequations}  
\label{eq_imex1}
\begin{align}
\label{eq_imex_qE}
   \BB \, \textbf{q}_E^{i}& = \BB \, \textbf{q}^n + \Delta t\,\sum_{j=1}^{\rm i-1}\tilde a_{i,j}\left(\Theta[\textbf{q}_E^j]\,\textbf{q}_I^j + {\bf \bar{g}}_h(t^n+\tilde{c}_j \Delta t) \right), \\
\label{eq_imex_qI}
 \BB \, \textbf{q}_I^{i} &= \BB \, \textbf{q}^n + \Delta t\,\sum_{j=1}^{\rm i} a_{i,j} \left(\Theta[\textbf{q}_E^j]\,\textbf{q}_I^j + {\bf \bar{g}}_h(t^n+c_j \Delta t) \right),
\end{align}
\end{subequations}
for $i=1,\ldots,s$, where $\rm s$ is the number of stages of the scheme.
The numerical solution is finally updated with
\begin{align}
\label{eq_imex_q_bi}
 \BB \,  \textbf{q}^{n+1} &= \BB \,\textbf{q}^n + \Delta t\sum_{i=1}^{\rm s} \left( \tilde{b}_i \left( \Theta[\textbf{q}_E^i]\,\textbf{q}_I^i+
 {\bf \bar{g}}_h(t^n+\tilde{c}_i \Delta t)
 \right) +
 b_i \left( \Theta[\textbf{q}_E^i]\,\textbf{q}_I^i+
 {\bf \bar{g}}_h(t^n+c_i \Delta t)
 \right) \right).
\end{align}
The coefficients $\tilde a_{i,j}$, $a_{i,j}$, $\tilde c_{i}$, $c_{i}$, $\tilde{b}_i$, $b_i$ are taken from double Butcher tableaux
\[
\begin{array}
{c|c}
\tilde c & \widetilde{A}\\
\hline
& \tilde b 
\end{array} \qquad 
\begin{array}
{c|c}
 c & A\\
\hline
& b
\end{array}
\]
In our numerical tests we choose ${\rm s} = 2$ or ${\rm s} = 4$ for a second or third order accuracy, respectively, and adopt the following Butcher tableaux~\cite{boscarino2024implicit}.

\subsubsection*{Two-stage second-order IMEX scheme}
For a second-order numerical scheme (${\rm s} = 2$):
\[
\begin{array}
{c|cc}
0 & 0 & 0\\
1 & 1 &0 \\
\hline
& 1/2 & 1/2 
\end{array}
\qquad 
\begin{array}
{c|cc}
1-\eta & 1-\eta & 0\\
\eta & \eta - \delta & \delta \\
\hline
& 1/2 & 1/2
\end{array}
\]
with $\eta = 1/\sqrt{2}$ and $\delta = 1- 1/(2\eta)$.

\subsubsection*{Four-stage third-order IMEX scheme}
For a third-order numerical scheme (${\rm s} = 4$):
\[
\begin{array}{c|cccc}
     0&  0 & 0 &0 & 0 \\
	\gamma&  \gamma &0& 0 & 0 \\
	0.717933260754 & 1.243893189 &-0.5259599287 &0 & 0 \\
	1& 0.6304125582 &0.7865807402 &-0.4169932983& 0 \\
    \hline
    & 0& 1.208496649176& -0.644363170684& \gamma
\end{array}
\]
\[
\begin{array}{c|cccc}
   \gamma &\gamma & 0 & 0 & 0\\
	\gamma & 0 &\gamma &0 &0 \\
	 0.717933260754& 0 &0.282066739245 &\gamma & 0 \\
	  1 & 0 & 1.208496649176 &-0.644363170684 & \gamma \\
      \hline
      & 0& 1.208496649176& -0.644363170684& \gamma
\end{array}
\]
with $\gamma = 0.435866521508$.

\subsubsection{Extrapolation technique and moving domains}
When dealing with moving domains, 
we need a strategy to extrapolate
functions beyond the current domain.
In order to solve Eqs.~\eqref{eq_imex1} and~\eqref{eq_imex_q_bi} on $\Omega^{n+1}_\text{act}$, we need to extend ${\bf q}^n$, originally defined on $\Omega^{n}_\text{act}$, to the active nodes of $\Omega^{n+1}_\text{act}$ not contained in $\Omega^{n}_\text{act}$.
We remark that the motion of the domain is prescribed {\it a priori} (one-way coupling).

Here, we adapt the strategy proposed in \cite{aslam2004partial} for finite difference schemes to the finite element framework. In Algorithm~\ref{alg_unn} we show how to obtain a quadratic extrapolation (to achieve third order accuracy) of a function $u$. 

\begin{figure}
  \centering
    \begin{subfigure}{0.32\textwidth}
        \includegraphics[width=\textwidth]{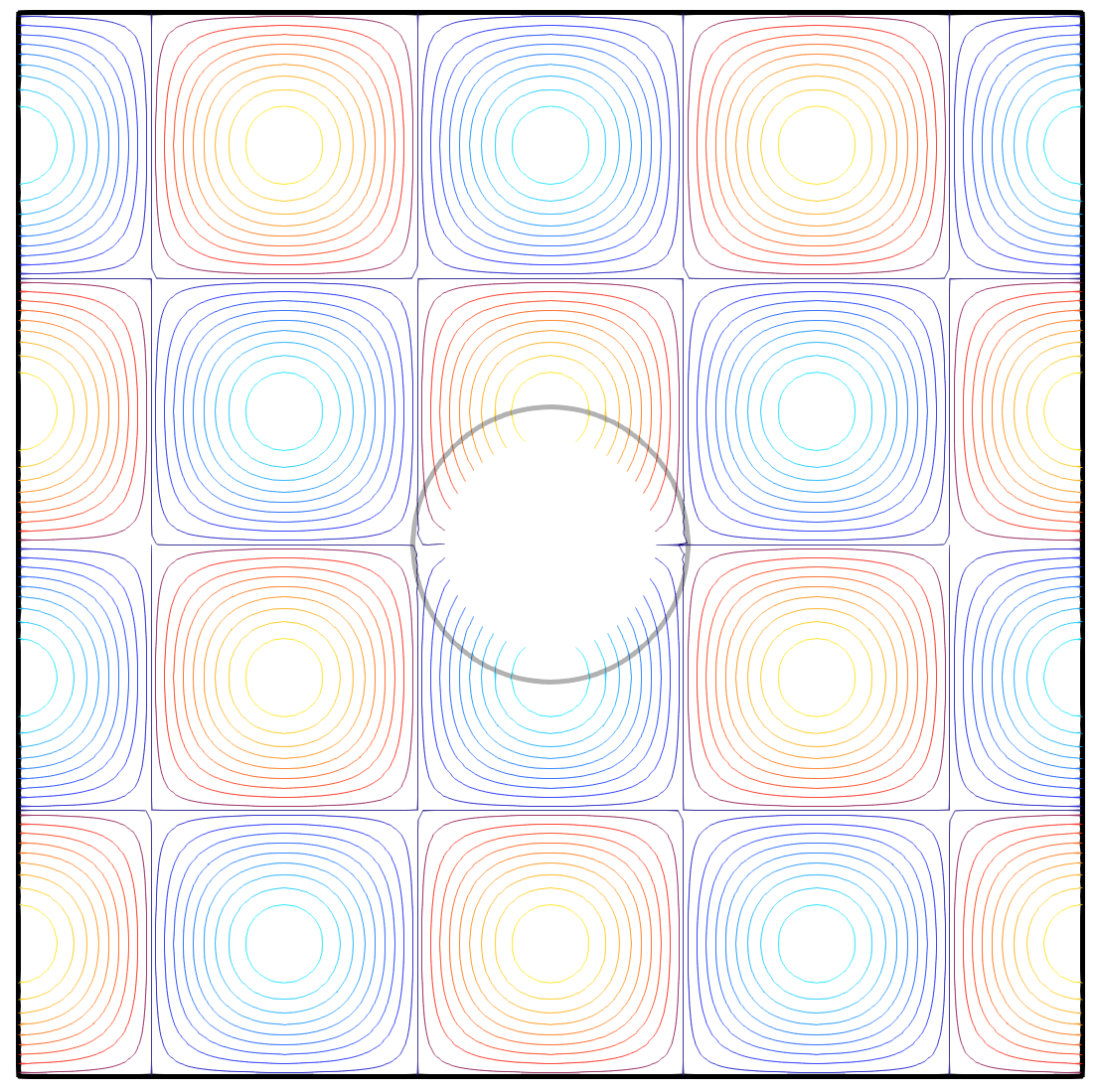}
        \caption{\textit{Disk}}
    \end{subfigure}
\hfill
\begin{subfigure}{0.32\textwidth}
        \includegraphics[width=\textwidth]{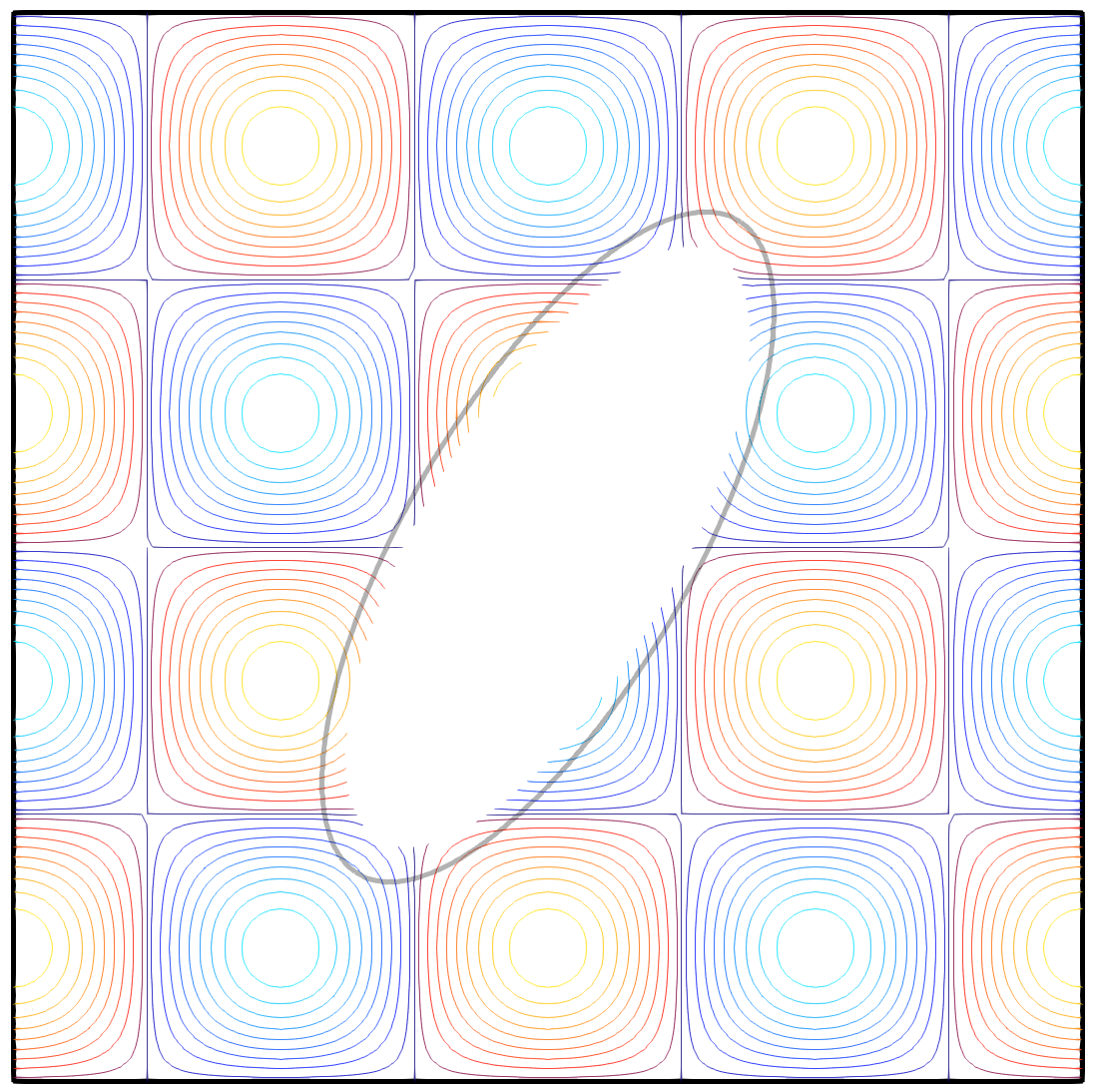}
        \caption{\textit{Ellipse}}
    \end{subfigure}
    \hfill
\begin{subfigure}{0.32\textwidth}
    \includegraphics[width=\textwidth]{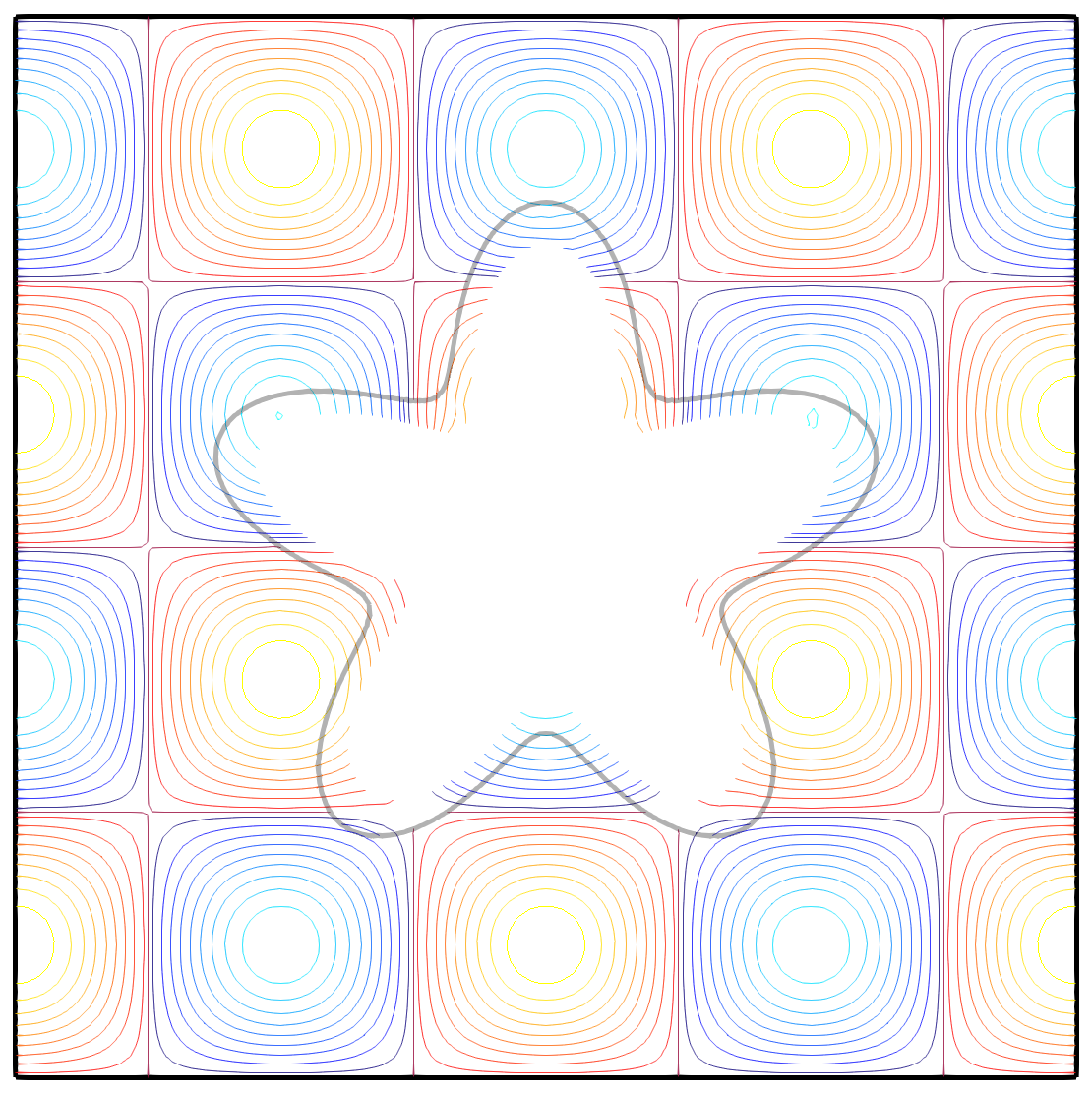}
    \caption{\textit{Flower}}
\end{subfigure}
    \centering
    \caption{\textit{Plot of extrapolated solutions using quadratic extrapolation on a narrow band near the boundary for different geometries.}}
    \label{fig:aslam}
\end{figure}

\begin{figure}
    \centering
    \begin{subfigure}{0.32\textwidth}
        \includegraphics[width=\textwidth]{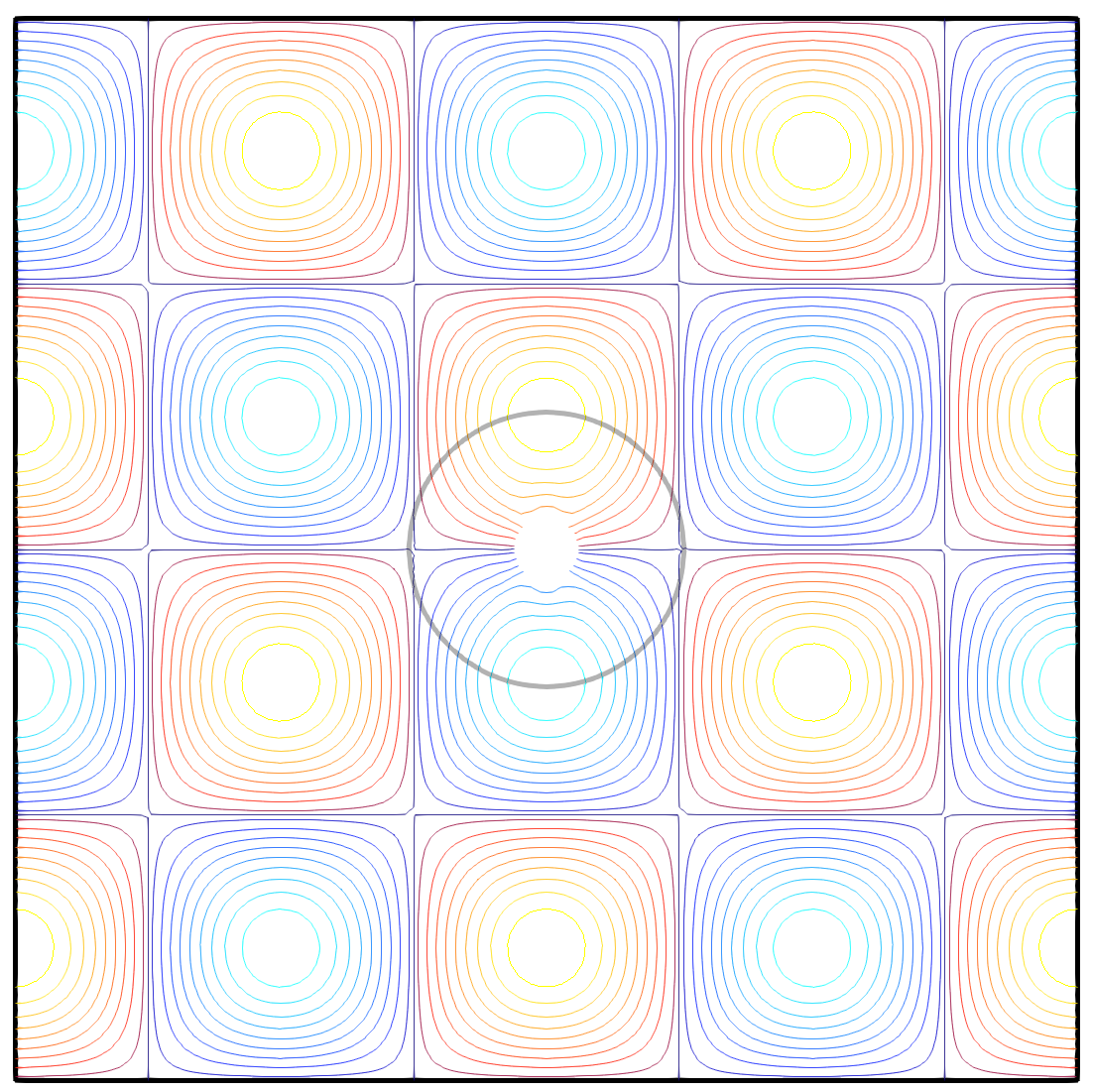}
        \caption{\textit{Disk}}
    \end{subfigure}
\hfill
\begin{subfigure}{0.32\textwidth}
        \includegraphics[width=\textwidth]{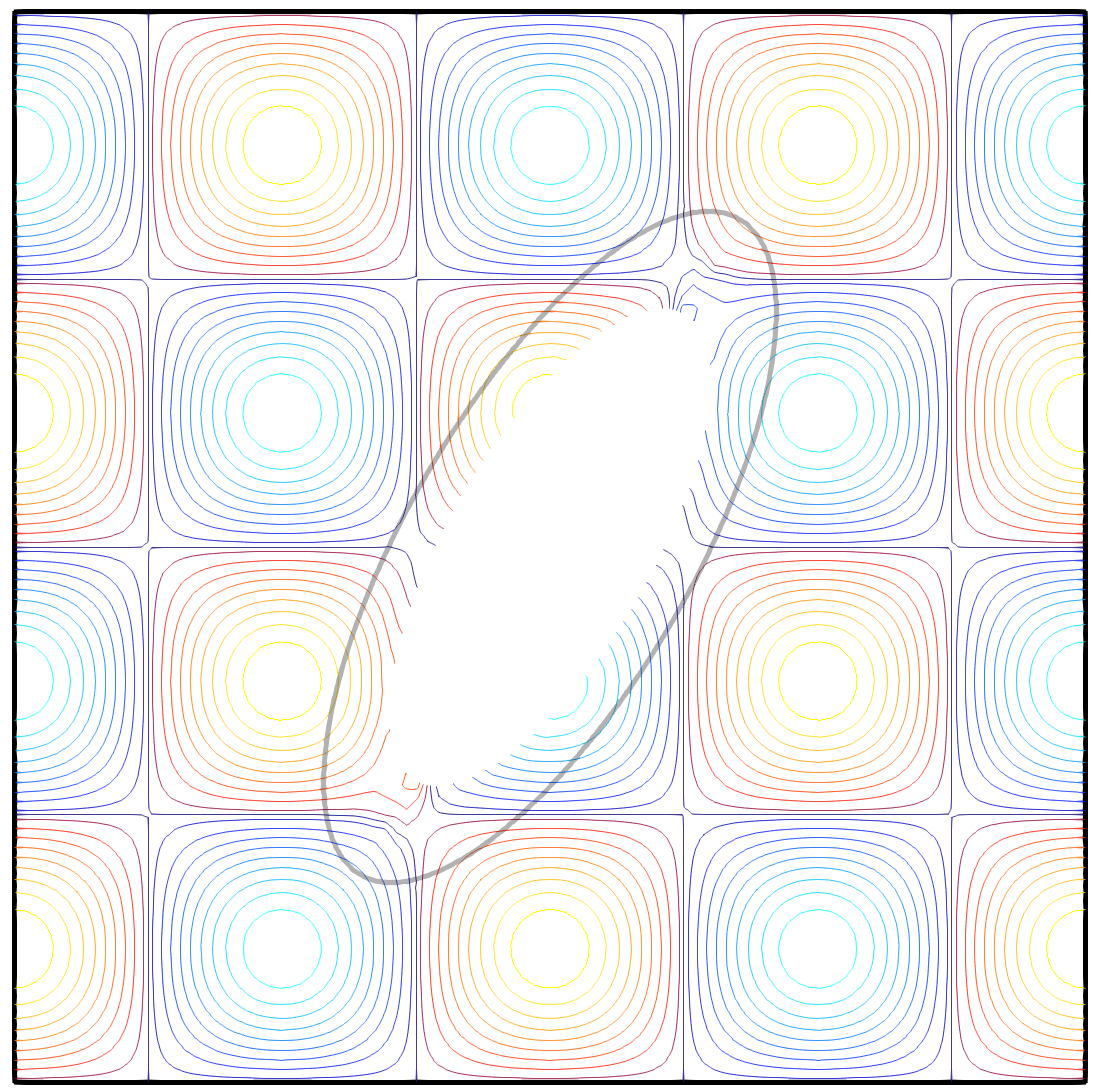}
        \caption{\textit{Ellipse}}
    \end{subfigure}
    \hfill
\begin{subfigure}{0.32\textwidth}
    \includegraphics[width=\textwidth]{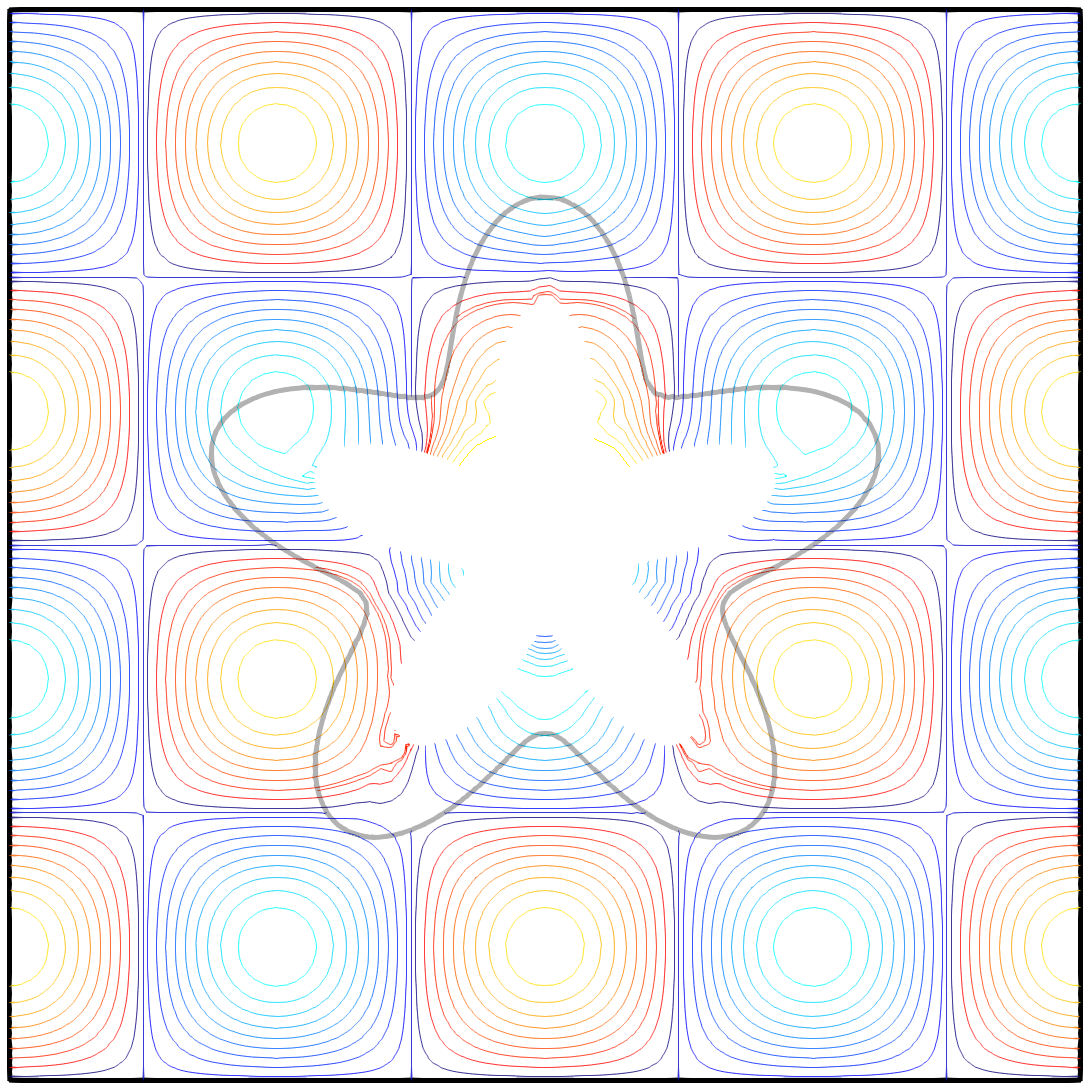}
    \caption{\textit{Flower}}
\end{subfigure}

    \caption{\textit{Plot of extrapolated solutions using quadratic extrapolation on a wider band near the boundary for different geometries, using upwind discretization for the convective term.}}
    \label{fig:aslam_xl}
\end{figure}

We start defining the characteristic function of 
$\Omega^n_\text{com} \coloneqq R \setminus\Omega^n_\text{act}$, as
\begin{equation}
    \chi(\bf x) = \begin{cases}
    0 & \text{if } {\bf x} \in \Omega^n_\text{act}, \\
    1 &  \text{ otherwise}.
  \end{cases}
\end{equation}
Making use of it, we preserve the values of the solution at the nodes of $\Omega^n_\text{act}$ (active nodes), while evolving the function in $ \Omega^n_\text{com}$, according to the Algorithm~\ref{alg_unn}.
\begin{algorithm}[H]
\caption{ Quadratic extrapolation of a function $u$ to $ \Omega^n_\text{com}$. Eqs~\eqref{eq:extr_unn}, \eqref{eq:extr_un} and \eqref{eq:extr_u} are solved
in a fictitious time $\tau$ up to steady-state in a narrow boundary layer.}\label{alg_unn}
\begin{algorithmic}
\State Compute $u_{\bf n} :=  {\bf n} \cdot  \nabla u $
\State Compute $u_{\bf nn} :=  {\bf n} \cdot  \nabla u_{\bf n} $
\State Solve 
\begin{equation}
\label{eq:extr_unn}
    \pad{ {u_{\bf nn}}}{ \tau} + \chi \, {\bf n} \cdot \nabla u_{\bf nn} = 0
\end{equation}
\State Solve 
\begin{equation}
\label{eq:extr_un}
    \pad{ {u_{\bf n}}}{ \tau} + \chi \left( {\bf n} \cdot \nabla u_{\bf n} - u_{\bf nn} \right) = 0
\end{equation}
\State Solve 
\begin{equation}
\label{eq:extr_u}
    \pad{ u}{ \tau}  + \chi \left(  {\bf n} \cdot \nabla u - u_{\bf n} \right)  = 0
\end{equation}
\end{algorithmic}
\end{algorithm}
The strategy of Algorithm~\ref{alg_unn} is as follows: the Heaviside function ensures that the known values of our function remain unchanged in $ \Omega^n_\text{act}$. In the region $ \Omega^n_\text{com}$, Eq.~\eqref{eq:extr_unn} defines a linear hyperbolic PDE with characteristics along ${\bf n}$; we solve it until we reach steady-state in a narrow boundary layer of width $3h$, where we have ${\bf n} \cdot \nabla u_{\bf nn} = 0$, meaning that $u_{\bf nn}$ is constant along the characteristic direction. Once the second directional derivative $u_{\bf nn}$ is available in boundary layer from Eq.~\eqref{eq:extr_unn}, we first recover the first directional derivative $u_{\bf n}$ using Eq.~\eqref{eq:extr_un}, and then proceed with a second integration step in Eq.~\eqref{eq:extr_u} to obtain the extrapolated solution $u$. The extrapolated solution retains the order of accuracy, corresponding to the order of extrapolation, only within the prescribed narrow band near the boundary. Solutions obtained using quadratic extrapolation for various geometries  are shown in Figure~\ref{fig:aslam}. In our numerical simulations, we extrapolate the solutions in a band of width $3h$, while the extrapolated solutions in Figure~\ref{fig:aslam} are displayed in a region larger than this band, for graphical purposes. However, when a wider region is considered, central differences are not stable and one has to use upwind discretization to get a smooth extrapolation.
We show in Fig.~\ref{fig:aslam_xl} the extrapolated solution in a wider layer by using upwind discretization of the convective term.

To solve Eqs.~(\ref{eq:extr_unn}-\ref{eq:extr_u}), we first consider the variational formulation of the respective equations. As we did in Section~\ref{sec:model}, multiplying (\ref{eq:extr_unn}-\ref{eq:extr_u}) 
by a test function ${\bf w} \in {\bf V}$ and integrating over $\Omega$, we obtain the following system
\begin{align*} \int_\Omega \left(\pad{ {u_{\bf nn}}}{ \tau} + \chi \, {\bf n} \cdot \nabla u_{\bf nn} \right) {\bf w} \,{\rm d}\Omega & = 0 \\
\int_\Omega \left(\pad{ {u_{\bf n}}}{ \tau} + \chi \, \left({\bf n} \cdot \nabla u_{\bf n} - u_{\bf nn}\right) \right) {\bf w} \,{\rm d}\Omega & = 0 \\
\int_\Omega \left(\pad{ {u_{\bf nn}}}{ \tau} + \chi \, \left({\bf n} \cdot \nabla u_{\bf } - u_{\bf n}\right) \right) {\bf w} \,{\rm d}\Omega & = 0
\end{align*}
Analogously to variational formulation introduced in \eqref{pro:variational}, the integrals defined over $\Omega$ are then evaluated in $\Omega_h$.

\section{Numerical results}
\label{sec:results}
In this section, we test our numerical schemes. We begin with simulations on a stationary domain. Subsequently, we present numerical results in the context of moving domains. 

All numerical experiments have been implemented using \texttt{Gridap} (version 0.17.23), a Julia-based open-source software designed for approximating the solution of PDEs \cite{Badia2020}. \texttt{Gridap} leverages the Julia JIT (just-in-time) compiler and has a user-interface that closely resembles that whiteboard-style notation of the weak formulation \cite{Verdugo2022}. In particular, we used GridapEmbedded (version 0.9.5) \cite{GridapEmbedded-jl}, a plug-in that provides tools for implementing embedded finite element methods.

\subsection{Accuracy tests: Stationary domain}
\label{section:accuracy_stat}

\begin{figure}
    \centering
    \begin{subfigure}{0.32\textwidth}
        \includegraphics[width=\textwidth]{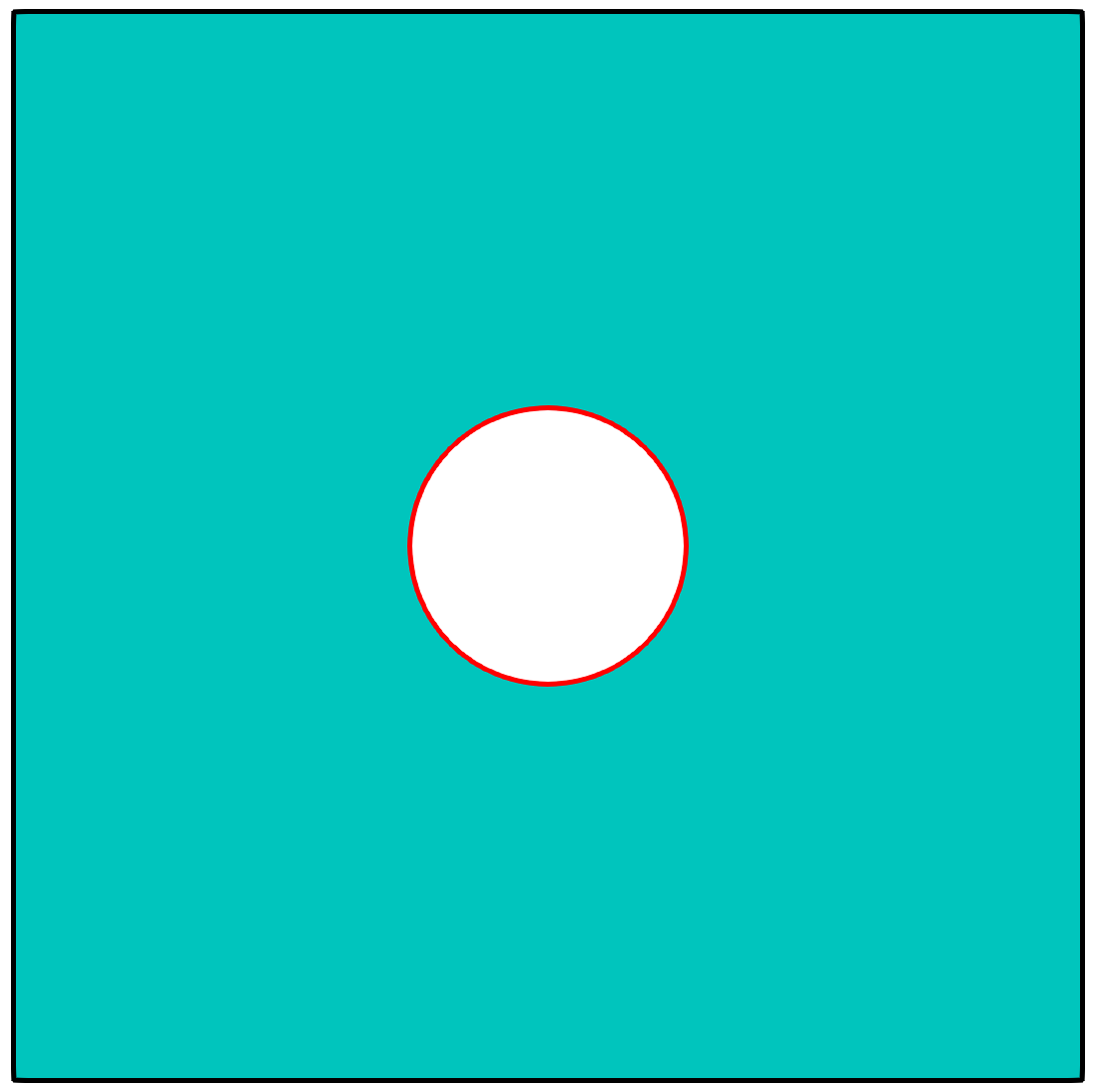}
        \caption{\textit{Disk}}
    \end{subfigure}
    \hfill 
    \begin{subfigure}{0.32\textwidth}
        \includegraphics[width=\textwidth]{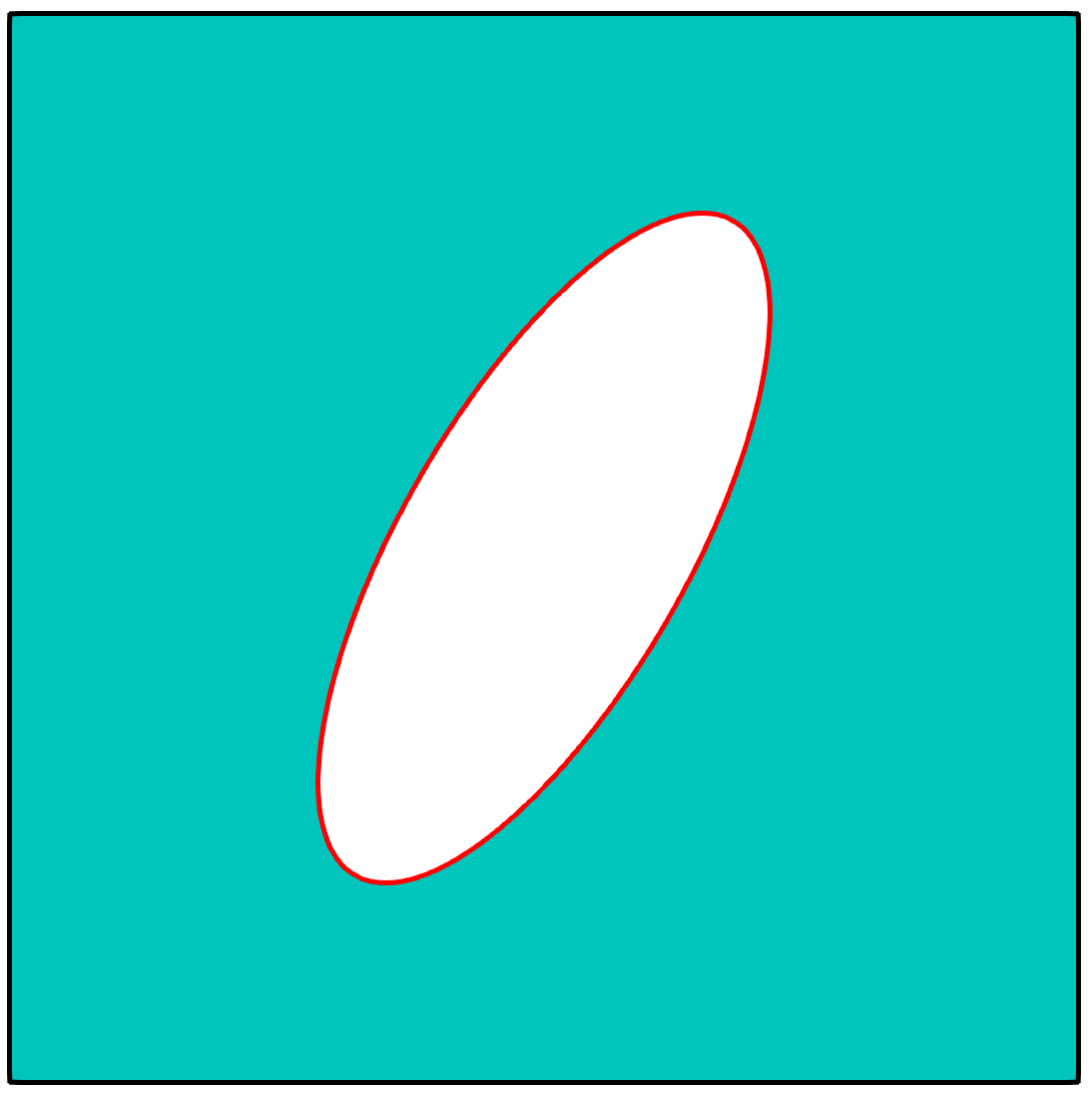}
        \caption{\textit{Ellipse}}
    \end{subfigure}
    \hfill
    \begin{subfigure}{0.32\textwidth}
        \includegraphics[width=\textwidth]{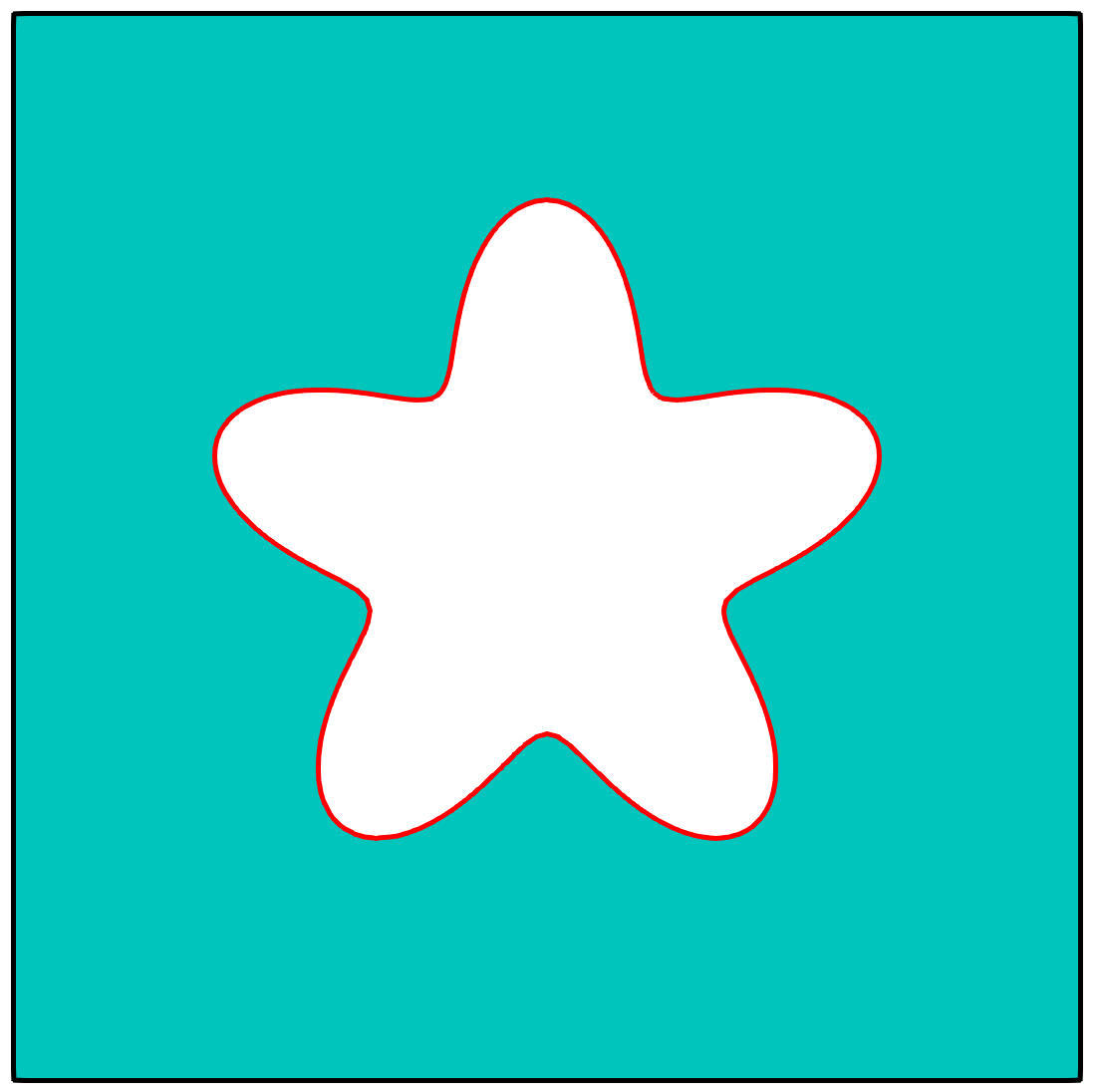}
        \caption{\textit{Flower}}
    \end{subfigure}
    \caption{\textit{Geometries used for tests.}}
    \label{fig:geos}
\end{figure}

In this section, we consider a time independent domain, i.e. $\Omega(t) = \Omega$. We consider the Kim and Moin problem, for which an exact solution is available. Let us consider the square region  $R = [-1,1]^2$ and the Navier-Stokes problem with Dirichlet conditions, whose exact solution is
\begin{subequations}
    \begin{align}
        u^{\rm exa} & = \sin(2\pi x) \cos(2\pi y)e^{(-2\nu  t)} \\
        v^{\rm exa} & = -\cos(2\pi x) \sin(2\pi y)e^{(-2\nu  t)} \\
        p^{\rm exa} & = -0.25 \left( \cos(4\pi x) + \cos(4\pi y) \right)e^{(-4\nu  t)}
    \end{align}
\end{subequations}
We take $\nu=1$, corresponding to $Re=1$. The domains we consider are defined by the following level set function (see Fig.~\ref{fig:geos}) : 
\begin{enumerate}
    \item[(a)] Disk: 
    \begin{equation}\label{phi:disk}
        \phi(x,y) = \frac{1}{\sqrt{15}} -\sqrt{x^2 + y^2} 
    \end{equation}
    \item[(b)] Ellipse: 
    \begin{align}\label{phi:ellipse}
        A = \frac{1}{\sqrt{14}}, \quad B = \frac{1}{\sqrt{2}}, \quad \theta = \frac{\pi}{6}, \\
      \phi(x,y) = 1-\frac{(\cos(\theta) x-\sin(\theta)y)^2}{A^2}  - \frac{(\sin(\theta)x+\cos(\theta)y)^2}{B^2} 
    \end{align}
    \item[(c)] Flower: 
    \begin{align}\label{phi:flower}
         A = 0.5, \quad B = 0.15,\quad  R = \sqrt{x^2 + y^2}, \\
			\phi(x,y) =   A -R + \frac{B(y^5 + 5x^4y - 10x^2y^3)}{  R^5}   
    \end{align}
\end{enumerate}

\begin{figure}
    \centering
    \begin{subfigure}{0.32\textwidth}
        \includegraphics[width=\textwidth]{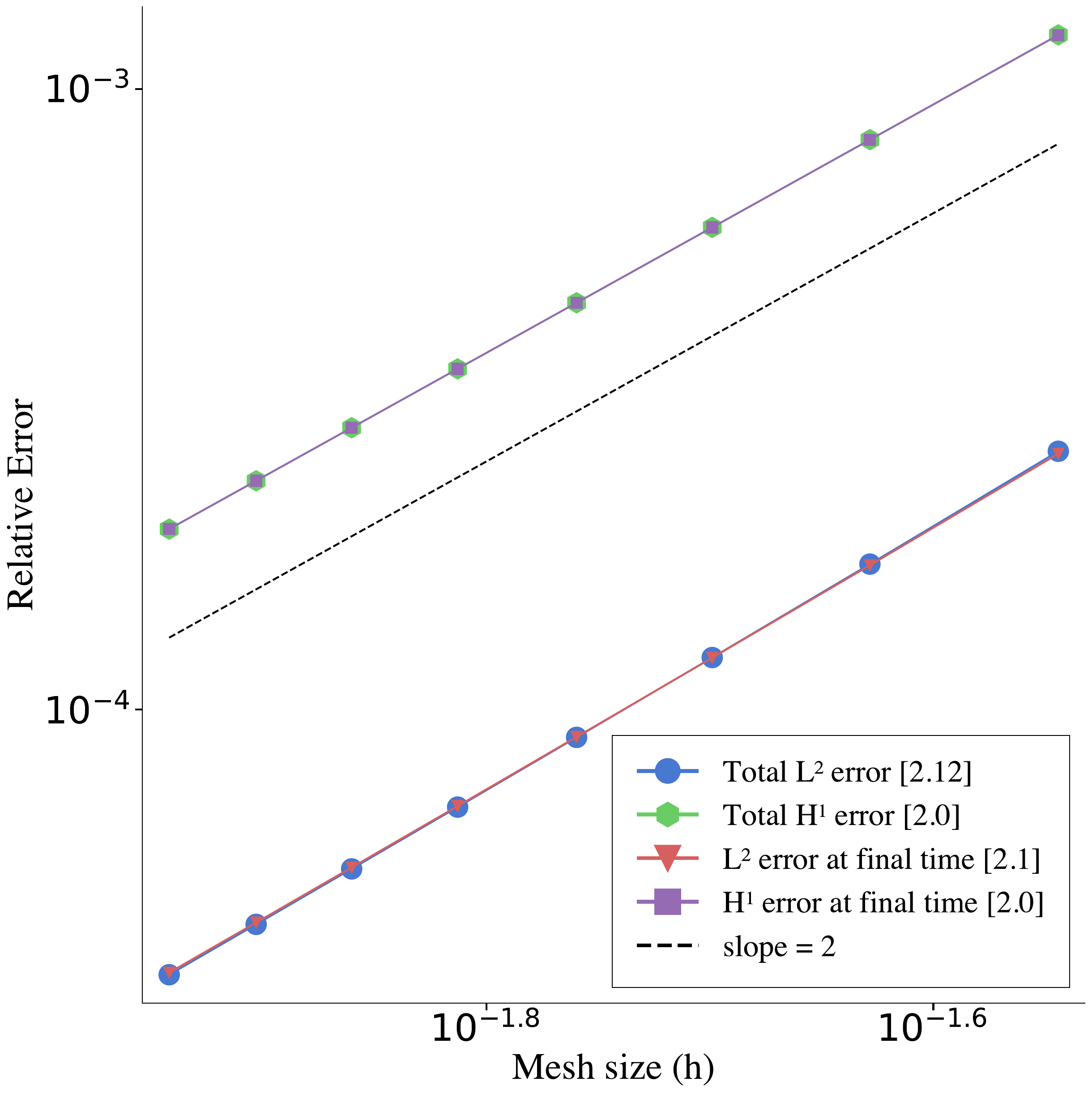}
        \caption{\textit{Disk}}
    \end{subfigure}
    \hfill
\begin{subfigure}{0.32\textwidth}
        \includegraphics[width=\textwidth]{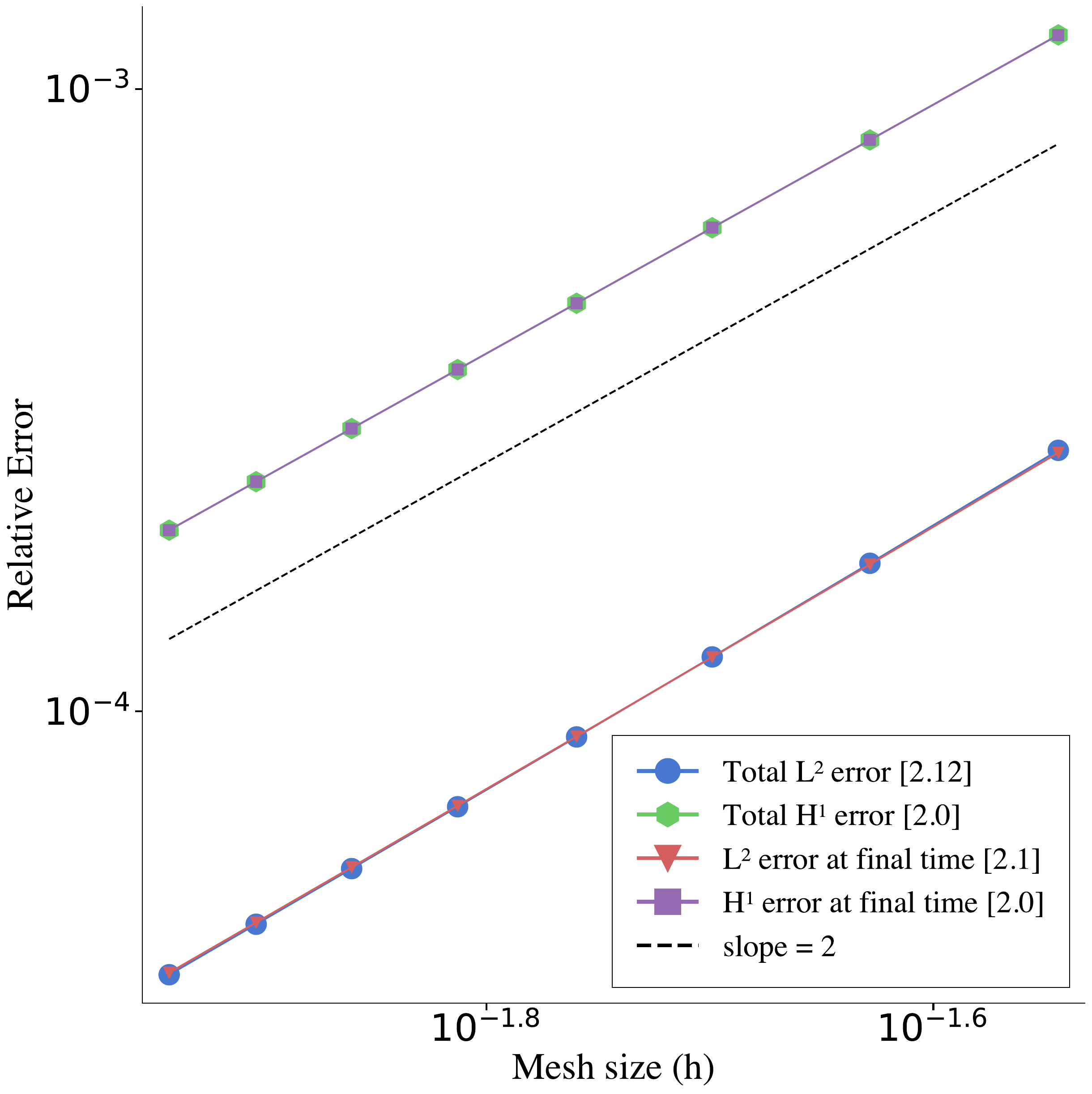}
        \caption{\textit{Ellipse}}
    \end{subfigure}
    \hfill
\begin{subfigure}{0.32\textwidth}
    \includegraphics[width=\textwidth]{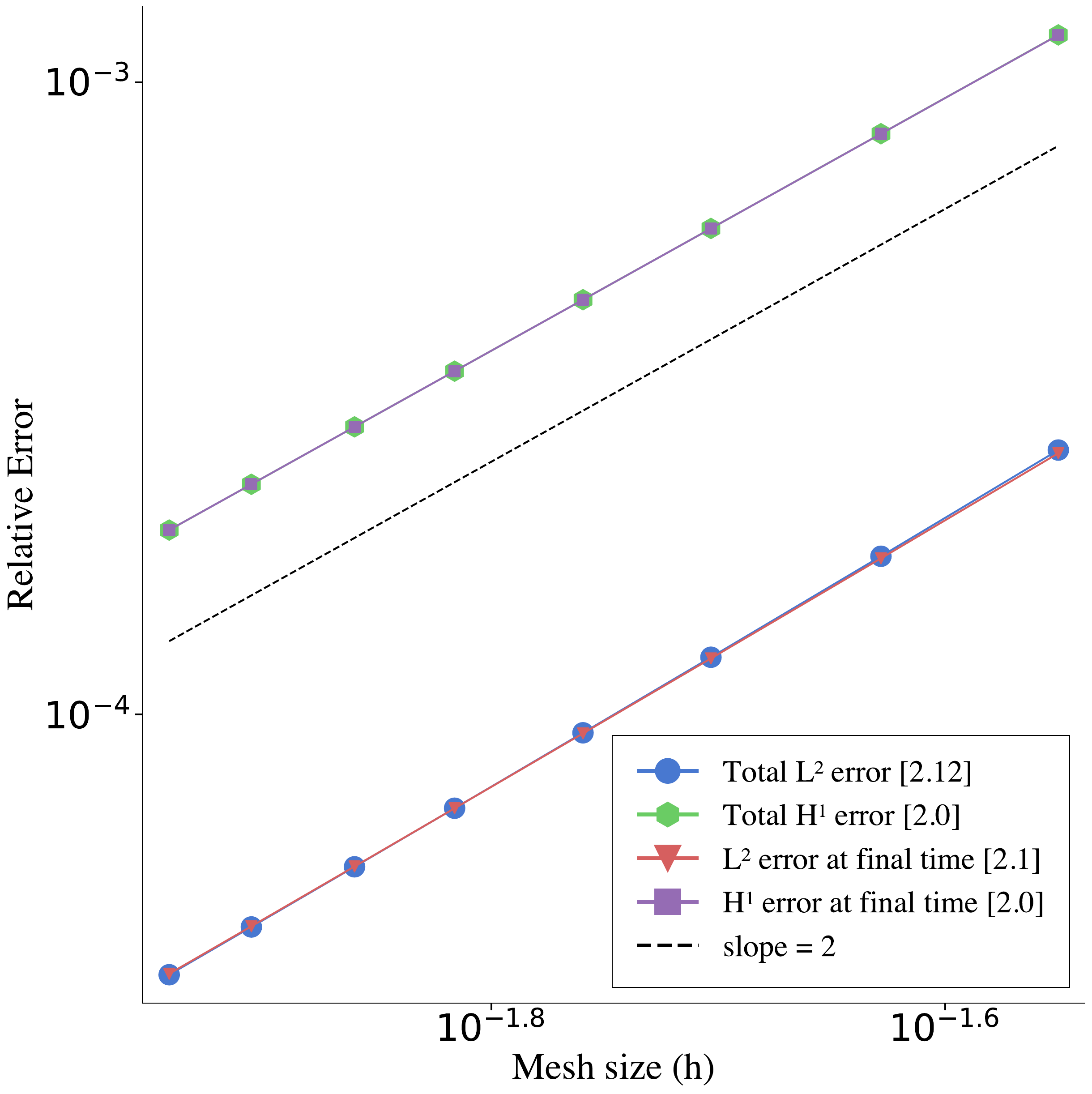}
    \caption{\textit{Flower}}
\end{subfigure}
    \caption{\textit{Plots of relative errors against spatial mesh size $h$ for a two-stage second-order IMEX scheme for Navier-Stokes equations on various stationary geometries.}} 
    \label{fig:imex_ns_2}
\end{figure}

\begin{figure}
  \centering
    \begin{subfigure}{0.32\textwidth}
        \includegraphics[width=\textwidth]{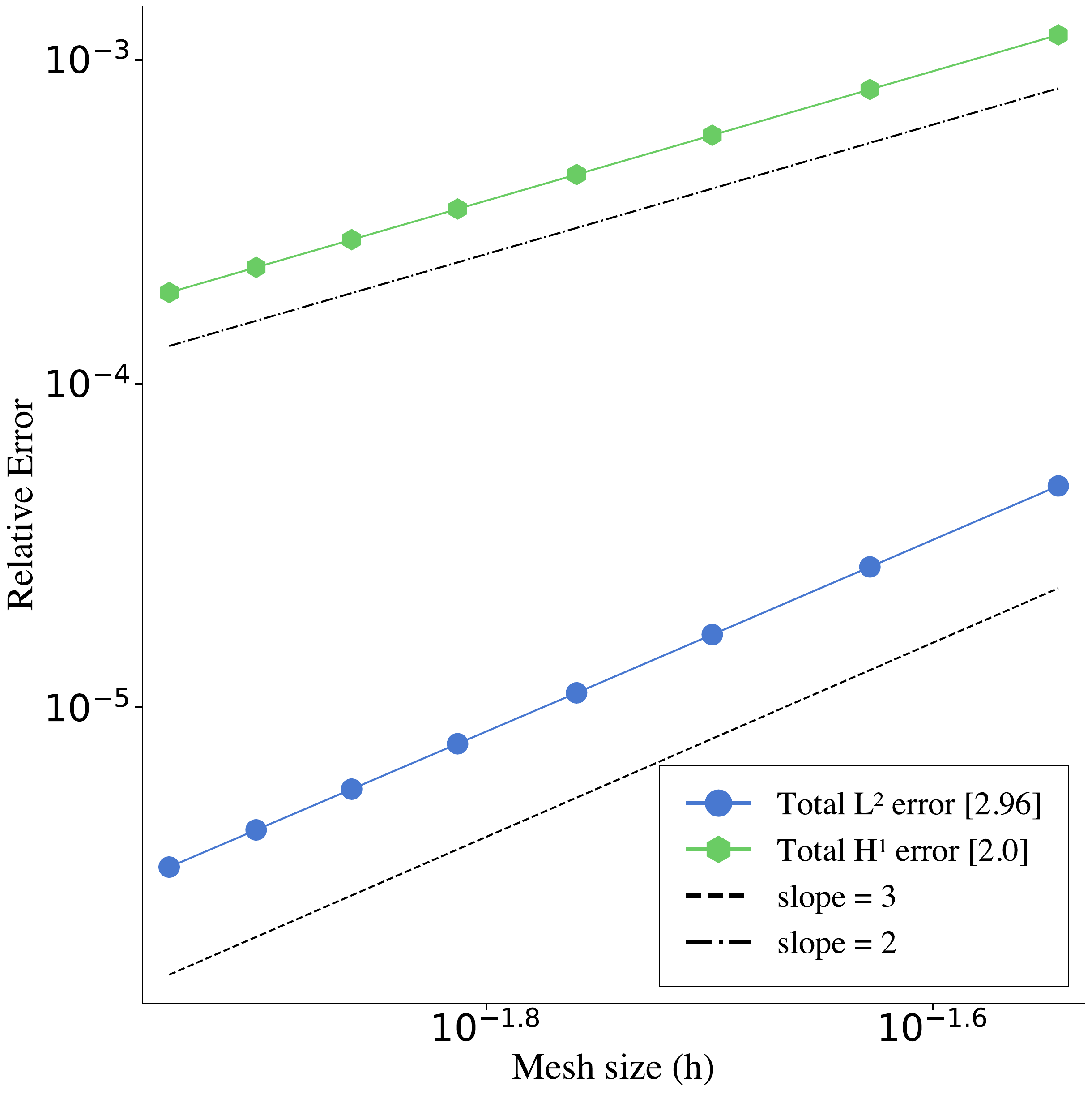}
        \caption{\textit{Disk}}
    \end{subfigure}
    \hfill
\begin{subfigure}{0.32\textwidth}
        \includegraphics[width=\textwidth]{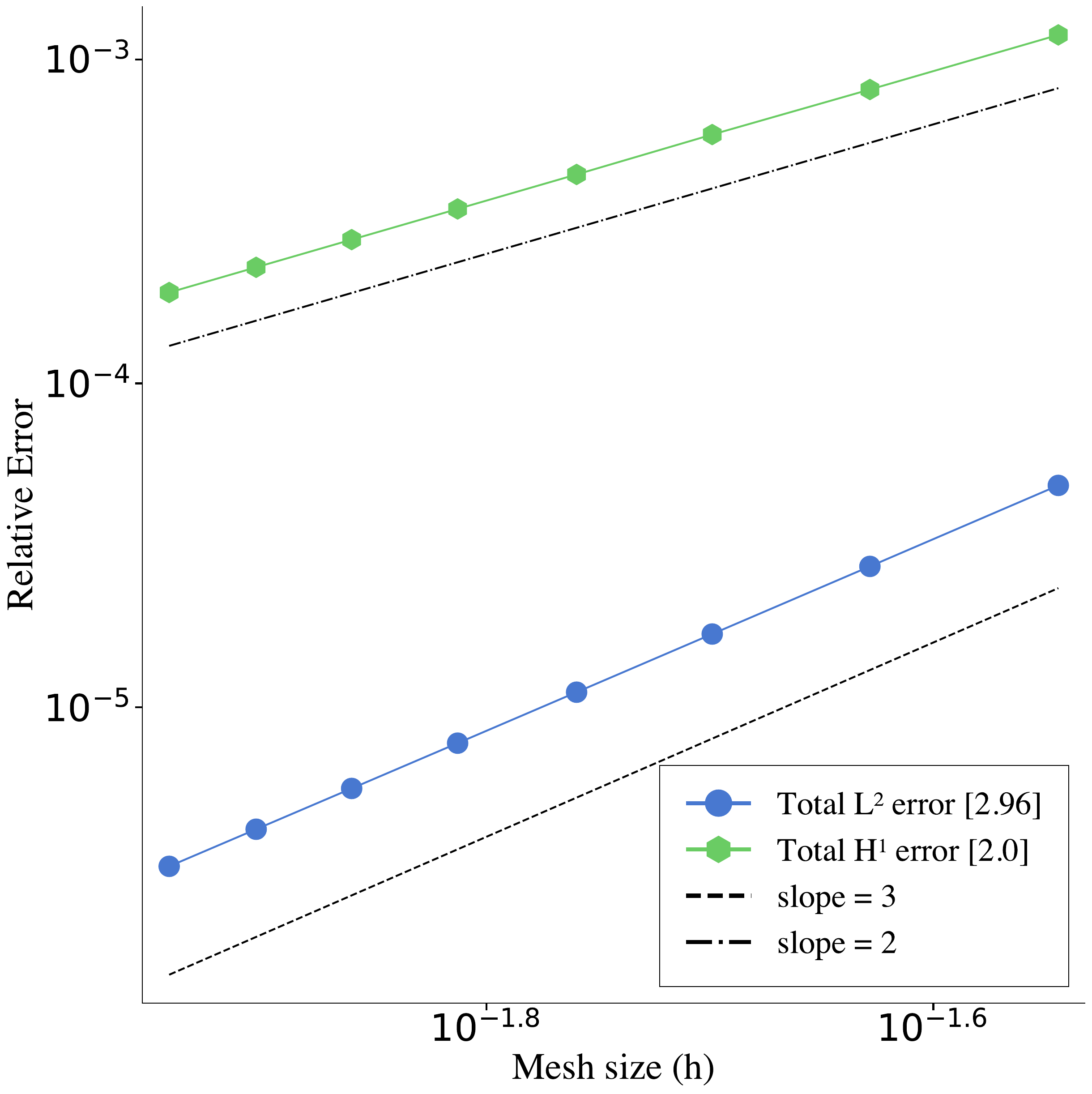}
        \caption{\textit{Ellipse}}
    \end{subfigure}
    \hfill
\begin{subfigure}{0.32\textwidth}
    \includegraphics[width=\textwidth]{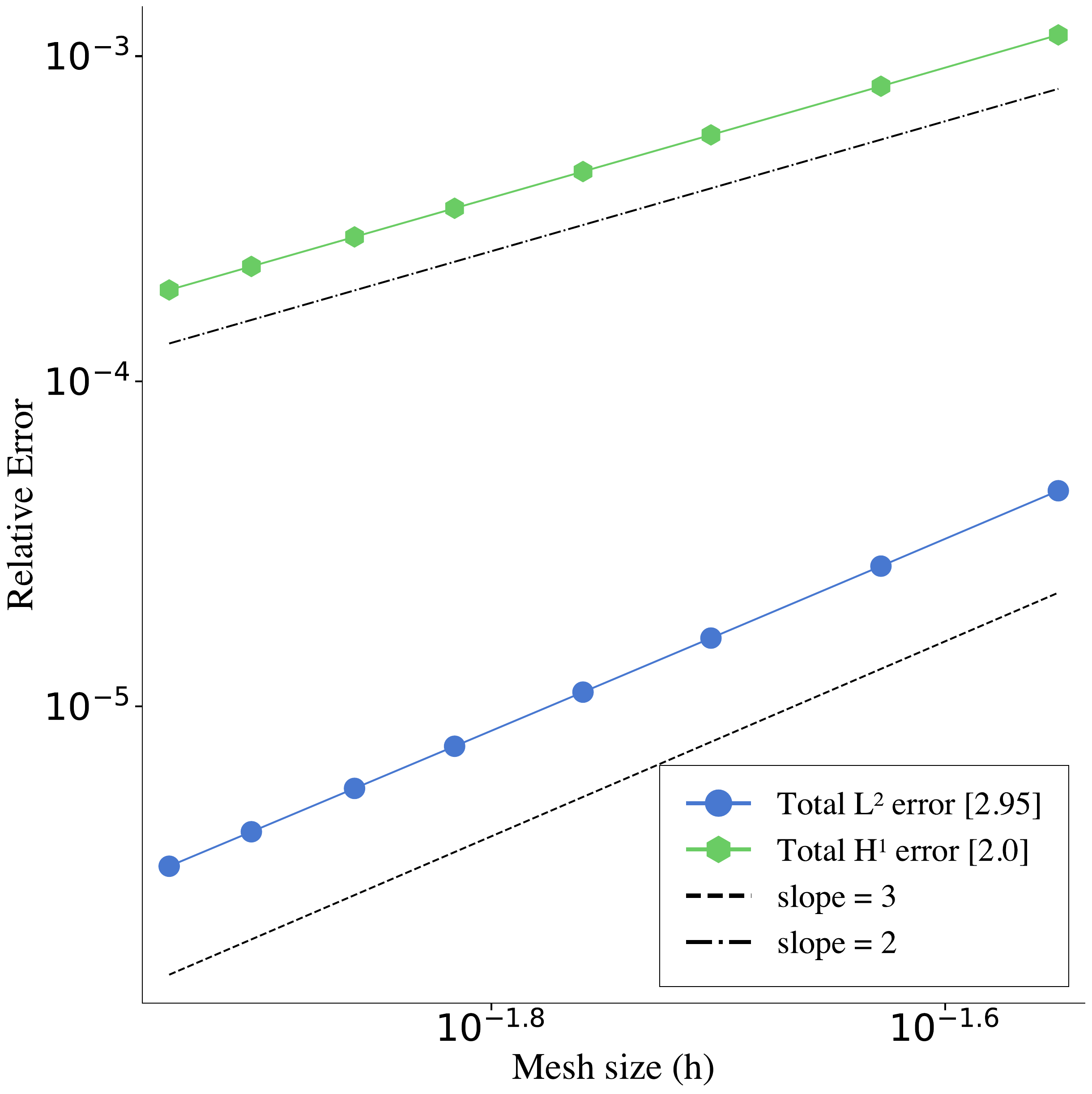}
    \caption{\textit{Flower}}
\end{subfigure}

    \caption{\textit{Plots of relative errors against spatial mesh size $h$ for a four-stage third-order IMEX scheme for Navier-Stokes equations on different stationary geometries.}}
    \label{fig:imex_ns_3}
\end{figure}

In Fig.~\ref{fig:imex_ns_2}, we show the plot of relative errors against spatial mesh size $h$ for various stationary geometries, choosing the two-step IMEX scheme. The time step size is chosen as $\Delta t = h$.
The following quantities are computed and plotted:
\begin{align*}
   {\rm Total }\, L^2 \, {\rm error} & =  \frac{\left(\Sigma_{n-h=1}^{N_{ts}} \|\boldsymbol{u}_{exa}(t^n)-\boldsymbol{u}_h^n\|^2_{L^2(\Omega_h)} \right)^{1/2}}{\left(\Sigma_{n=1}^{N_{ts}} \|\boldsymbol{u}_{exa}(t^n)\|^2_{L^2(\Omega_h)}\right)^{1/2}}, \quad  {\rm Total }\, H^1 \, {\rm error} = \frac{\left(\Sigma_{n=1}^{N_{ts}} \|\boldsymbol{u}_{exa}(t^n)-\boldsymbol{u}_h^n\|^2_{H^1(\Omega_h)}  \right)^{1/2}}{\left(\Sigma_{n=1}^{N_{ts}} \|\boldsymbol{u}_{exa}(t^n)\|^2_{H^1(\Omega_h)}\right)^{1/2}}, \\
   L^2 \, {\rm error \, at \, final\, time} & = \frac{ \|\boldsymbol{u}_{exa}(t^{N_{ts}})-\boldsymbol{u}_h^{N_{ts}}\|_{L^2(\Omega_h)} }{ \|\boldsymbol{u}_{exa}(t^{N_{ts}})\|_{L^2(\Omega_h)}}, \quad  H^1 \, {\rm error \, at \, final\, time} = \frac{\|\boldsymbol{u}_{exa}(t^{N_{ts}})-\boldsymbol{u}_h^{N_{ts}}\|_{H^1(\Omega_h)} }{\|\boldsymbol{u}_{exa}(t^{N_{ts}})\|_{H^1(\Omega_h)}},
\end{align*} 
where $\|u\|_{L^2({\Omega_h})}^2=\int_{\Omega_h}(u\cdot u)d\Omega$ and $\|u\|_{H^1(\Omega_h)}^2=\|u\|_{L^2(\Omega_h)}^2+\int_{\Omega_h}( \nabla{u}\cdot\nabla{u})d\Omega$.
We observe from Fig.~\ref{fig:imex_ns_2} that
the final-time errors are indistinguishable from the total ones, meaning that the scheme is stable over time. Therefore, for all the following tests we show only the total errors (in $L^2$ and $H^1$ norms). 
In Fig.~\ref{fig:imex_ns_3}, we repeat the same test, using the four-step IMEX scheme. For this choice of IMEX scheme, we set $\Delta t = h/2$. As shown in Fig.~\ref{fig:imex_ns_2}, the accuracy tests on stationary domains indicate that, for the two-stage second-order IMEX scheme, both the $L^2$ and $H^1$ errors decrease with a rate of 2, consistent with theoretical expectations. As shown in Fig.~\ref{fig:imex_ns_3}, for the four-stage third-order IMEX scheme on stationary domains, the $L^2$ error exhibits third-order convergence, while the $H^1$ error converges at second order, as expected. Overall, these results confirm that the numerical schemes achieve their expected orders of accuracy.

\subsection{Accuracy tests: Moving domain}
In this section, we perform numerical tests on domains that depend on time. The evolution of the moving domain is defined by time-dependent level set function $\phi(x,y,t)=\phi_0(\tilde{x}(t),\tilde{y}(t))$, where $\phi_0(x,y)$ represents the domain at time $t=0$, and
\begin{equation}\label{eq:rotxy}
\begin{pmatrix}
    \tilde{x}(t) \\
    \tilde{y}(t)
\end{pmatrix} = \begin{pmatrix}
    \cos(\omega t) &-\sin(\omega t) \\
    \sin(\omega t) & \cos(\omega t)
\end{pmatrix} \begin{pmatrix}
    x \\
    y
\end{pmatrix},
\end{equation}
with an angular velocity of $\omega = 2 \pi / 5$.
We test rotating ellipse and flower, then the initial level set $\phi_0(x,y)$ is given by Eq.~\eqref{phi:ellipse} (ellipse) or Eq.~\eqref{phi:flower} (flower).

\begin{figure}
  \centering
\begin{subfigure}{0.45\textwidth}
       \includegraphics[width=\textwidth]{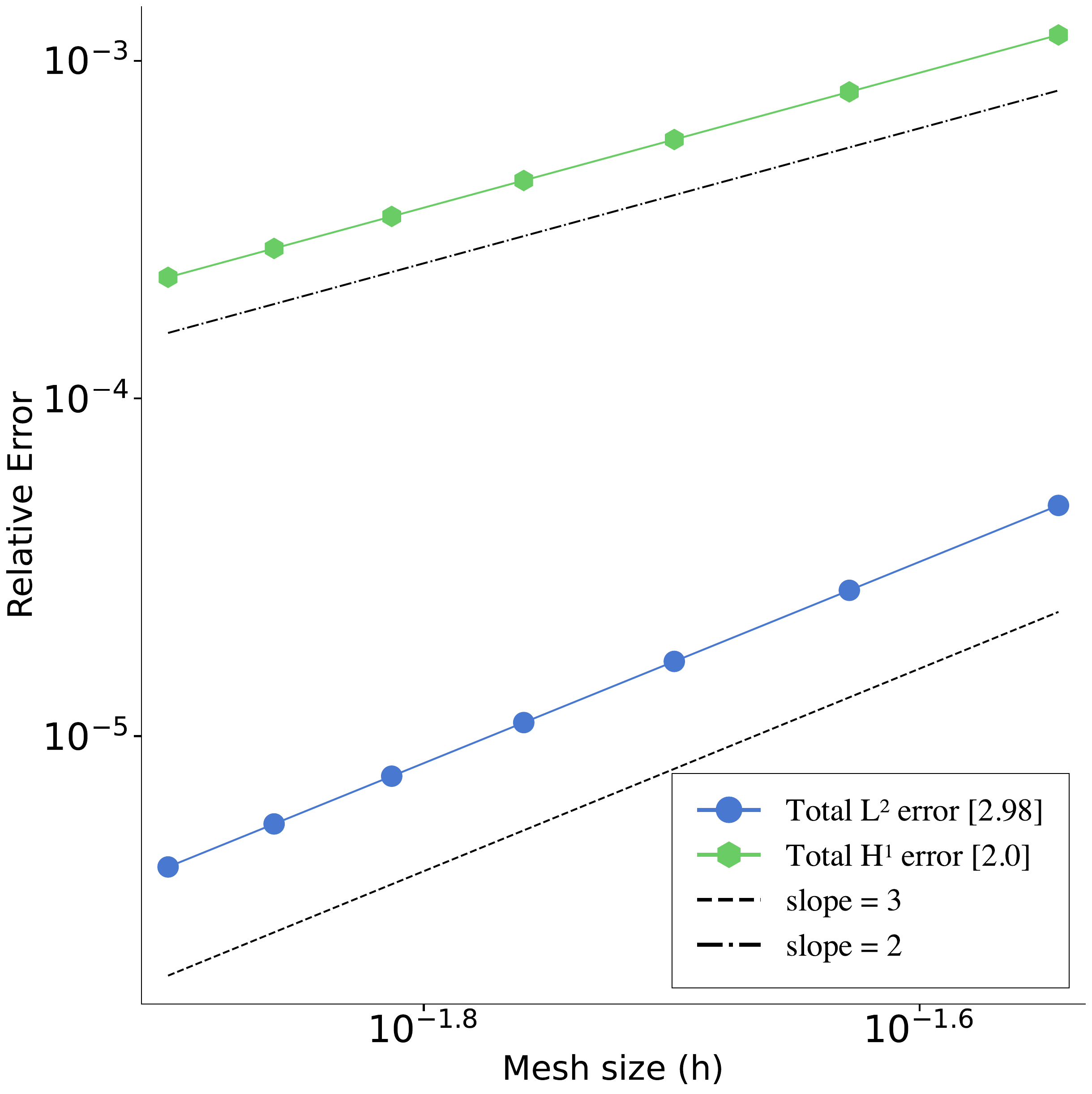}
    \caption{\textit{Ellipse}}
    \end{subfigure}
    \hfill
\begin{subfigure}{0.45\textwidth}
    \includegraphics[width=\textwidth]{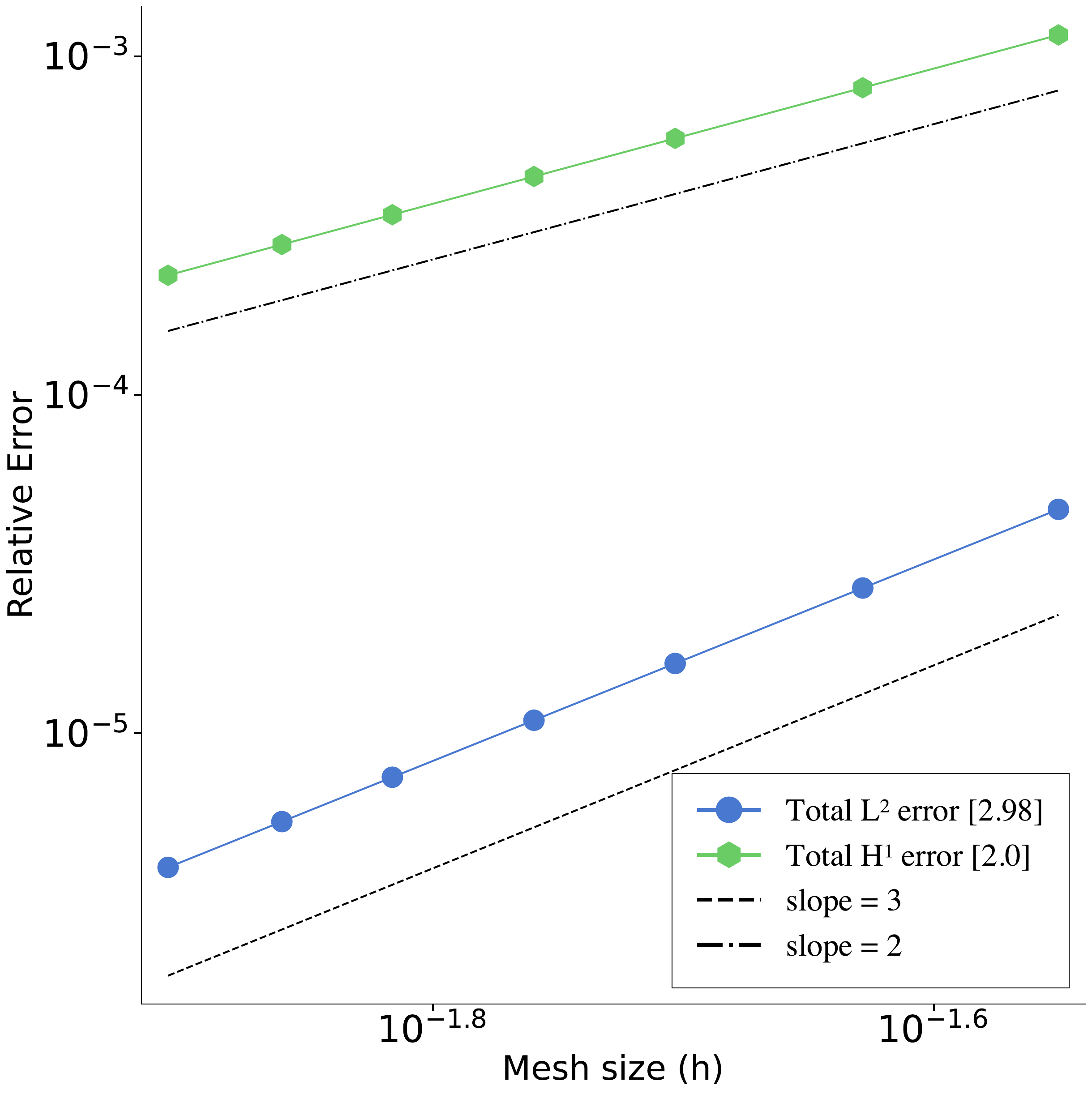}
    \caption{\textit{Flower}}
\end{subfigure}

    \caption{\textit{Plots of relative errors against spatial mesh size $h$ for a four-stage third-order IMEX scheme for Navier-Stokes equations on different moving geometries.}}
    \label{fig:imex_ns_mv}
\end{figure}

Relative errors against spatial mesh sizes $h$ are illustrated in Fig.~\ref{fig:imex_ns_mv}.
The four-step IMEX scheme attains third-order convergence in the 
$L^2$ norm and second-order convergence in the $H^1$ norm, in agreement with theoretical expectations, even on moving geometries.
\subsection{Flow around a disk}
\begin{figure}
    \centering
    \begin{subfigure}{\textwidth}
       \includegraphics[width=\textwidth]{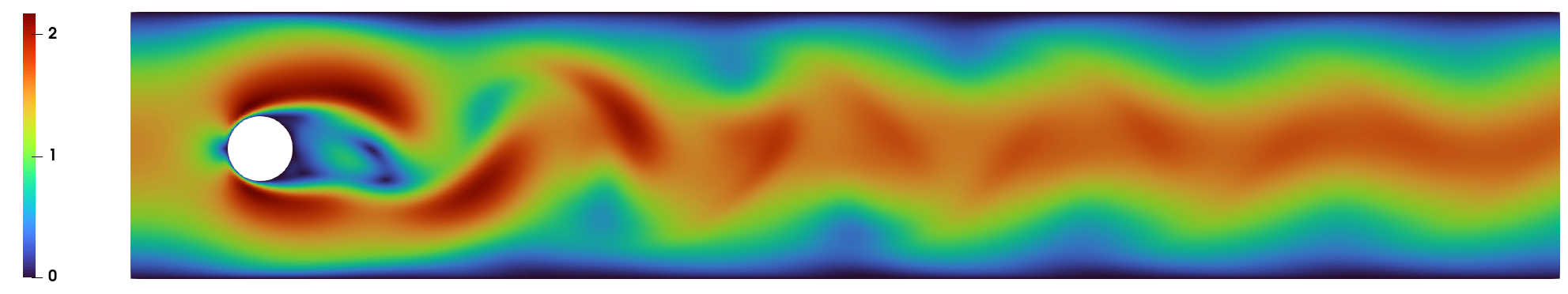}
        \caption{\textit{Magnitude of the velocity field}}
   \end{subfigure}
   \hfill
   \begin{subfigure}{\textwidth}
       \includegraphics[width=\textwidth]{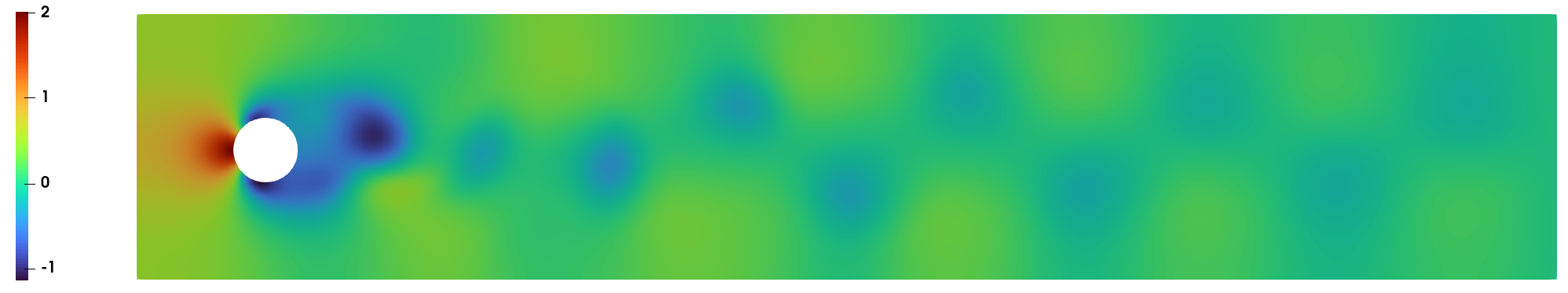}
        \caption{\textit{Pressure field}}
   \end{subfigure}
    \caption{\textit{The magnitude of the velocity and pressure fields for flow around a disk at time = $10s$ for $Re =100$.}}
    \label{fig:flow_cyl}
\end{figure}

In this test, we study the unsteady two-dimensional flow around a circular disk at Reynolds number $Re=100$ \cite{Schfer1996}. The computational domain is given by $\Omega =[0,2.2]\times [0,0.41]\setminus B_R$, where $B_R$ is a disk centered at (0.2,0.2) with radius $R_B=0.05$. 
No-slip boundary conditions are imposed on the lower and upper walls $\Gamma_1 = [0,2.2]\times \{0\}$ and $\Gamma_2 = [0,2.2]\times \{0.41\}$ as well as on the boundary of the disk $\Gamma_3 = \partial B_R$, that is 
\[ {\bf{u}}|_{ \Gamma_1} = {\bf{u}}|_{ \Gamma_2}= {\bf{u}}|_{ \Gamma_3} = 0.\] 
On the inlet $\Gamma_4 = \{0\} \times [0,0.41]$, a parabolic inflow is prescribed as ${\bf{u}}(0,y) = \left(\frac{4Uy(0.41-y)}{0.41^2},0 \right)$, where $U =1.5$ denotes the maximum inflow velocity. At the outlet $\Gamma_5 = \{2.2\} \times [0,0.41]$, the do-nothing boundary condition is applied $\nu \partial_{\bf{n}} {\bf{u}} - p{\bf{n}} =0,$ where ${\bf{n}}$ denotes the outward normal. 
Following \cite{Schfer1996} $\nu = 0.001, Re = \frac{2}{3}u(0,0.41/2)(2R_B)/\nu$.
The final simulation time is set to $T = 10s$, with a time step $\Delta t = 0.01s$. The velocity and pressure fields at the final time are shown in Fig.~\ref{fig:flow_cyl}. The results exhibit the formation of alternating vortices in the wake of the cylinder, which is consistent with the expected behavior at $Re=100$. 

\par The drag and lift forces on the cylinder are computed as 
\[ \left(F_D,F_L\right) = \int_{\Gamma_3} \left( \nu \nabla \mathbf{u} - p \mathbf{I} \right) \mathbf{n} \mathrm{~d}\mathbf{S},\] 
where \(\mathbf{n}\) denotes the outward unit normal on the disk boundary \(\Gamma_3\). The drag and lift coefficients are defined as 
\[ C_D = \frac{2}{U_{\mathrm{mean}}^2 L}F_D, \qquad C_L = \frac{2}{U_{\mathrm{mean}}^2 L}F_L,\]
where $U_{\mathrm{mean}} = 1$ is the mean velocity of the parabolic inflow velocity and $L = 2r = 0.1$ is the characteristic length of the flow configuration. \par
\begin{figure}
    \centering
    \begin{subfigure}{0.32\textwidth}
       \includegraphics[width=\textwidth]{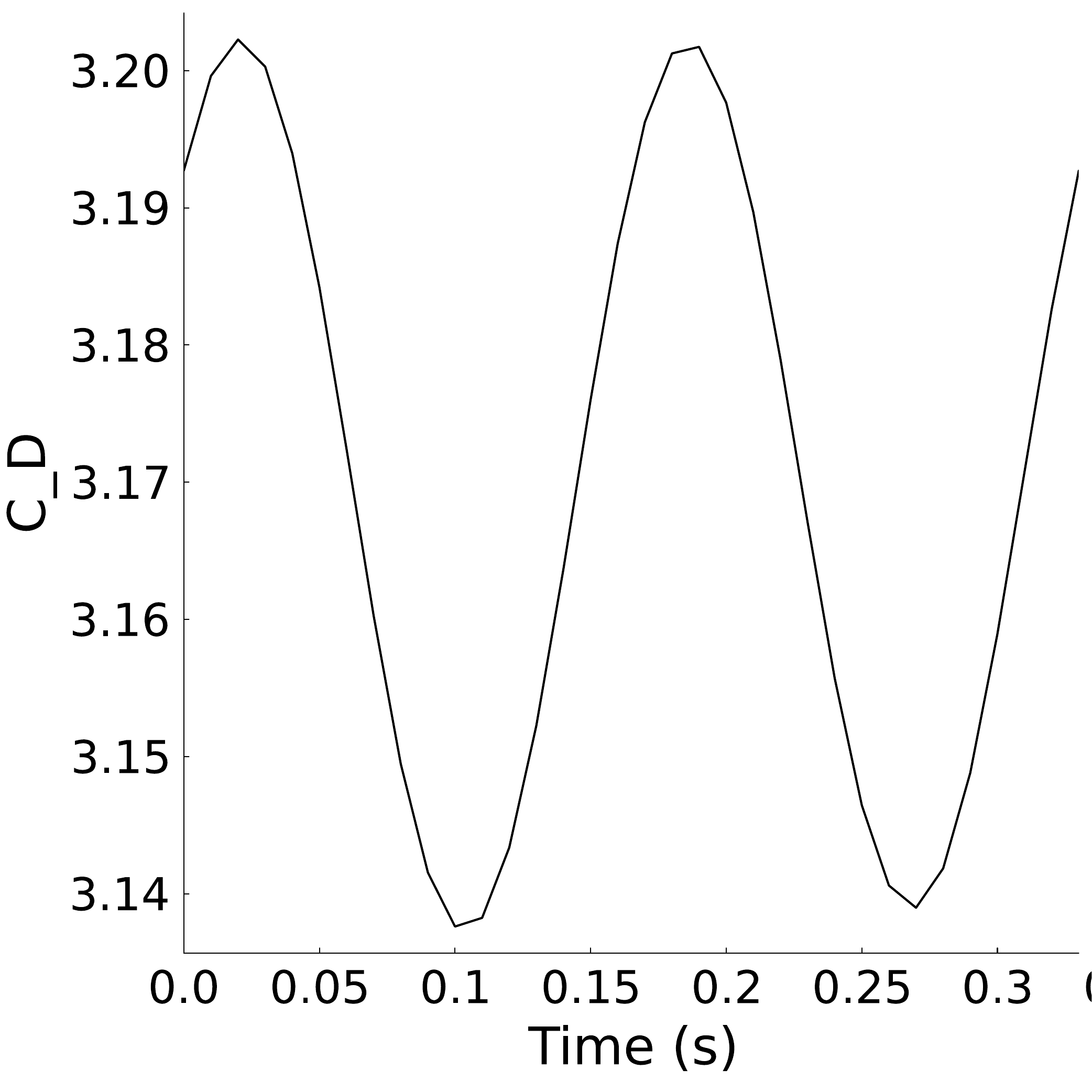}
        \caption{\textit{Drag coefficient}}
   \end{subfigure}
   \hfill
   \begin{subfigure}{0.32\textwidth}
       \includegraphics[width=\textwidth]{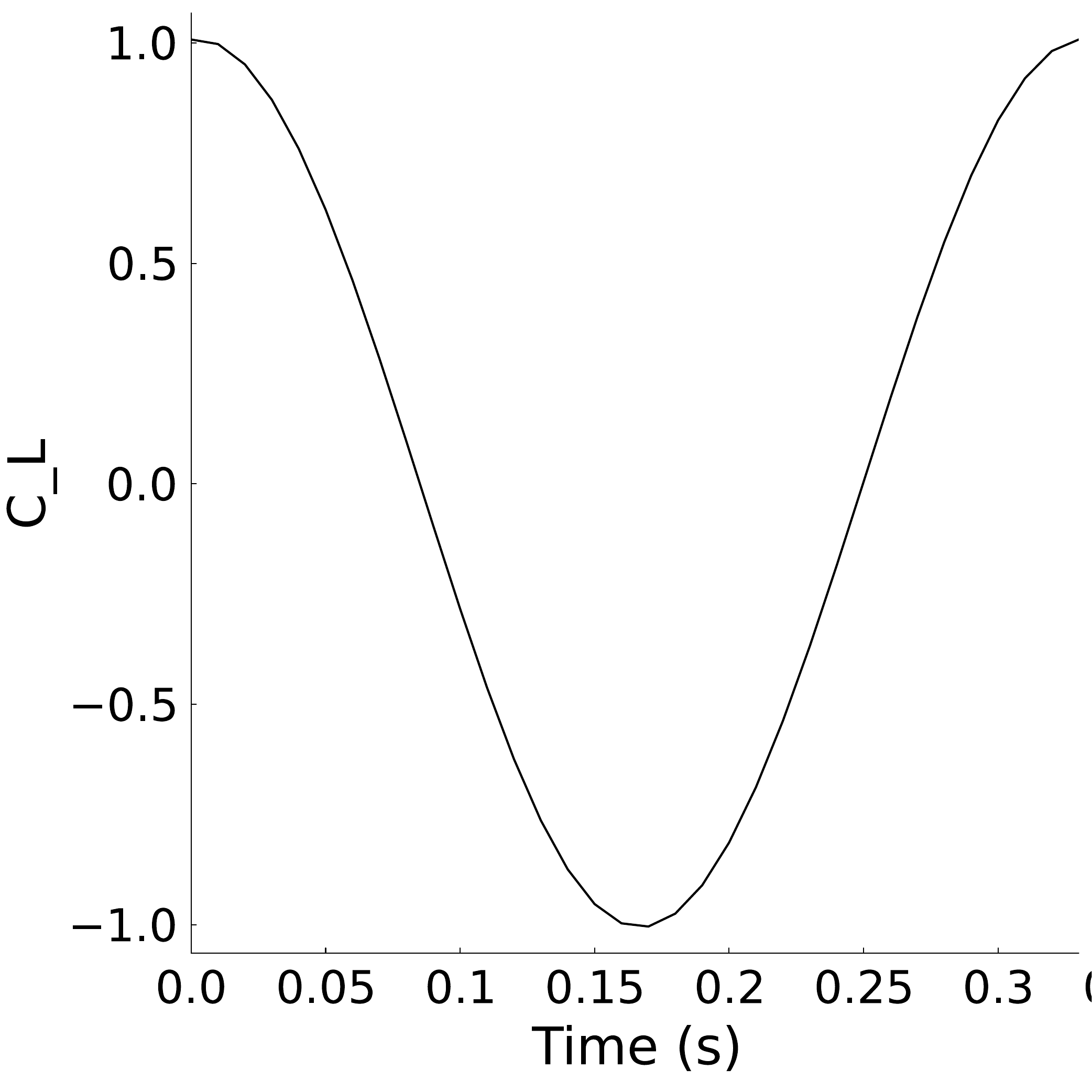}
        \caption{\textit{Lift coefficient}}
   \end{subfigure}
   \hfill
   \begin{subfigure}{0.32\textwidth}
       \includegraphics[width=\textwidth]{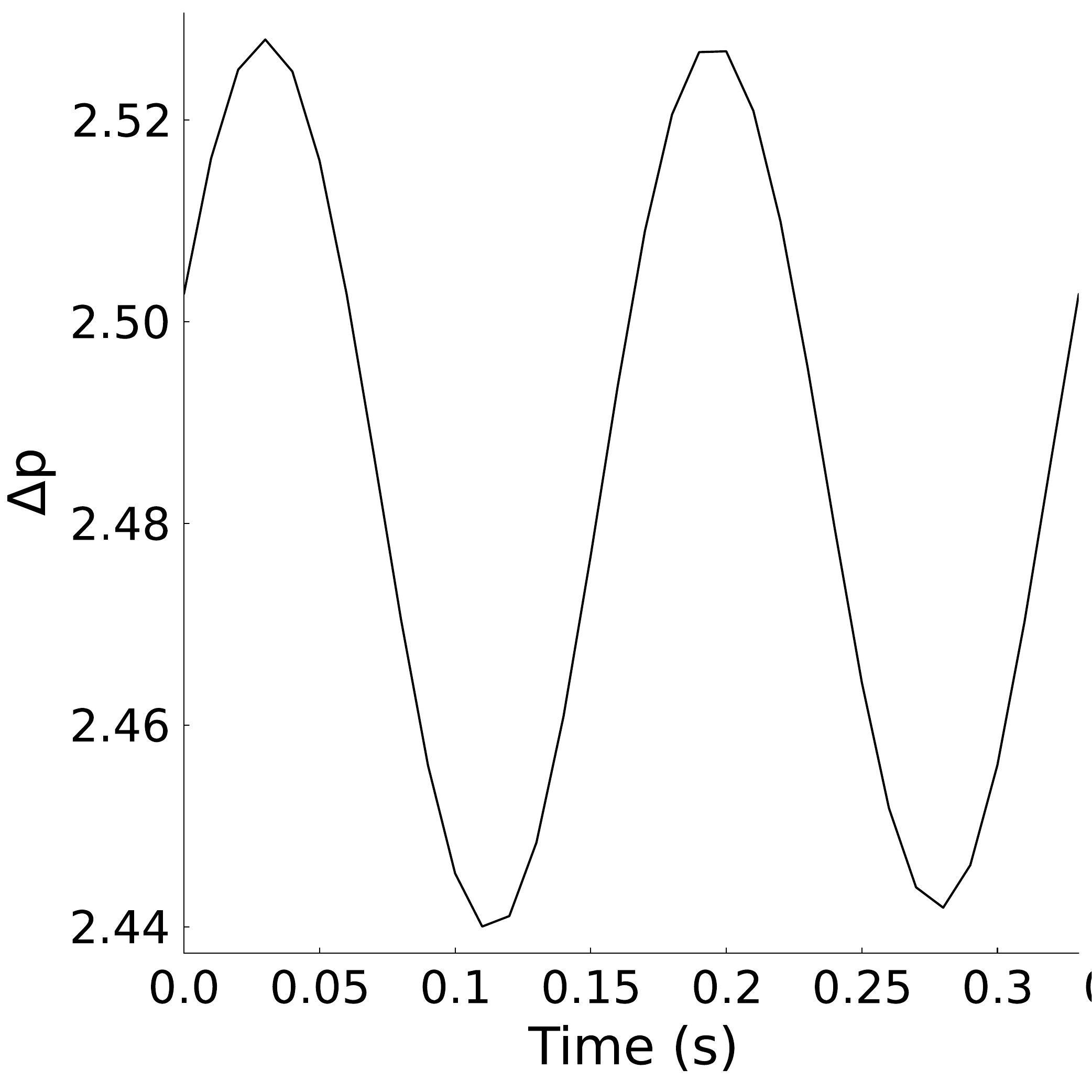}
        \caption{\textit{Pressure difference}}
   \end{subfigure}
    \caption{\textit{Plots of the measured numerical quantities over the time interval $[t_0, t_0 + 1/f]$ for a flow around a disk.}}
    \label{fig:drag_lift_pres_cyl}
\end{figure}

 The following quantities are computed: the drag coefficient $C_D$, the lift coefficient $C_L$, and the pressure difference $\Delta P$
as functions of time over one period $[t_0,\, t_0 + 1/f]$, where $f = f(C_L)$ denotes the frequency of the lift coefficient;
the maximum drag coefficient $C_{D,\max}$ and the maximum lift coefficient $C_{L,\max}$; the Strouhal number $\mathrm{St}$; the pressure difference $\Delta p(t)$ evaluated at
$t = t_0 + \frac{1}{2f}$. The initial data at $t = t_0$ corresponds to a flow state with $C_{L,\max}$. The numerical measured quantities are presented in Tab.~\ref{tab:flow_cyl_meas_quan}. Fig.~\ref{fig:drag_lift_pres_cyl} depicts the drag and lift coefficients as well as the pressure difference $\Delta p =  p(a_1)-p(a_2)$ between the points $a_1=(0.15,0.2)$ on the front of the disk and $a_2=(0.25,0.2)$ on the rear of the disk over the interval $[t_0,\, t_0 + 1/f]$. 

\begin{table}
\centering
\begin{tabular}{|c|c|c|c|c|c|c|}
\hline
Quantity & \# Elements & \# Unknowns & $C_{D,\max}$& $C_{L,\max}$ & $\text{St}$ & $\Delta p$  \\ 
\hline
Reference value & & & 3.22 - 3.24 &  0.99 - 1.01& 0.295 - 0.305& 2.46 - 2.50 \\
\hline
Measured value& 11668 &104407 & 3.20227 & 1.00773&0.30303&2.50133  \\
\hline

\end{tabular}
\caption{\textit{Numerically measured quantities for a flow past a disk}}
\label{tab:flow_cyl_meas_quan}
\end{table}

The numerically measured results are in good agreement with the reference ranges for the flow past a disk taken from \cite{Schfer1996},
and the periodic drag and lift coefficients as well as the pressure difference exhibit stable amplitudes, confirming that the simulation accurately captures the unsteady flow dynamics.

\subsection{Driven cavity}
\begin{figure}
    \centering
   \begin{subfigure}{0.45\textwidth}
       \includegraphics[width=\textwidth]{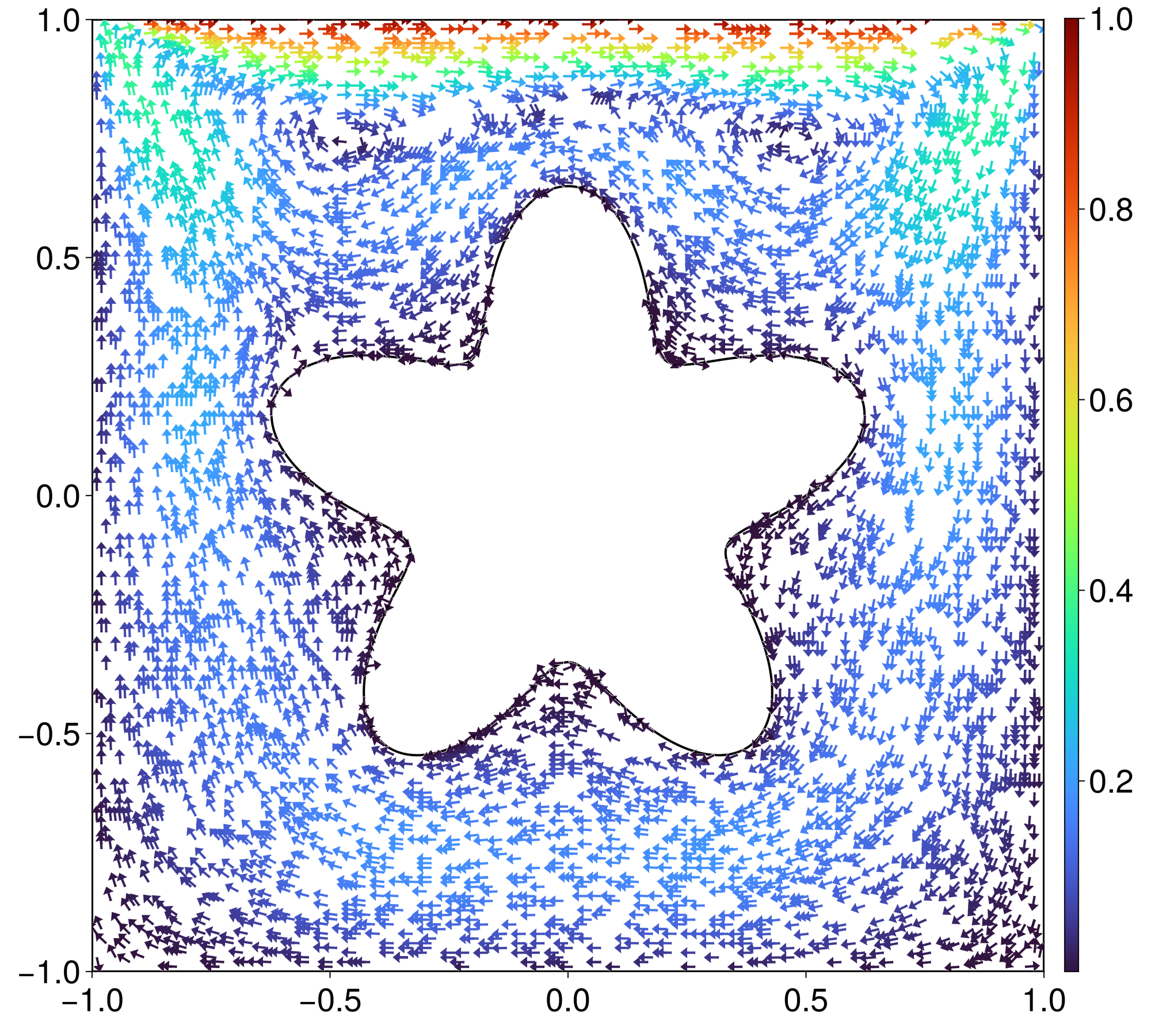}
       \caption{$Re = 1$,\quad $t = 10s$}
   \end{subfigure}
    \hfill 
   \begin{subfigure}{0.45\textwidth}
       \includegraphics[width=\textwidth]{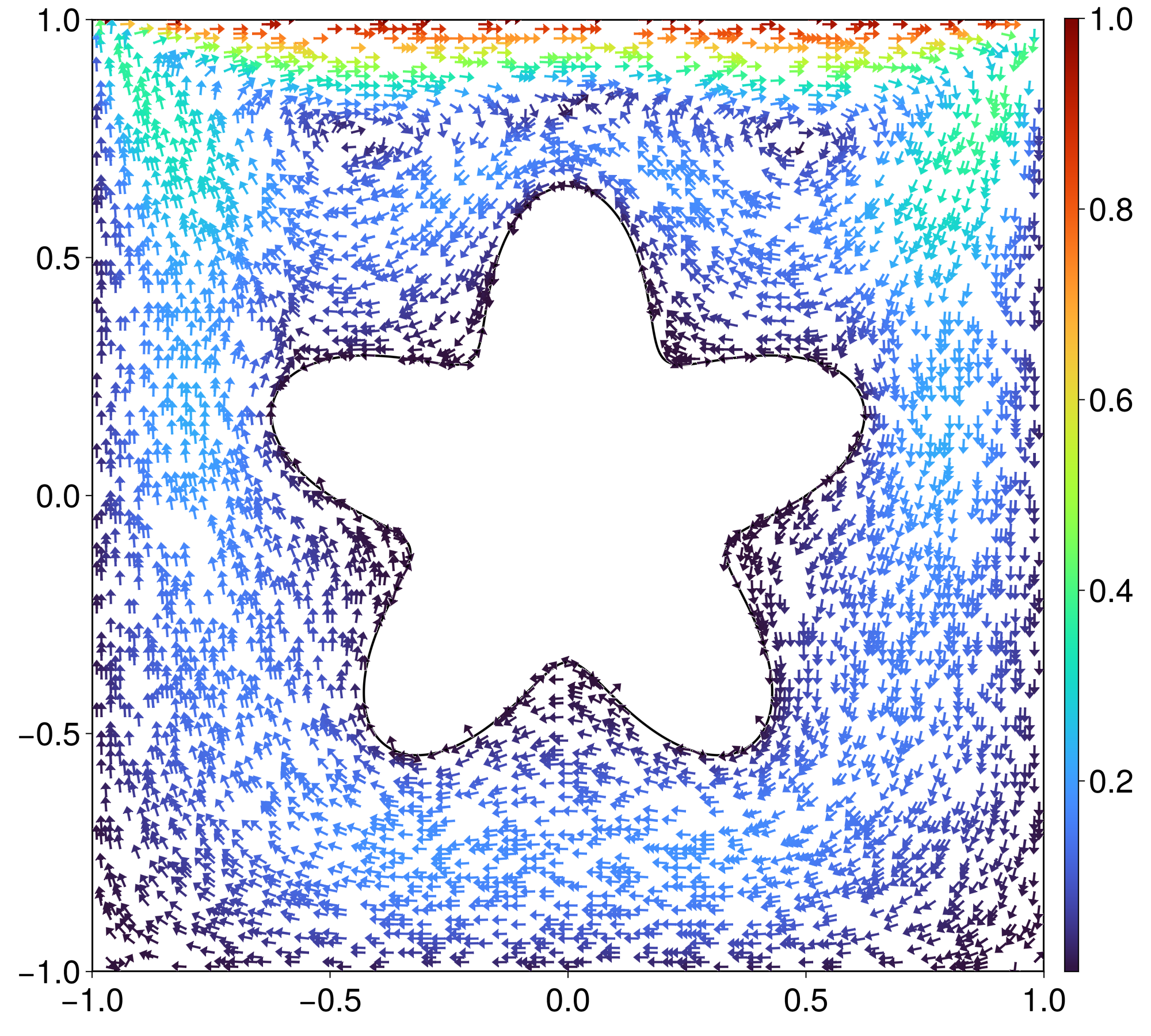}
        \caption{$Re = 10$,\quad $t = 10s$}
   \end{subfigure}
   \hfill
   \begin{subfigure}{0.45\textwidth}
       \includegraphics[width=\textwidth]{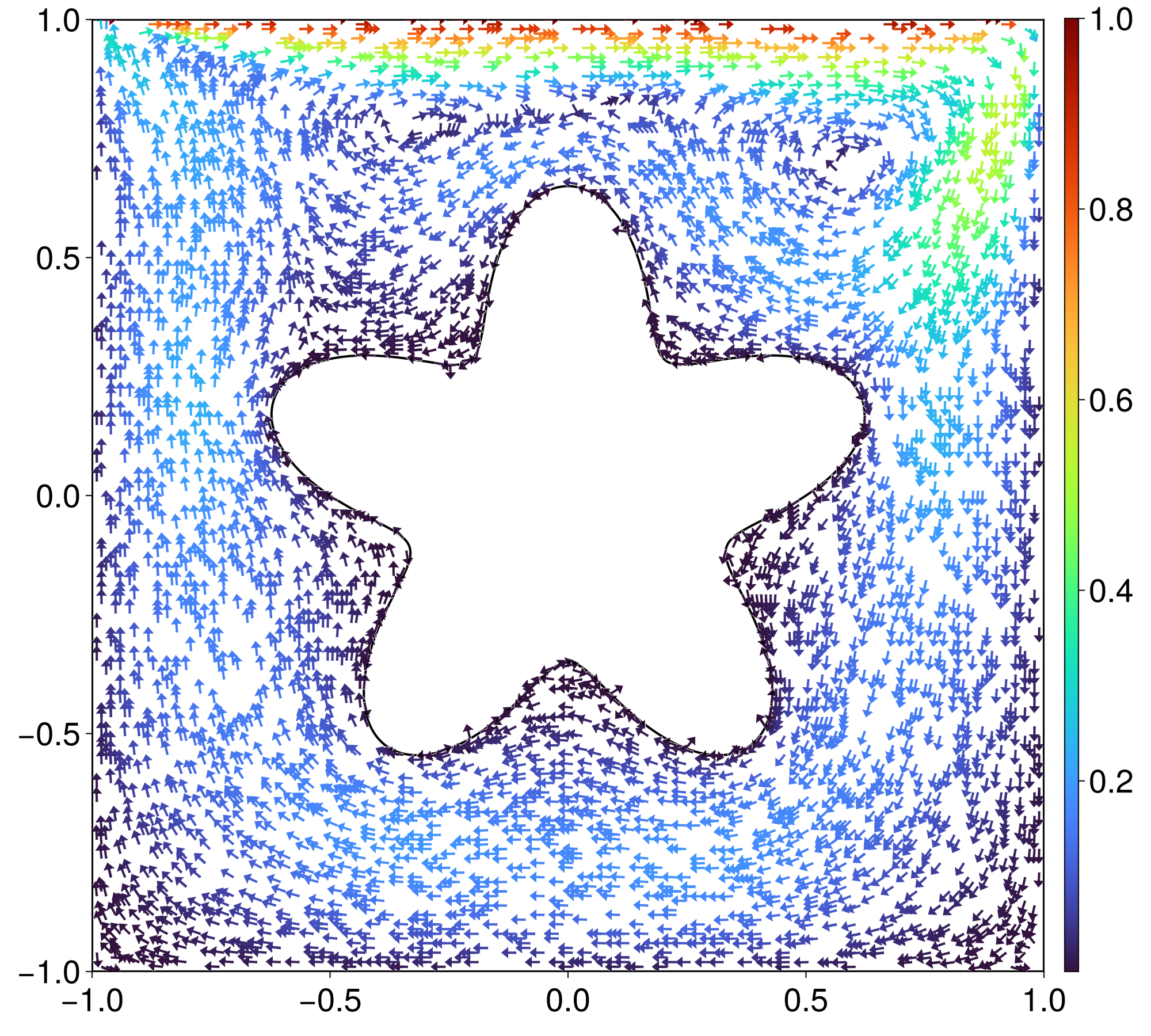}
        \caption{$Re = 100$,\quad $t = 10s$}
   \end{subfigure}
   \hfill
   \begin{subfigure}{0.45\textwidth}
       \includegraphics[width=\textwidth]{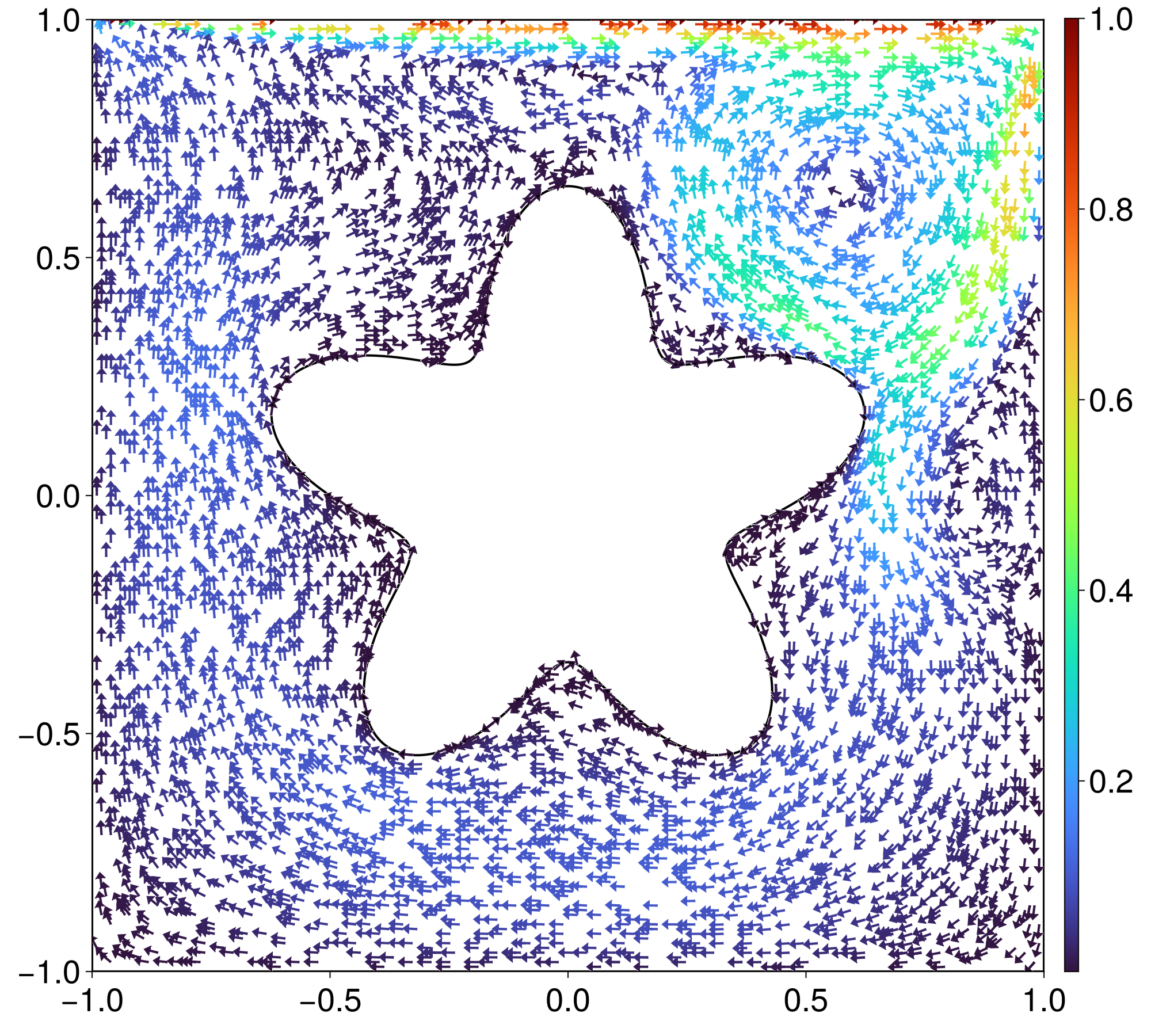}
        \caption{$Re = 1000$,\quad $t = 10s$}
   \end{subfigure}
    \caption{\textit{Velocity field for the driven cavity test with $N_x=100$}.}
    \label{fig:driven_cavity}
\end{figure}

In this section, we investigate a driven cavity flow around a stationary flower-shaped obstacle \cite{COCO2020109623}. No-slip conditions, $\mathbf{u} = (0, 0)$, are applied on the surface of the cavity and on all boundaries of the computational domain $(-1,1)^2$, except for the upper boundary, which represents a lid moving horizontally at a constant speed, $\mathbf{u} = (1, 0)$. External forces, including gravity, are ignored ($\mathbf{f}=0$). Fig.~\ref{fig:driven_cavity} depicts the velocity field at time $t = 10 s$ using 
$N_x=100$ for four Reynolds numbers: 
$Re=1, 10, 100, 1000$.  
\par The presence of the flower-shaped object inhibits the development of anticlockwise vortices in the lower corners, which are typically observed in the classical driven cavity configuration without obstacles. Furthermore, the large central clockwise vortex commonly found in the standard case is divided into two distinct vortical structures at low Reynolds numbers ($Re=1$ and $Re=10$). This behavior is primarily induced by the reentrant geometry of the flower-shaped object. The centers of these vortices are located approximately at $x=-0.5$ and $x=0.5$. As the Reynolds number increases to $Re=100$, the upper-left vortex becomes more asymmetric and shifts slightly to the right. At $Re=1000$, this vortex disappears, while a new vortex begins to form in the central-right region of the cavity. The resulting flow structures are in good agreement with those reported in \cite{COCO2020109623}.
\subsection{Flow around a rotating flower-shaped object}
\begin{figure}
    \centering
   \begin{subfigure}{0.32\textwidth}
       \includegraphics[width=\textwidth]{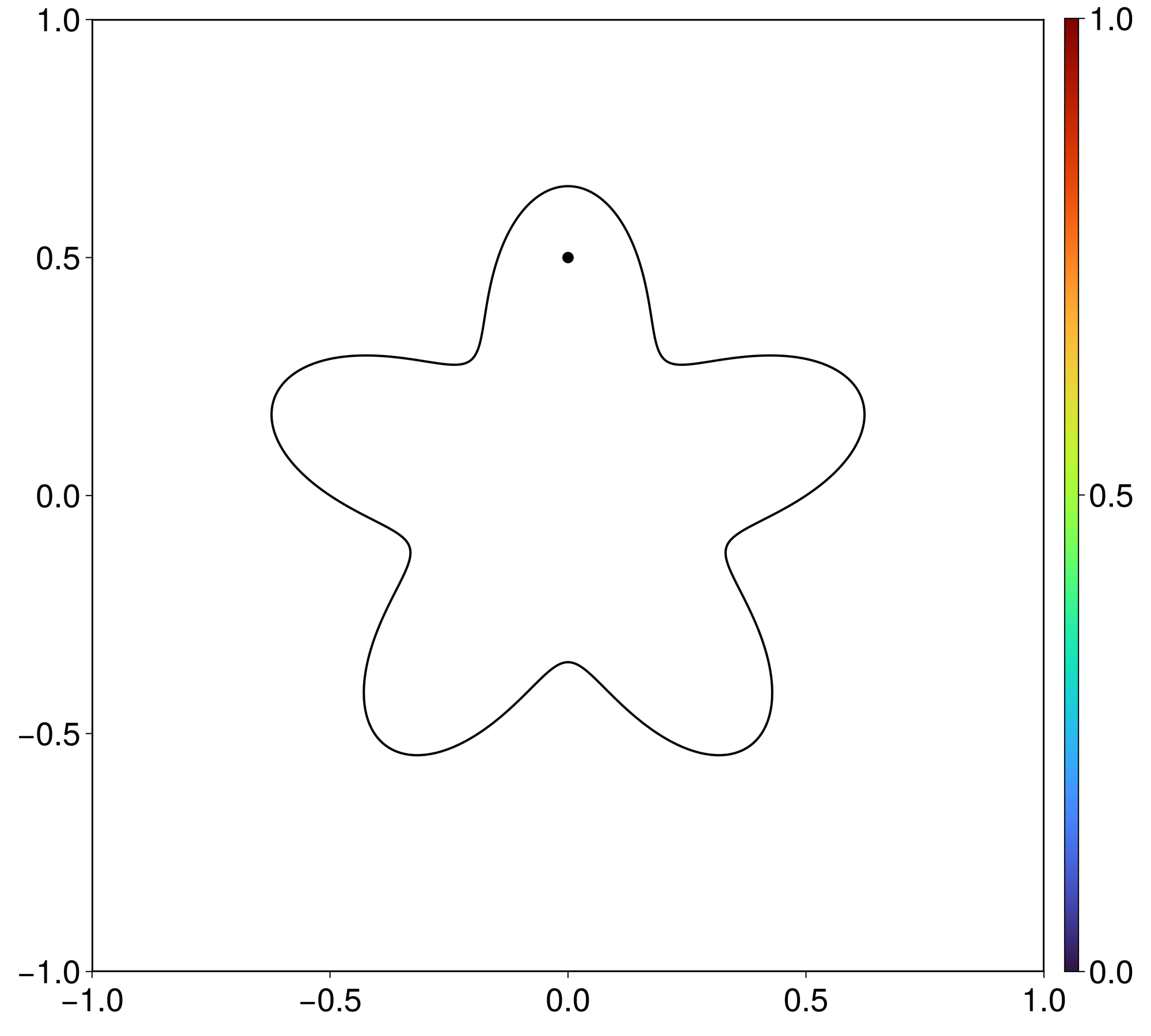}
       \caption{$Re = 100$,\quad $t = 0s$}
   \end{subfigure}
    \hfill 
   \begin{subfigure}{0.32\textwidth}
       \includegraphics[width=\textwidth]{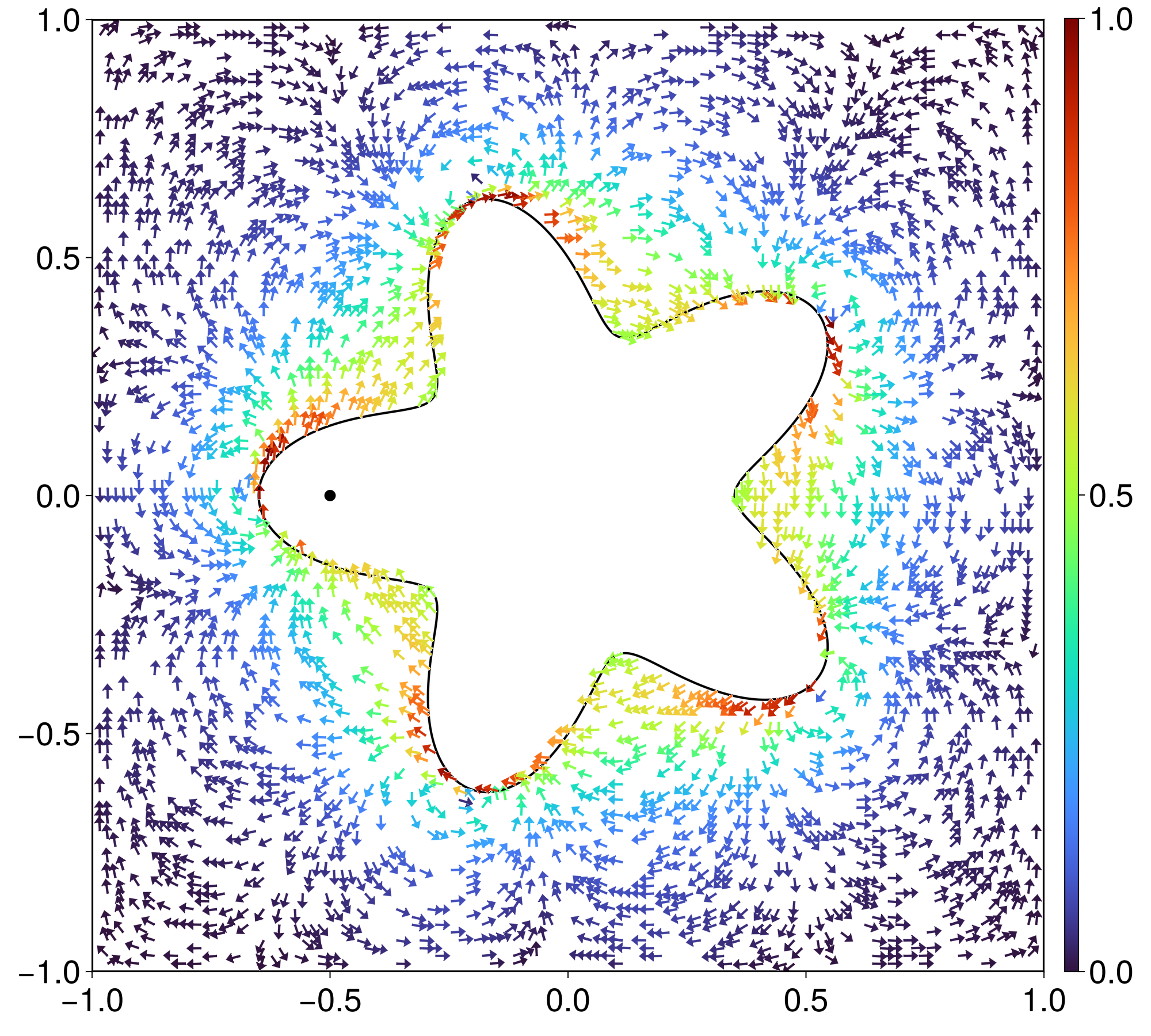}
        \caption{$Re = 100$,\quad $t = 1.25s$}
   \end{subfigure}
   \hfill
   \begin{subfigure}{0.32\textwidth}
       \includegraphics[width=\textwidth]{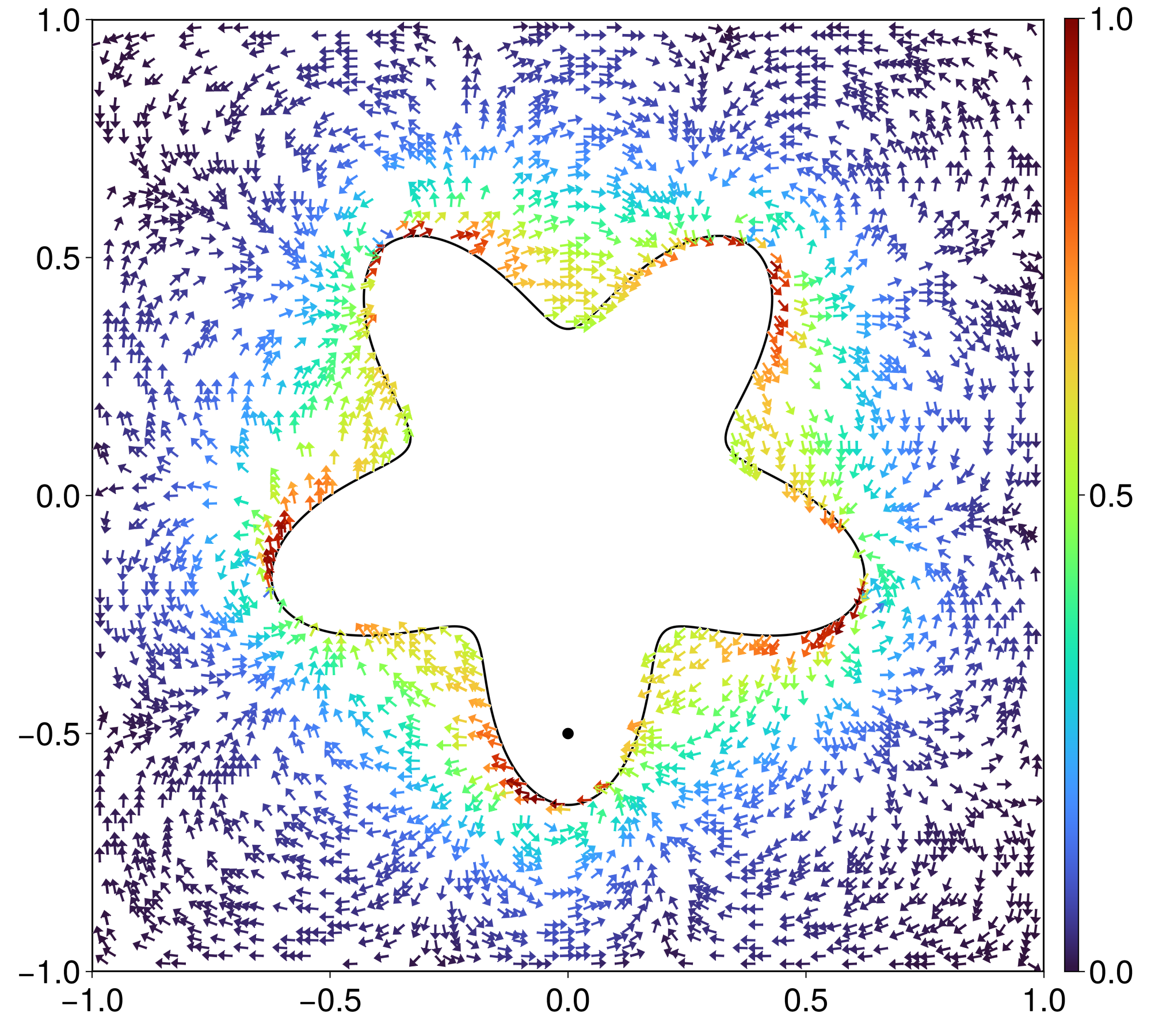}
        \caption{$Re = 100$,\quad $t = 2.50s$}
   \end{subfigure}
   \hfill
   \begin{subfigure}{0.32\textwidth}
       \includegraphics[width=\textwidth]{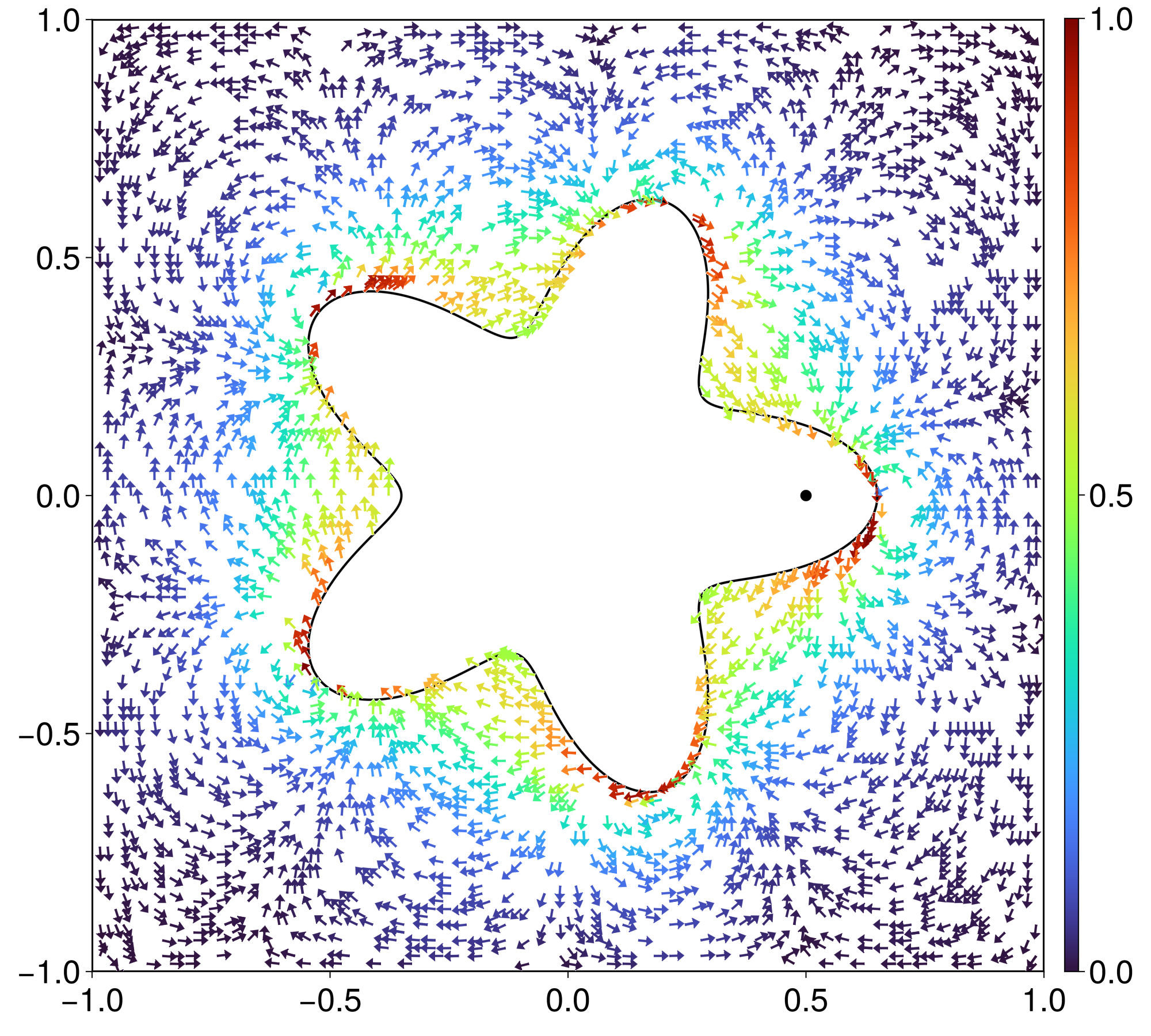}
        \caption{$Re = 100$,\quad $t = 3.75s$}
   \end{subfigure}
   \hfill
   \begin{subfigure}{0.32\textwidth}
       \includegraphics[width=\textwidth]{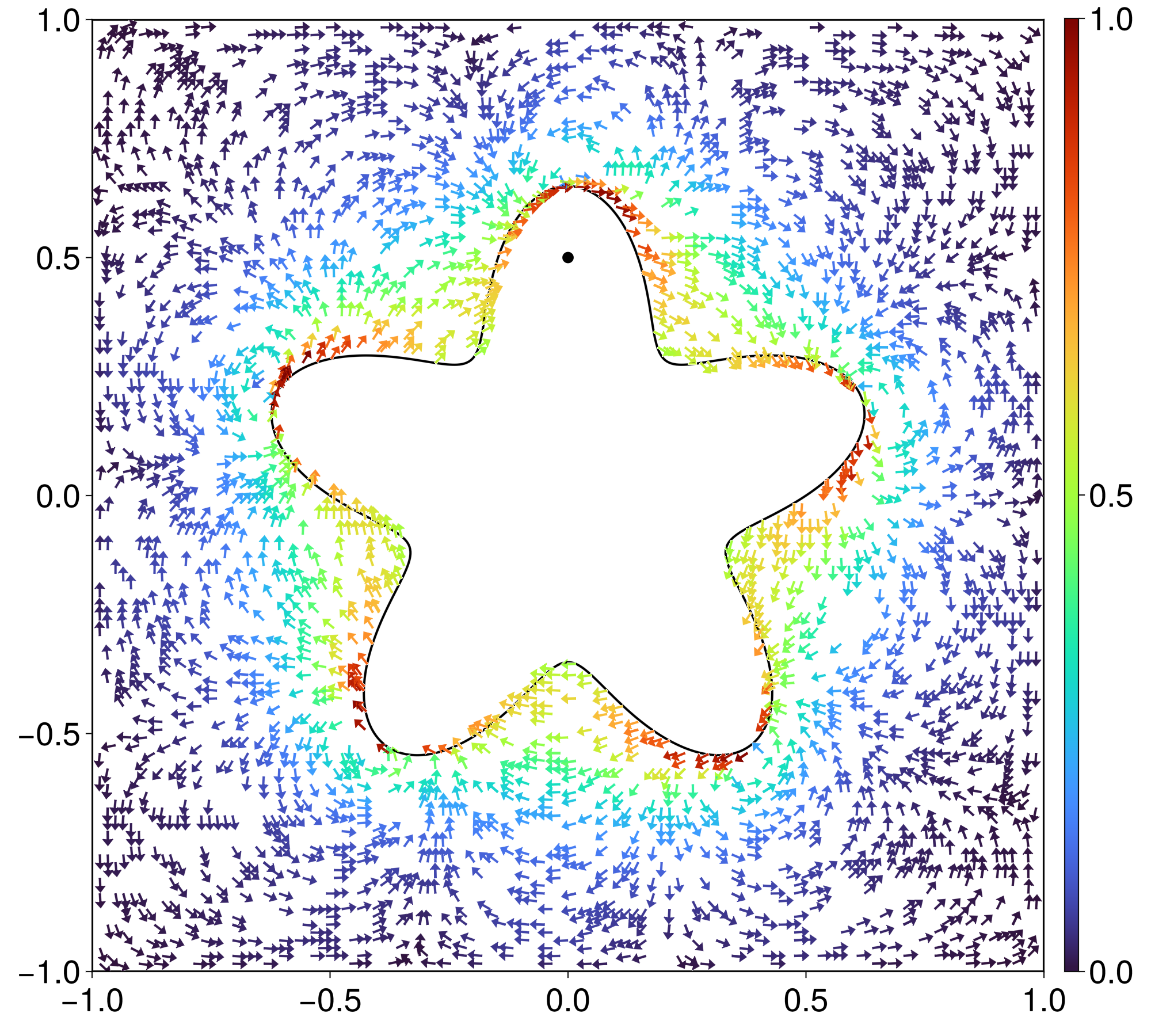}
       \caption{$Re = 100$,\quad $t = 5.00s$}
   \end{subfigure}
    \hfill 
   \begin{subfigure}{0.32\textwidth}
       \includegraphics[width=\textwidth]{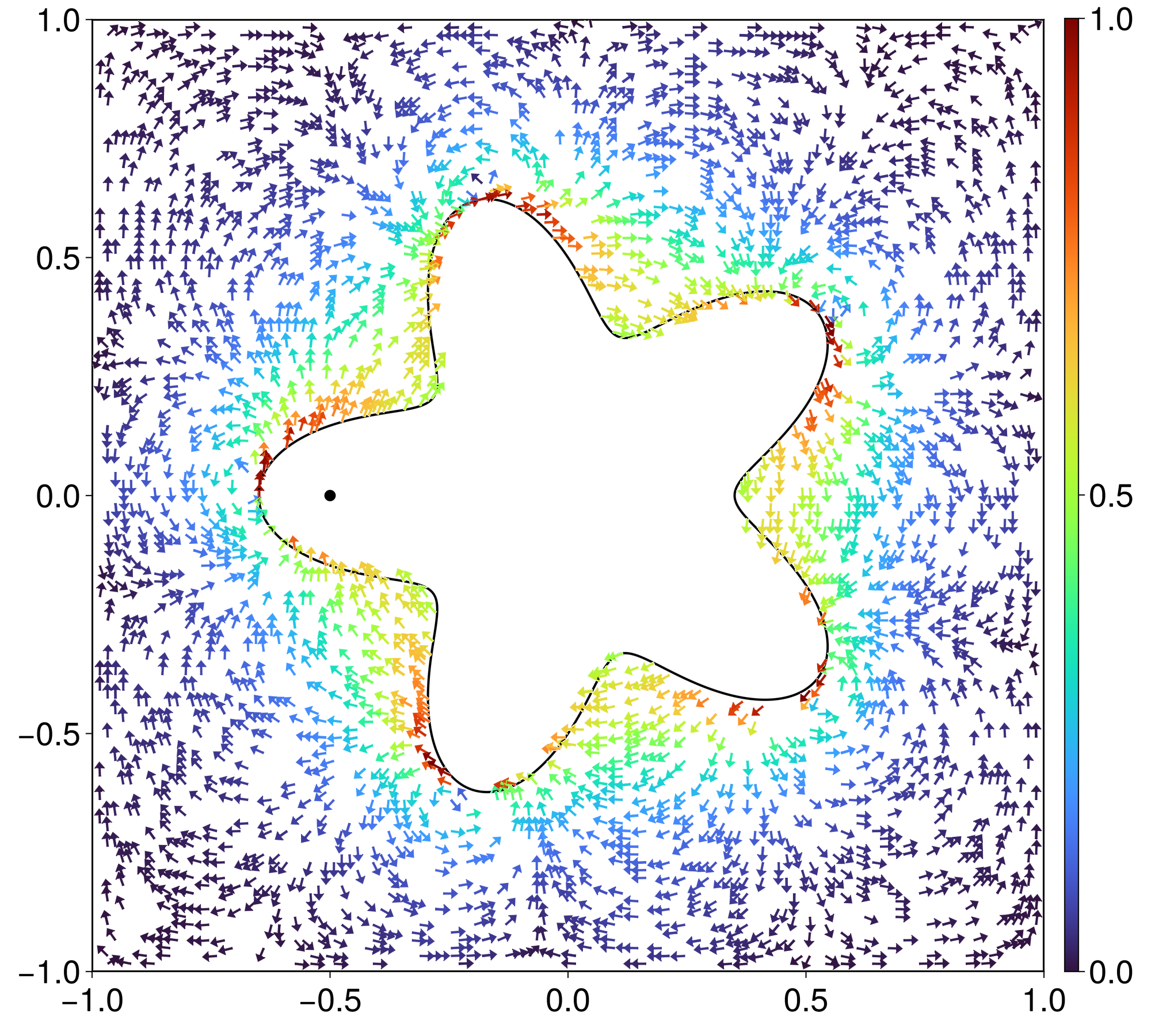}
        \caption{$Re = 100$,\quad $t = 6.25s$}
   \end{subfigure}
   \hfill
   \begin{subfigure}{0.32\textwidth}
       \includegraphics[width=\textwidth]{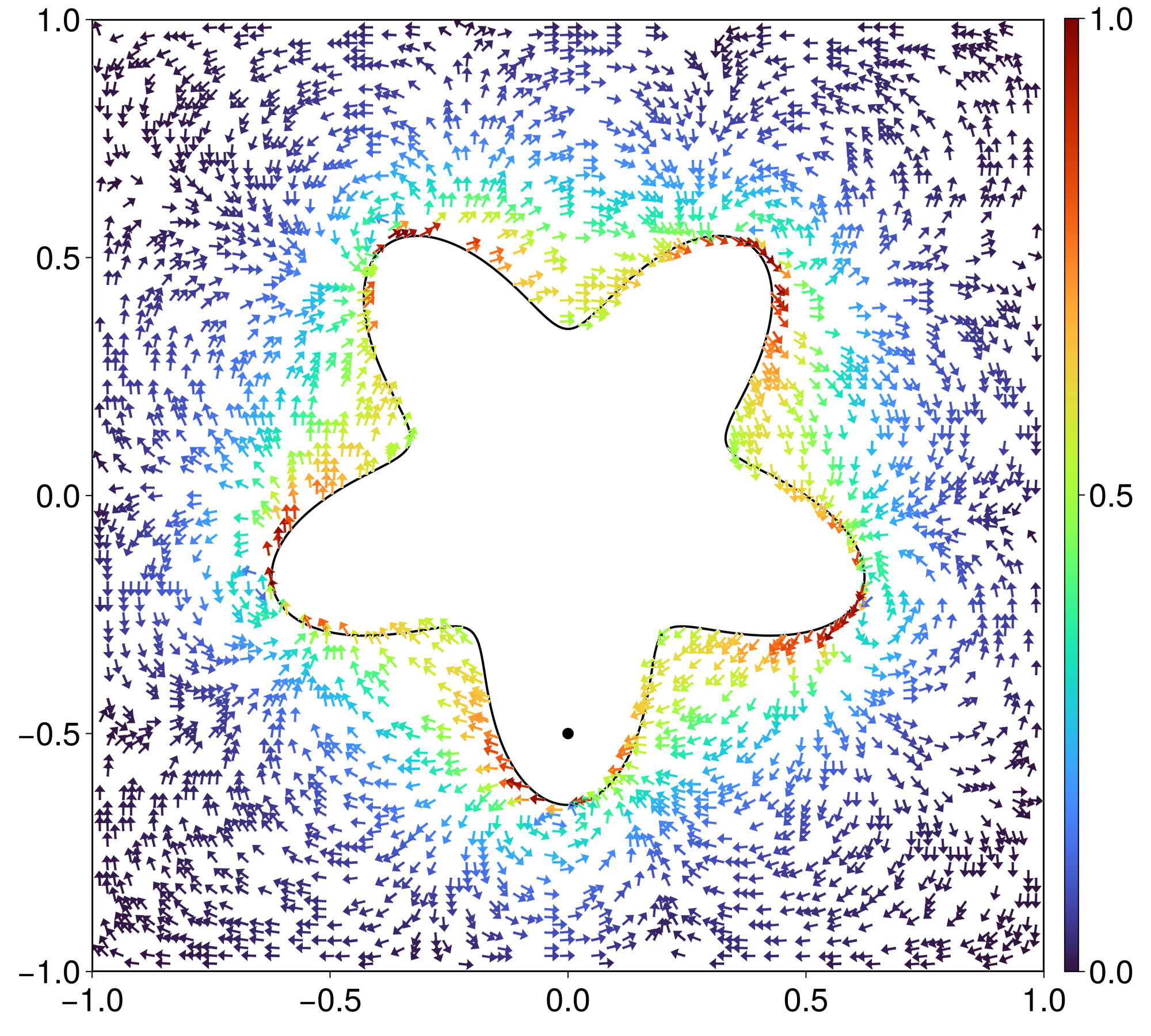}
        \caption{$Re = 100$,\quad $t = 7.50s$}
   \end{subfigure}
   \hfill
   \begin{subfigure}{0.32\textwidth}
       \includegraphics[width=\textwidth]{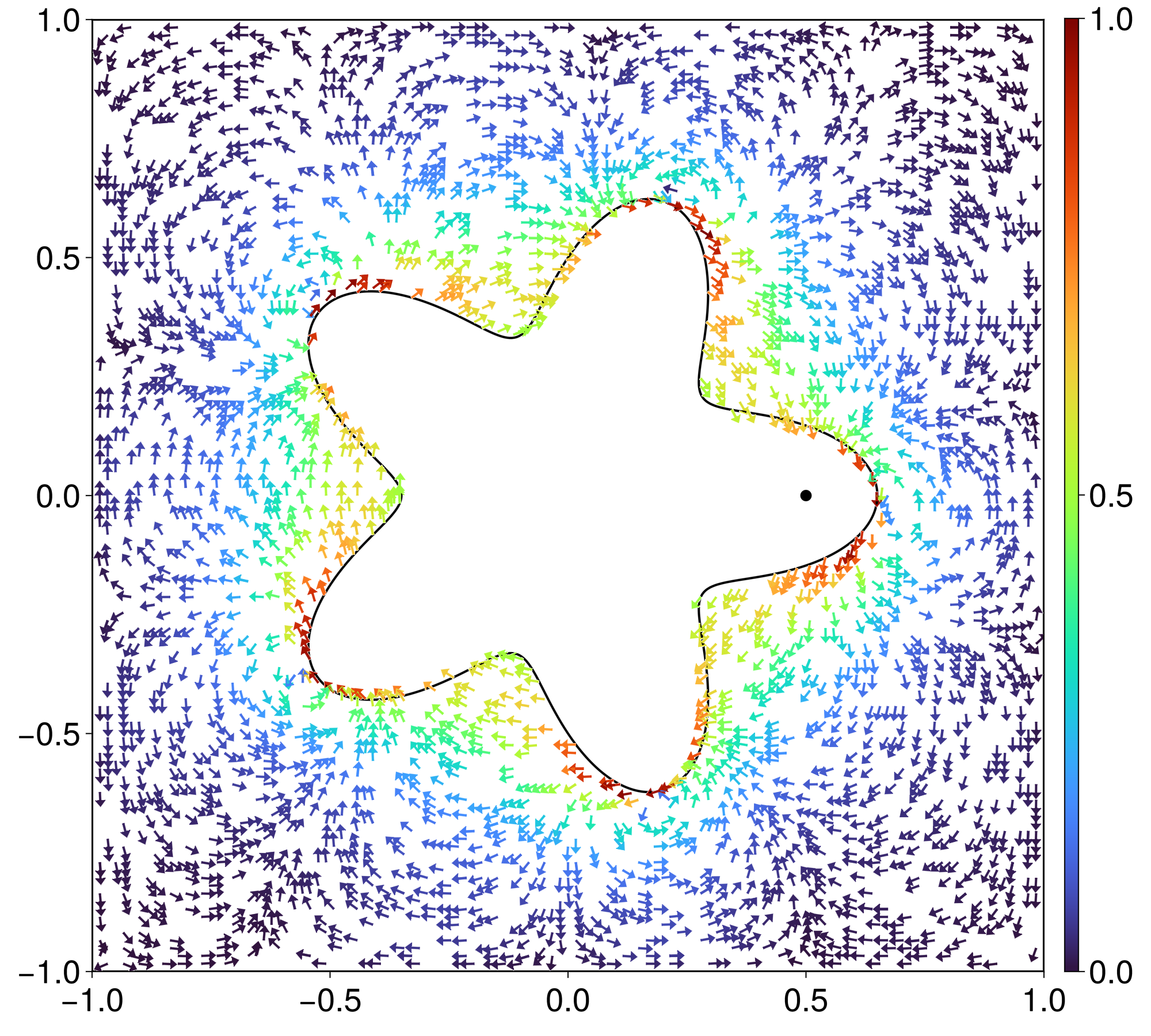}
        \caption{$Re = 100$,\quad $t = 8.75s$}
   \end{subfigure}
   \hfill
   \begin{subfigure}{0.32\textwidth}
       \includegraphics[width=\textwidth]{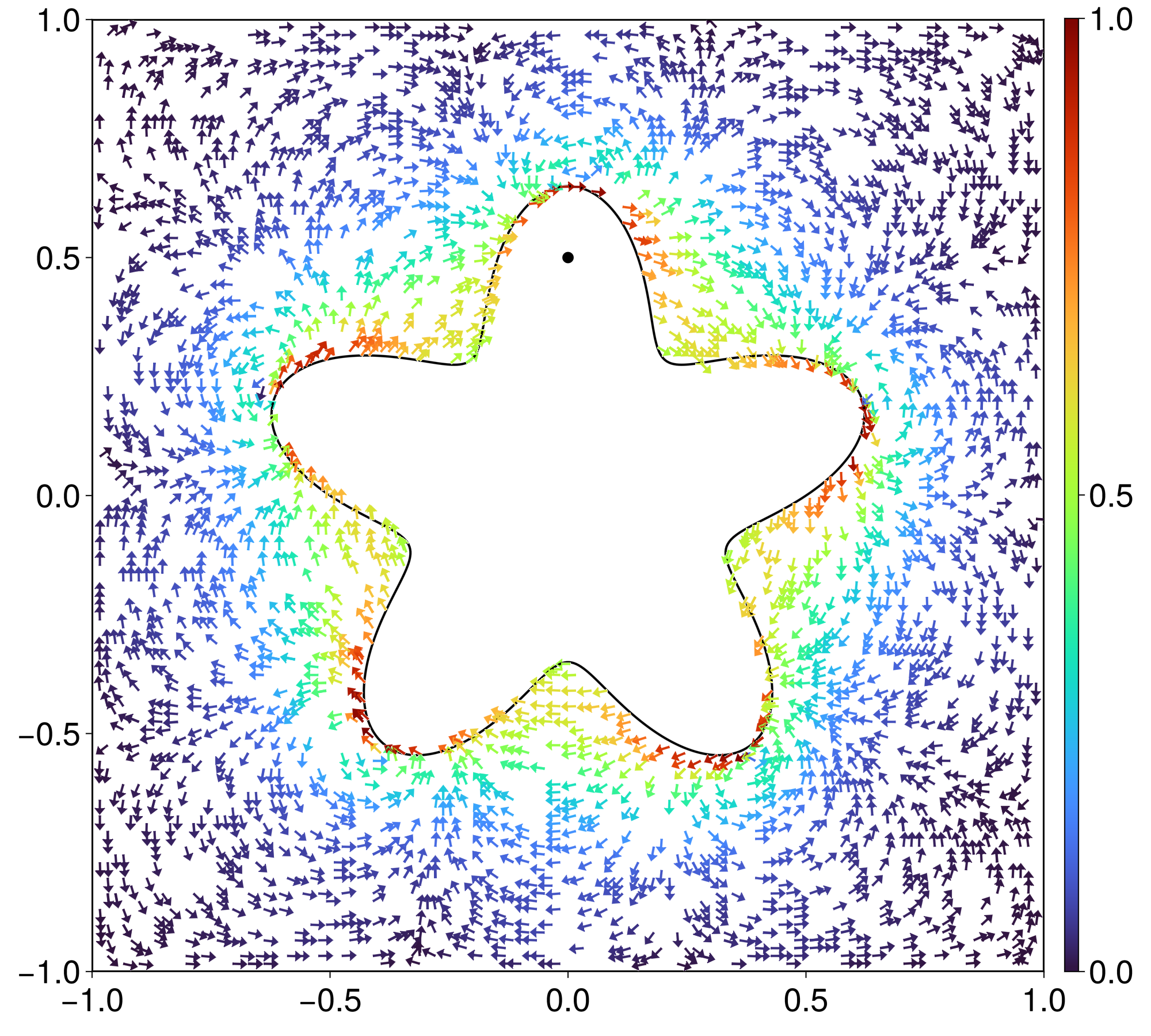}
        \caption{$Re = 100$,\quad $t = 10.0s$}
   \end{subfigure}
    \caption{Velocity field at different time instants for the rotating flower-shaped object test with $N_x=60$ at $Re=100$. Vortices are formed at the tips of the petals at all selected time instants.}
    \label{fig:rot_flower}
\end{figure}

\begin{figure}
    \centering
   \begin{subfigure}{0.32\textwidth}
       \includegraphics[width=\textwidth]{figures/rot_flower_quiver_t_0.png}
       \caption{$Re = 1$,\quad $t = 0s$}
   \end{subfigure}
    \hfill 
   \begin{subfigure}{0.32\textwidth}
       \includegraphics[width=\textwidth]{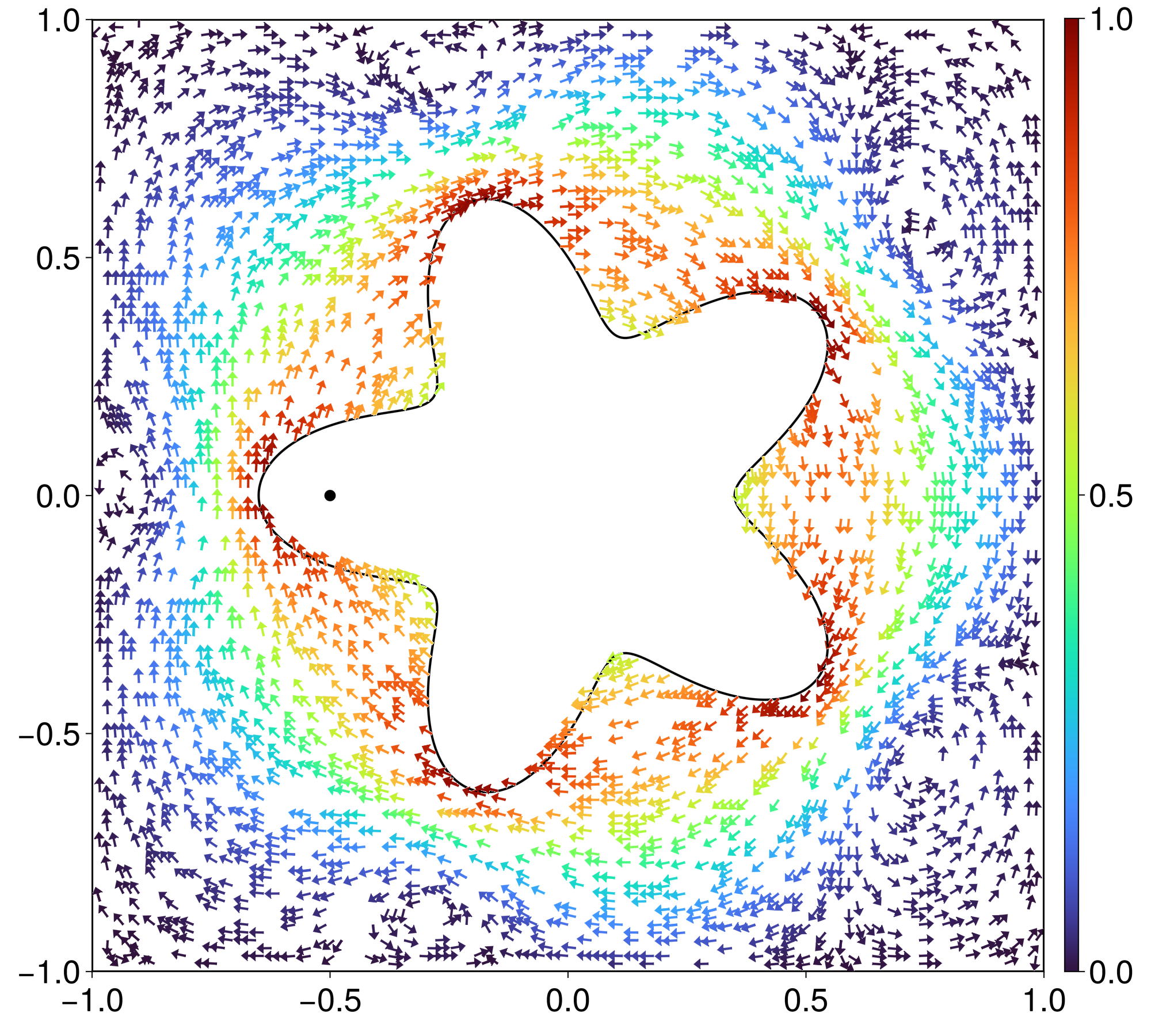}
        \caption{$Re = 1$,\quad $t = 1.25s$}
   \end{subfigure}
   \hfill
   \begin{subfigure}{0.32\textwidth}
       \includegraphics[width=\textwidth]{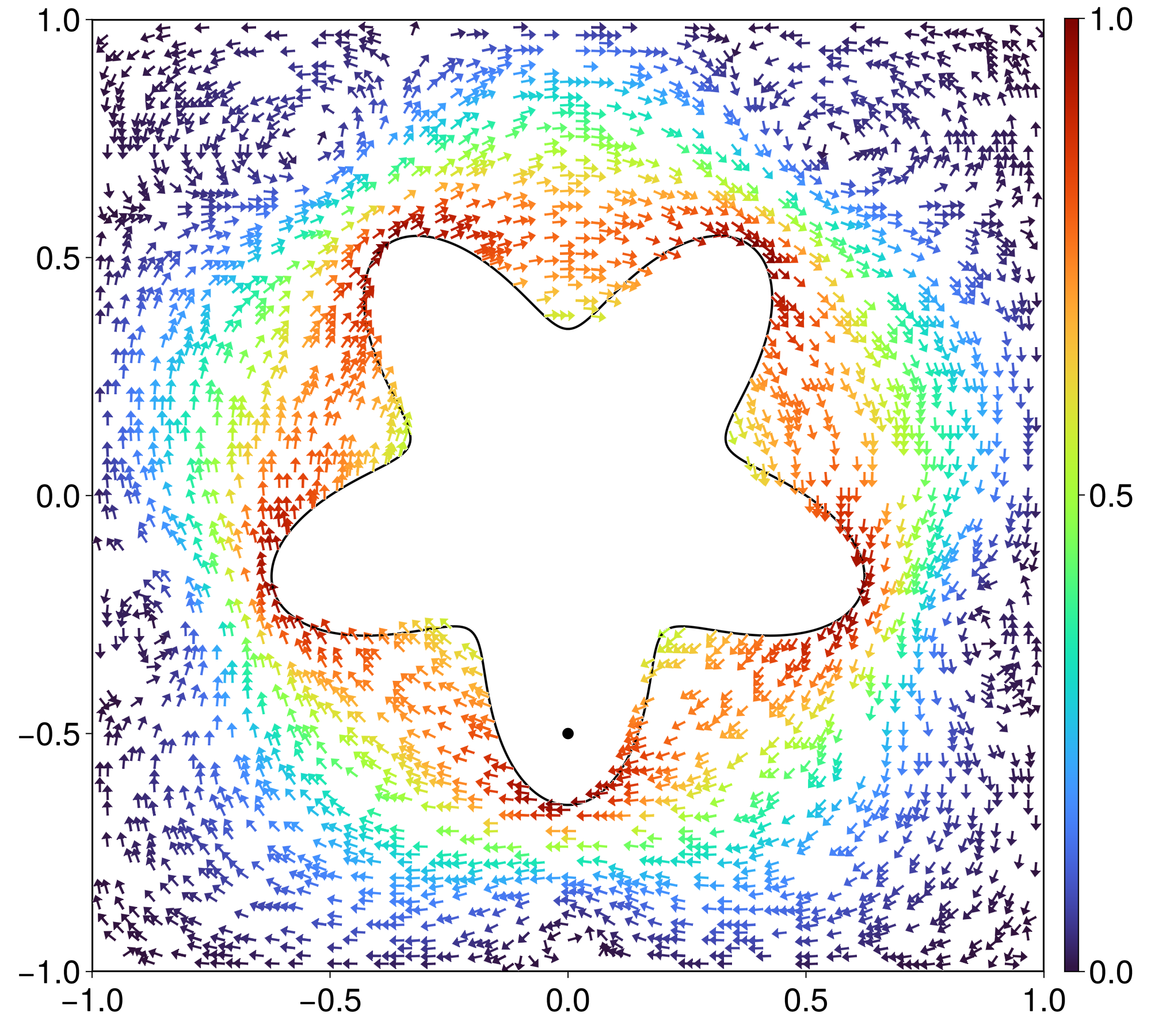}
        \caption{$Re = 1$,\quad $t = 2.50s$}
   \end{subfigure}
   \hfill
   \begin{subfigure}{0.32\textwidth}
       \includegraphics[width=\textwidth]{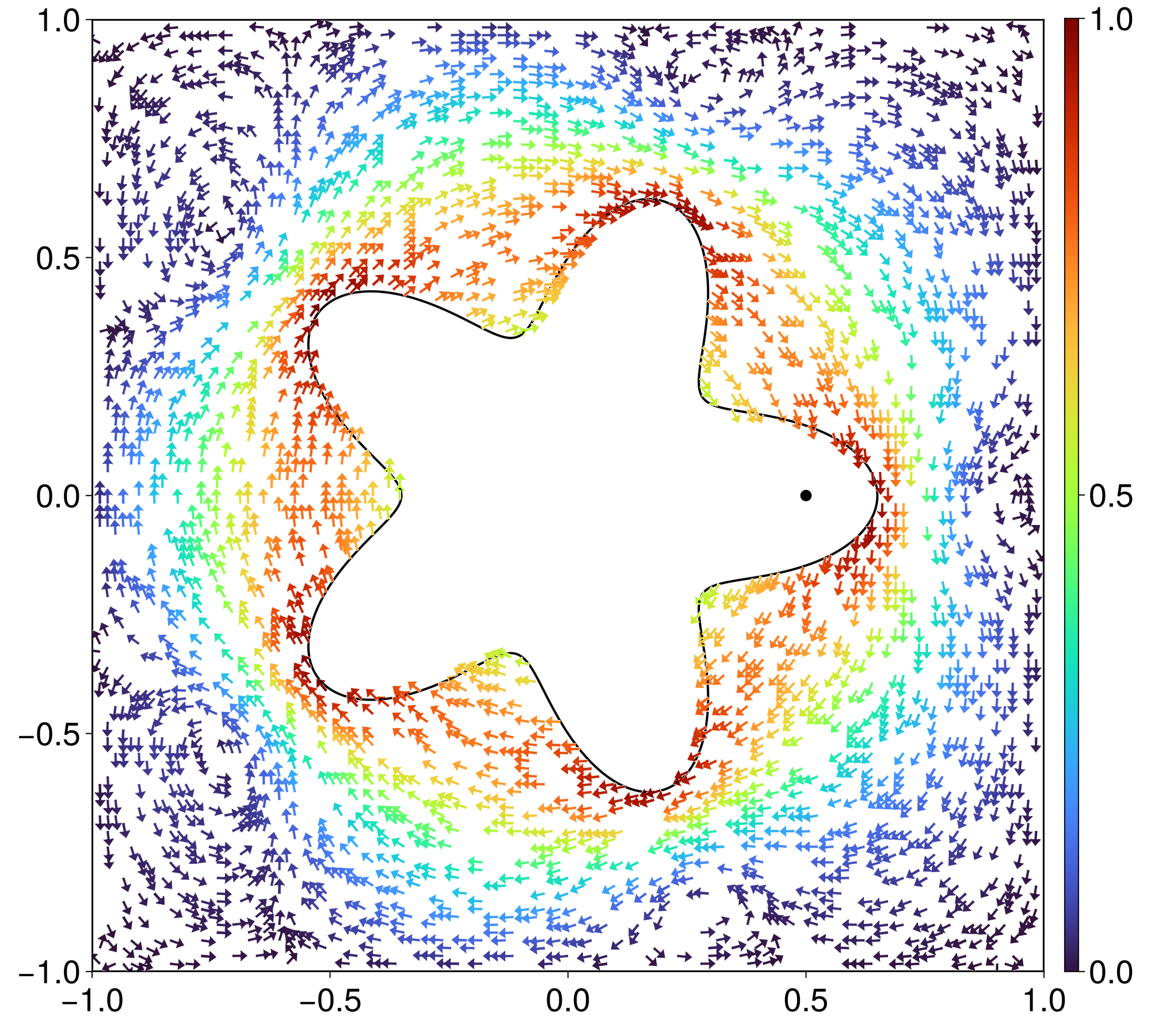}
        \caption{$Re = 1$,\quad $t = 3.75s$}
   \end{subfigure}
   \begin{subfigure}{0.32\textwidth}
       \includegraphics[width=\textwidth]{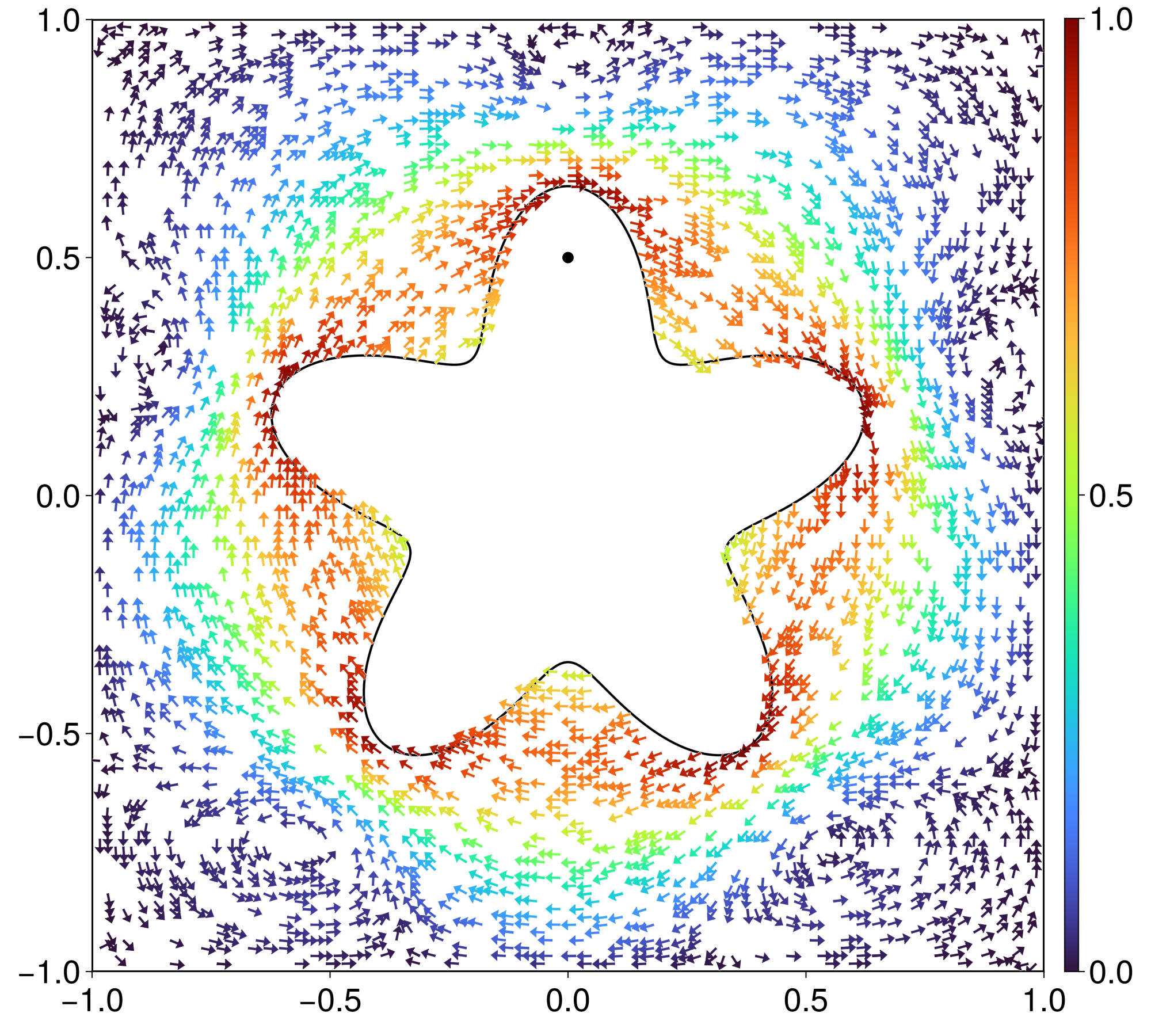}
       \caption{$Re = 1$,\quad $t = 5.00s$}
   \end{subfigure}
    \hfill 
   \begin{subfigure}{0.32\textwidth}
       \includegraphics[width=\textwidth]{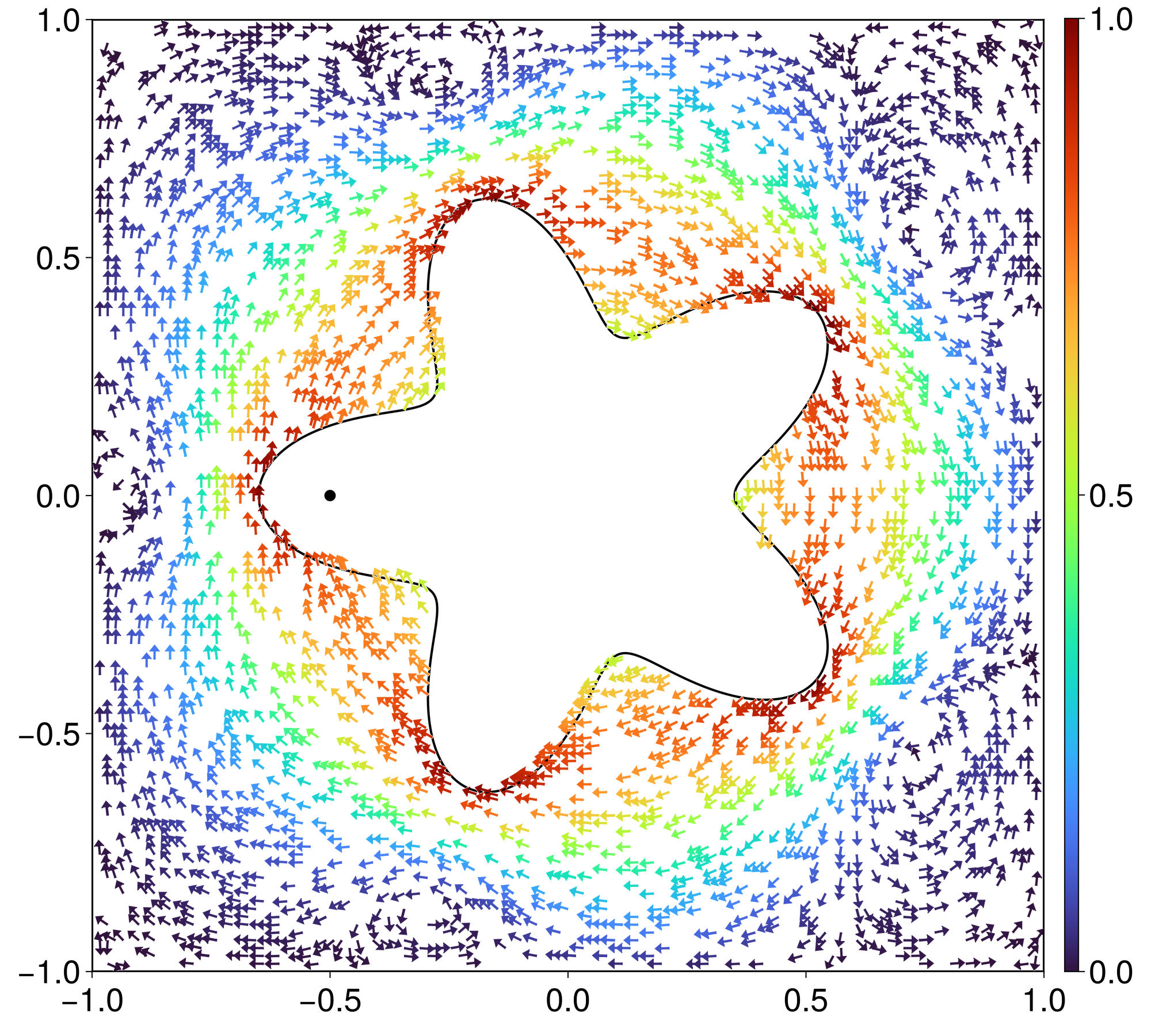}
        \caption{$Re = 1$,\quad $t = 6.25s$}
   \end{subfigure}
   \hfill
   \begin{subfigure}{0.32\textwidth}
       \includegraphics[width=\textwidth]{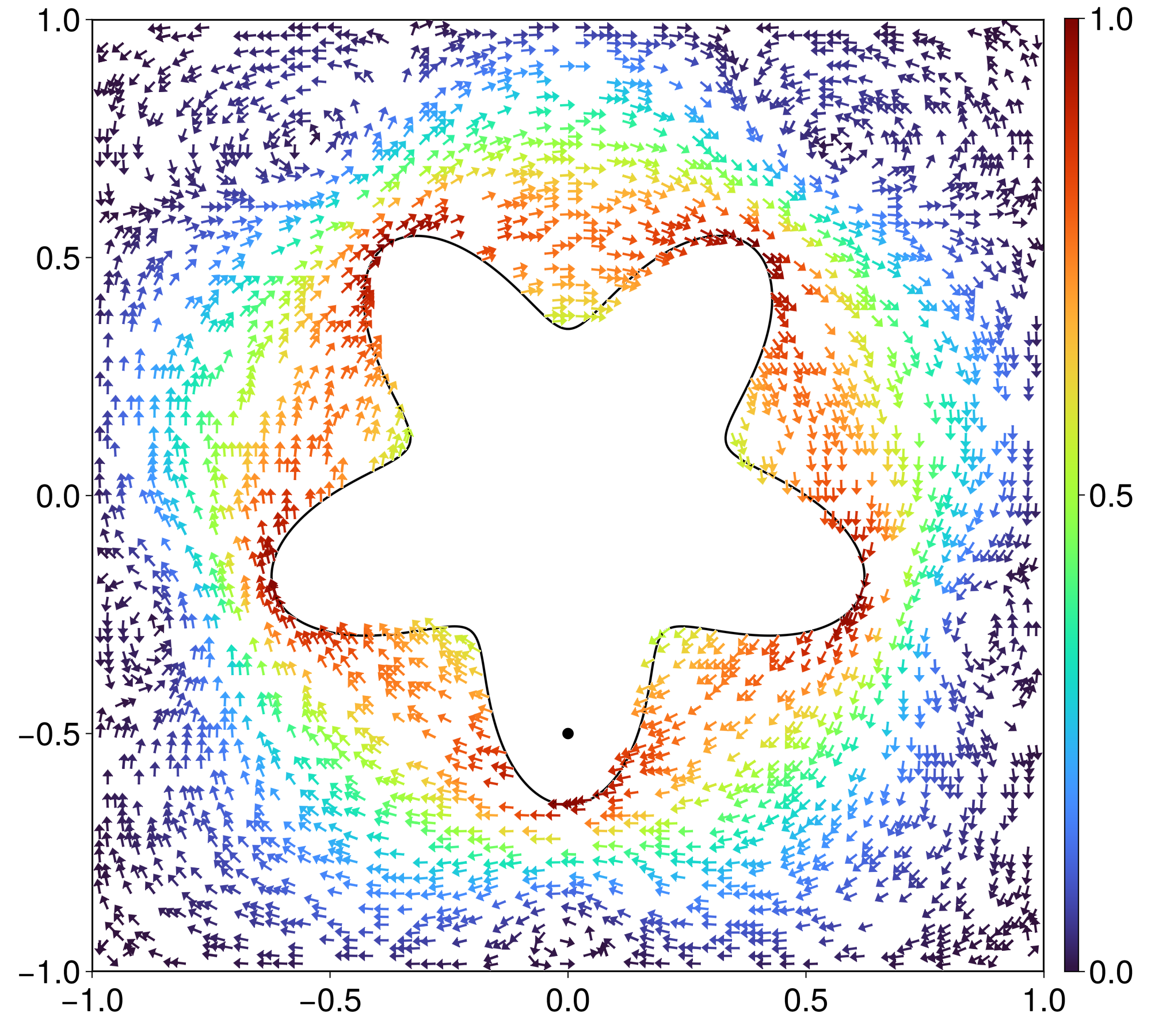}
        \caption{$Re = 1$,\quad $t = 7.50s$}
   \end{subfigure}
   \hfill
   \begin{subfigure}{0.32\textwidth}
       \includegraphics[width=\textwidth]{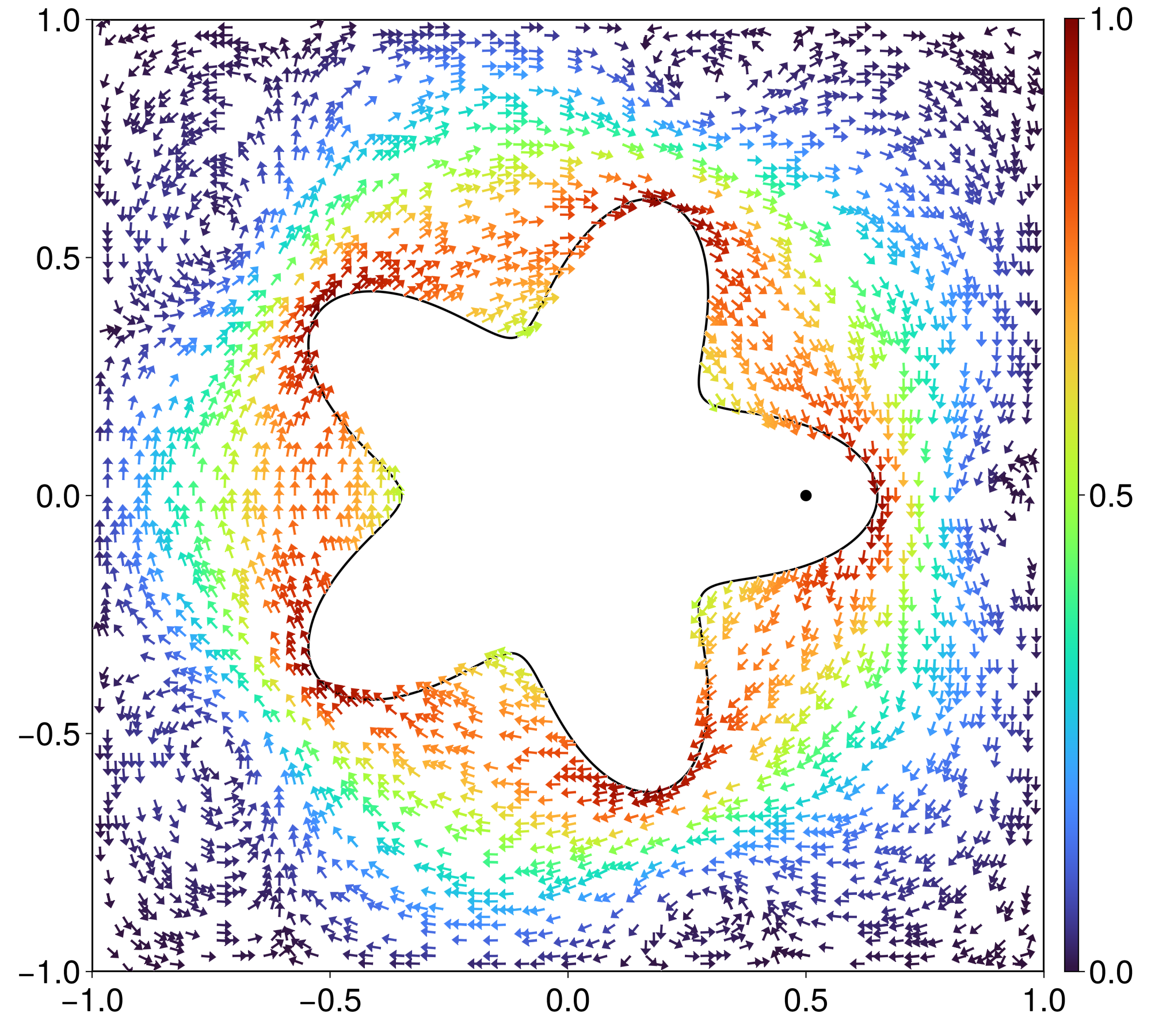}
        \caption{$Re = 1$,\quad $t = 8.75s$}
   \end{subfigure}
   \hfill
   \begin{subfigure}{0.32\textwidth}
       \includegraphics[width=\textwidth]{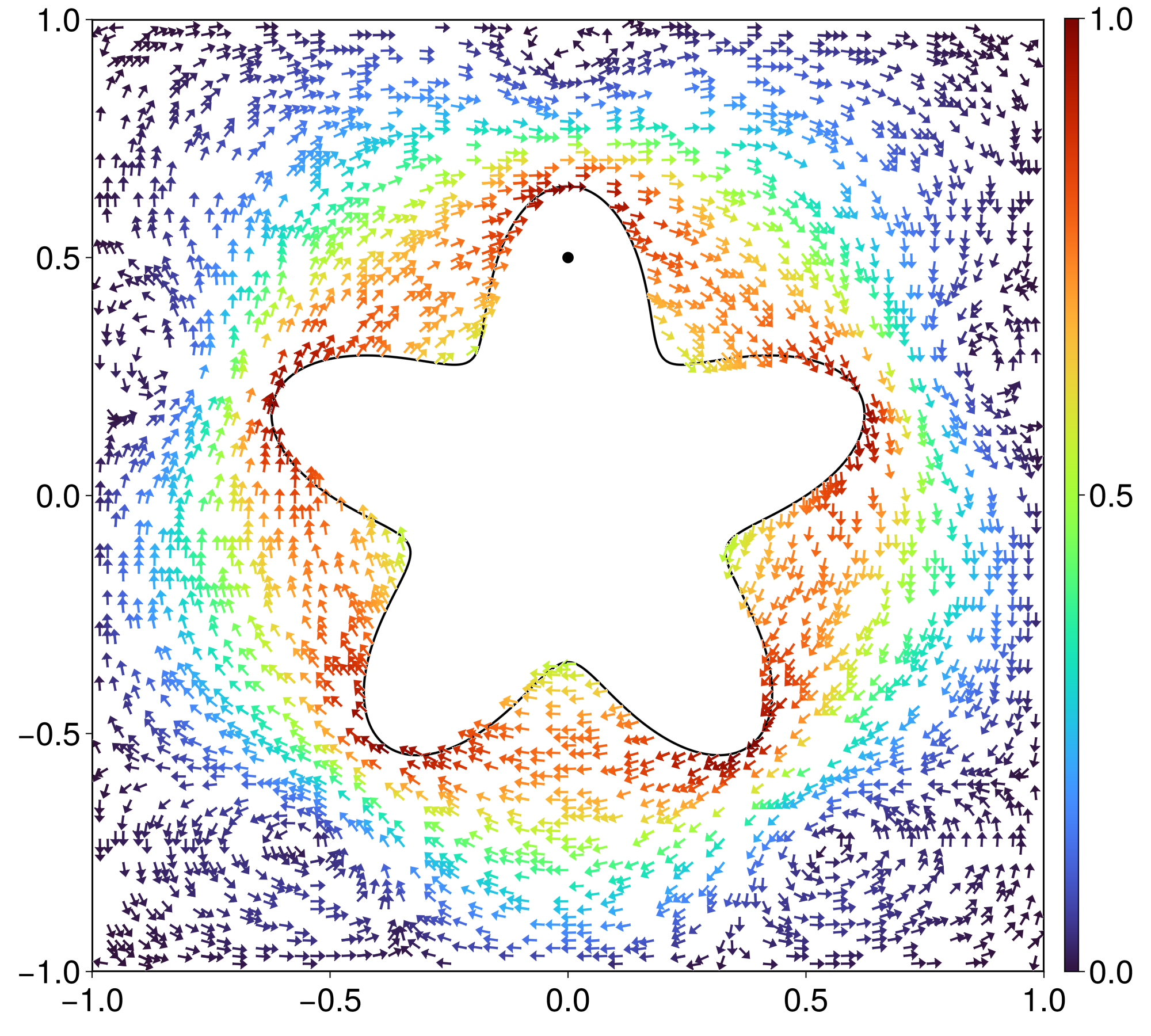}
        \caption{$Re = 1$,\quad $t = 10.0s$}
   \end{subfigure}
    \caption{Velocity field at different time instants for the rotating flower-shaped object test with $N_x=60$ at $Re=1$. Vortices are observed in the wake of the petals farther from the tips at all time instants.}
    \label{fig:rot_flower_re_1}
\end{figure}

In this test, a flower-shaped object rotates anticlockwise about the origin with a constant angular velocity of $\omega = 2 \pi /5$ \cite{COCO2020109623}. The fluid is initially at rest, and no-slip boundary conditions are applied on the walls of the computational domain $(\mathbf{u} =0)$ as well as on the boundary of the object $(\mathbf{u}(x,y) = (0,0,\omega) \times (x,y,0))$ throughout the simulation. Fig.~\ref{fig:rot_flower} and Fig.~\ref{fig:rot_flower_re_1} illustrate the velocity field at nine different time steps for a grid resolution of $N_x=60$, at $Re =100$ and $Re = 1$, respectively.
The moving domain is defined by the time-dependent level set function $\phi(x,y,t)=\phi_0(\tilde{x}(t),\tilde{y}(t))$, where $\phi_0(x,y)$ represents the flower at time $t=0$, defined in~\eqref{phi:ellipse}, and $(\tilde{x}(t),\tilde{y}(t))$ are given by~\eqref{eq:rotxy}.

\par When $Re=1$, vortices are formed in the wake of the petals farther from the tips at all time instants. This behavior arises because, at low Reynolds numbers, diffusion dominates the flow. In contrast, when $Re=100$, vortices are observed at the tips of the petals of the flower-shaped object at all selected time instants. At this Reynolds number, the flow is no longer diffusion-dominated; consequently, vortices remain closer to the petal tips. 
\par Compared with the results reported in \cite{COCO2020109623} for $Re=100$, no corner vortices are observed in the present simulations. Moreover, the flow exhibits periodic behavior and does not attain a steady state.

\subsection{Flow past oscillating bubbles}
\begin{figure}
    \centering
   \begin{subfigure}{0.245\textwidth}
       \includegraphics[width=\textwidth]{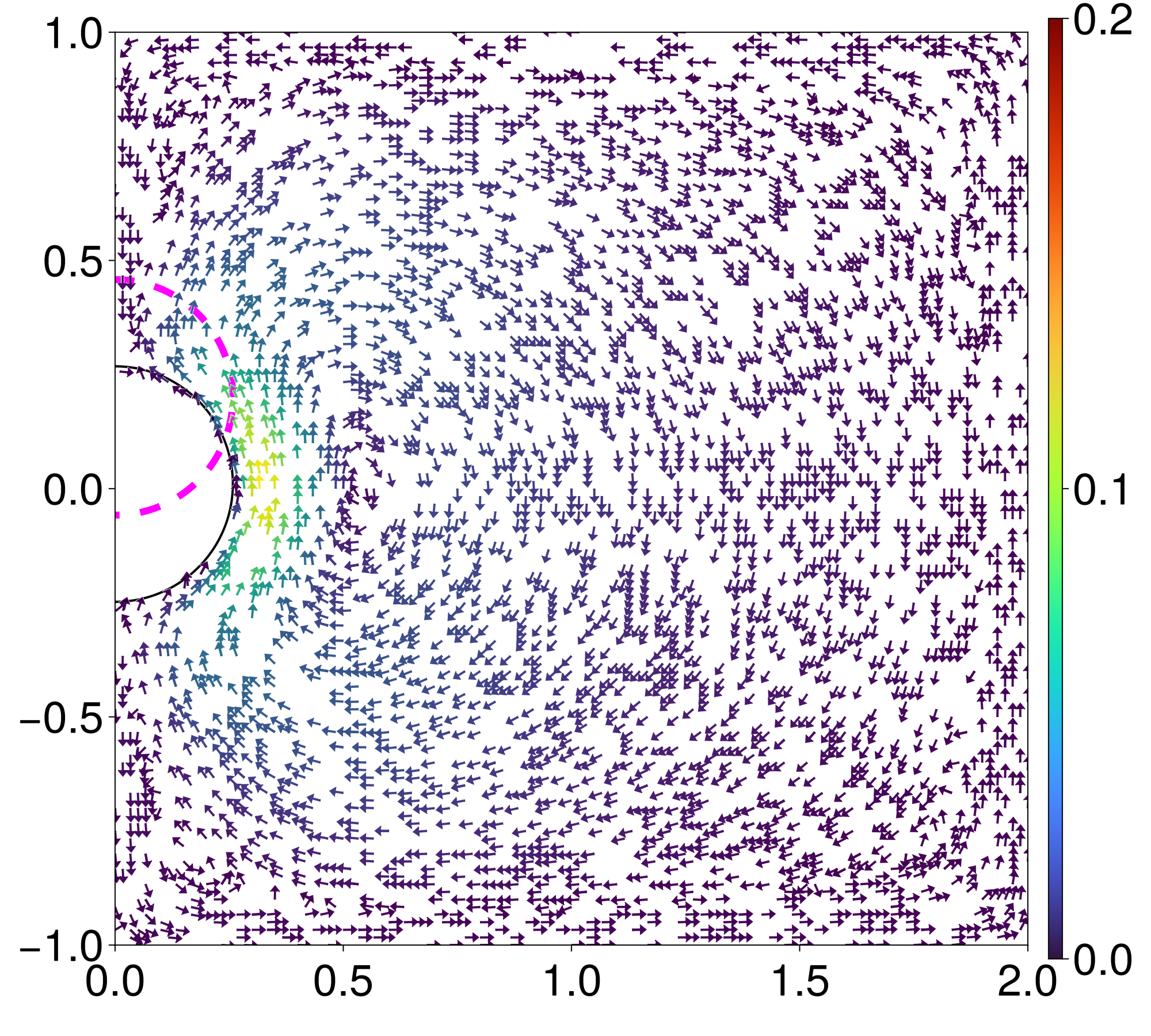}
       \caption{$\nu \cdot t = 0.25$}
   \end{subfigure}
    \hfill 
   \begin{subfigure}{0.245\textwidth}
       \includegraphics[width=\textwidth]{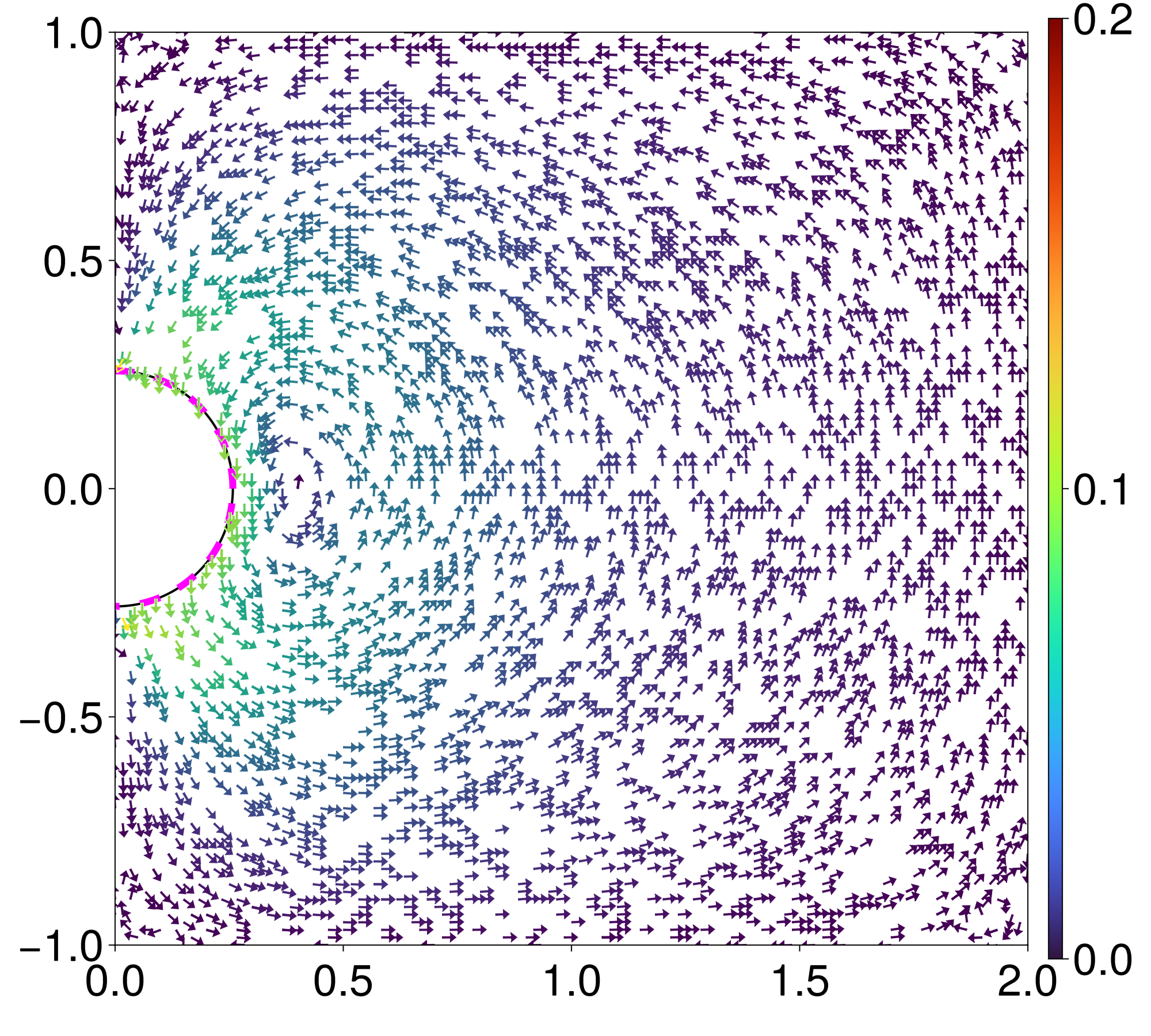}
        \caption{$\nu \cdot t = 0.5$}
   \end{subfigure}
   \hfill
   \begin{subfigure}{0.245\textwidth}
       \includegraphics[width=\textwidth]{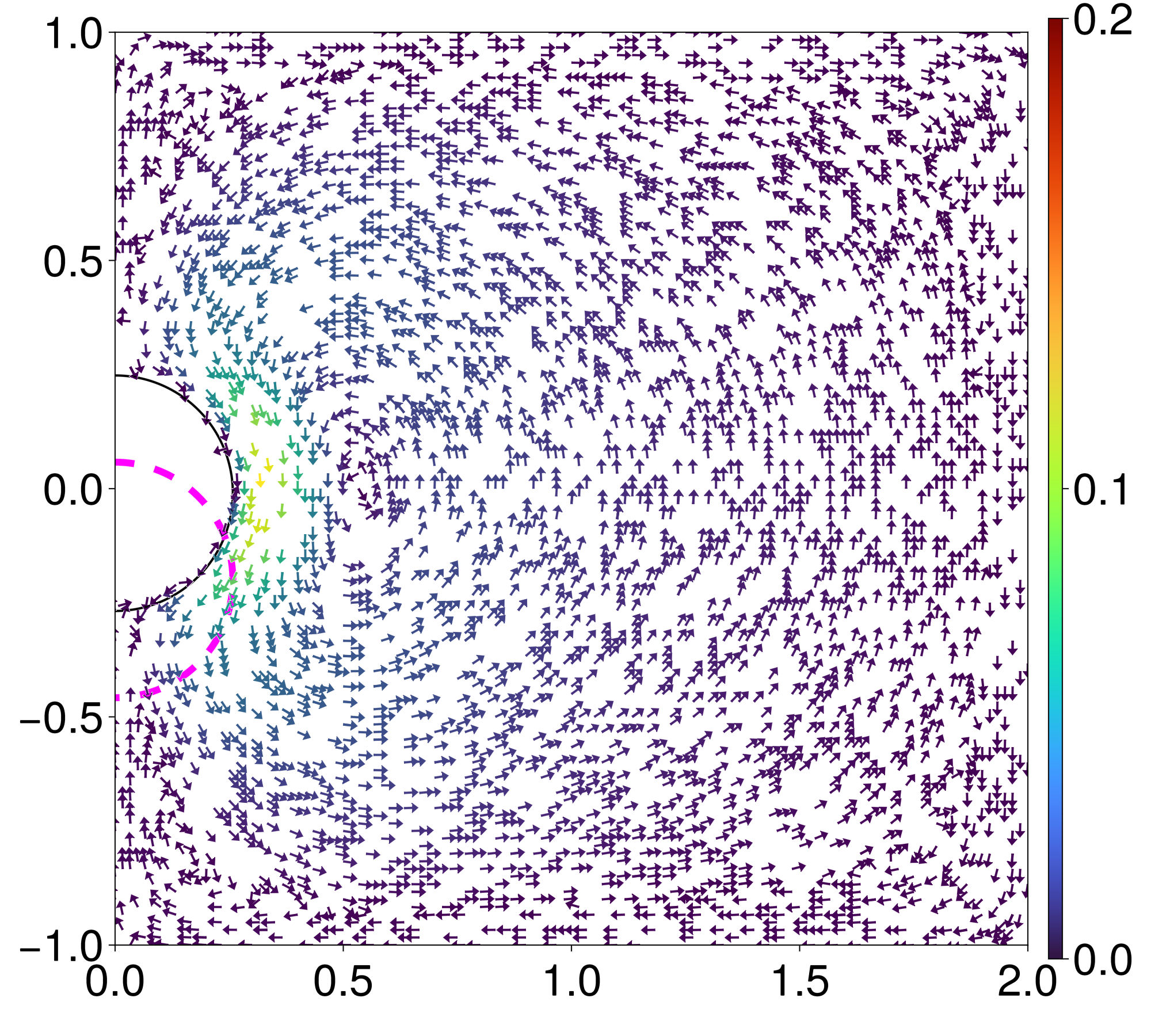}
        \caption{$\nu \cdot t = 0.75$}
   \end{subfigure}
   \hfill
   \begin{subfigure}{0.245\textwidth}
       \includegraphics[width=\textwidth]{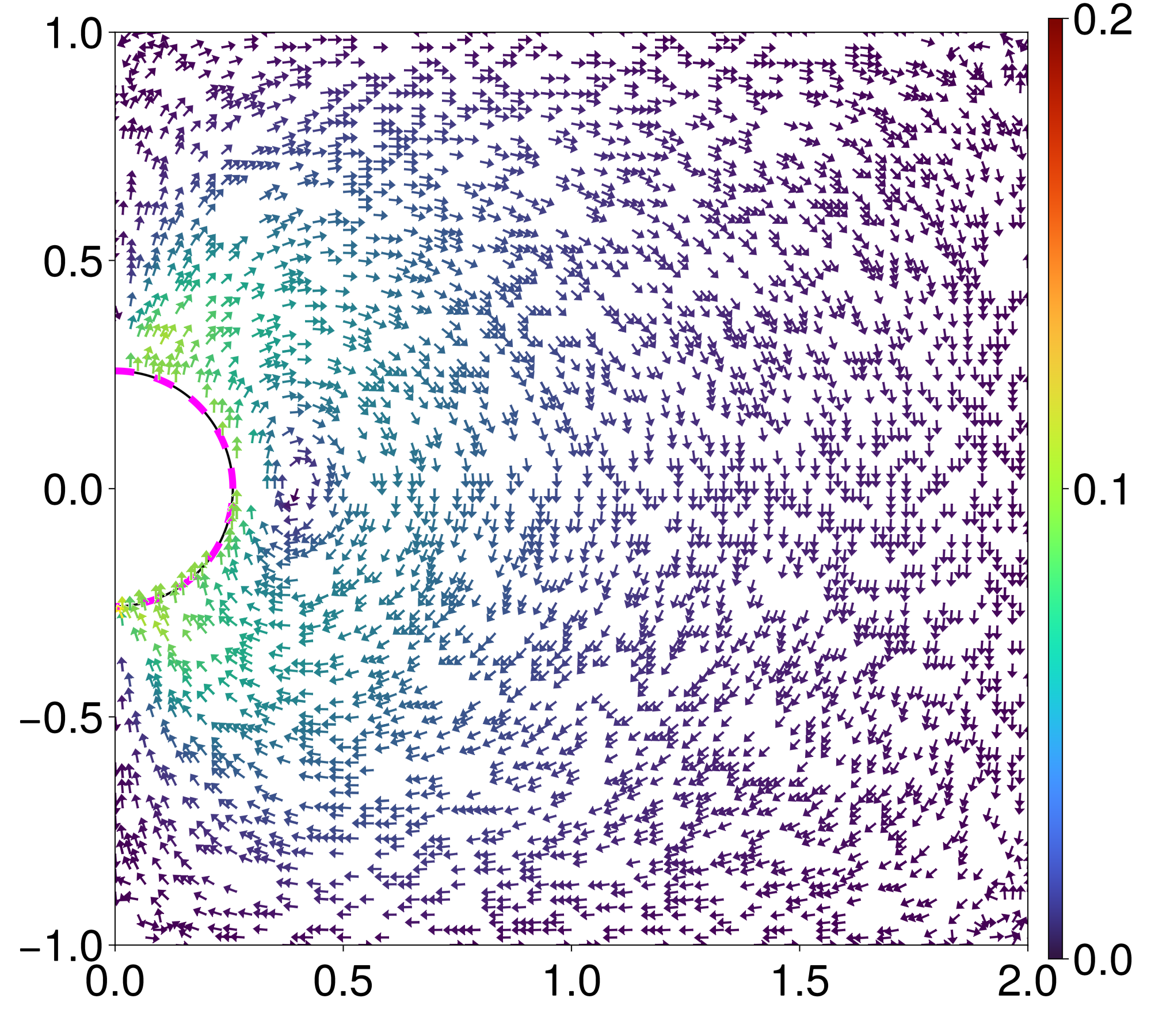}
        \caption{$\nu \cdot t = 1.0$}
   \end{subfigure}
    \caption{Velocity field for the flow past a transversely oscillating bubble at $Re=0.1$ at selected fractions of the first oscillation period. The pink dashed line represents a fictitious depiction of the bubble, where the oscillation amplitude $A$ has been amplified by a factor of 20 for graphical clarity.}
    \label{fig:osc_bubble}
\end{figure}

\begin{figure}
    \centering
   \begin{subfigure}{0.245\textwidth}
       \includegraphics[width=\textwidth]{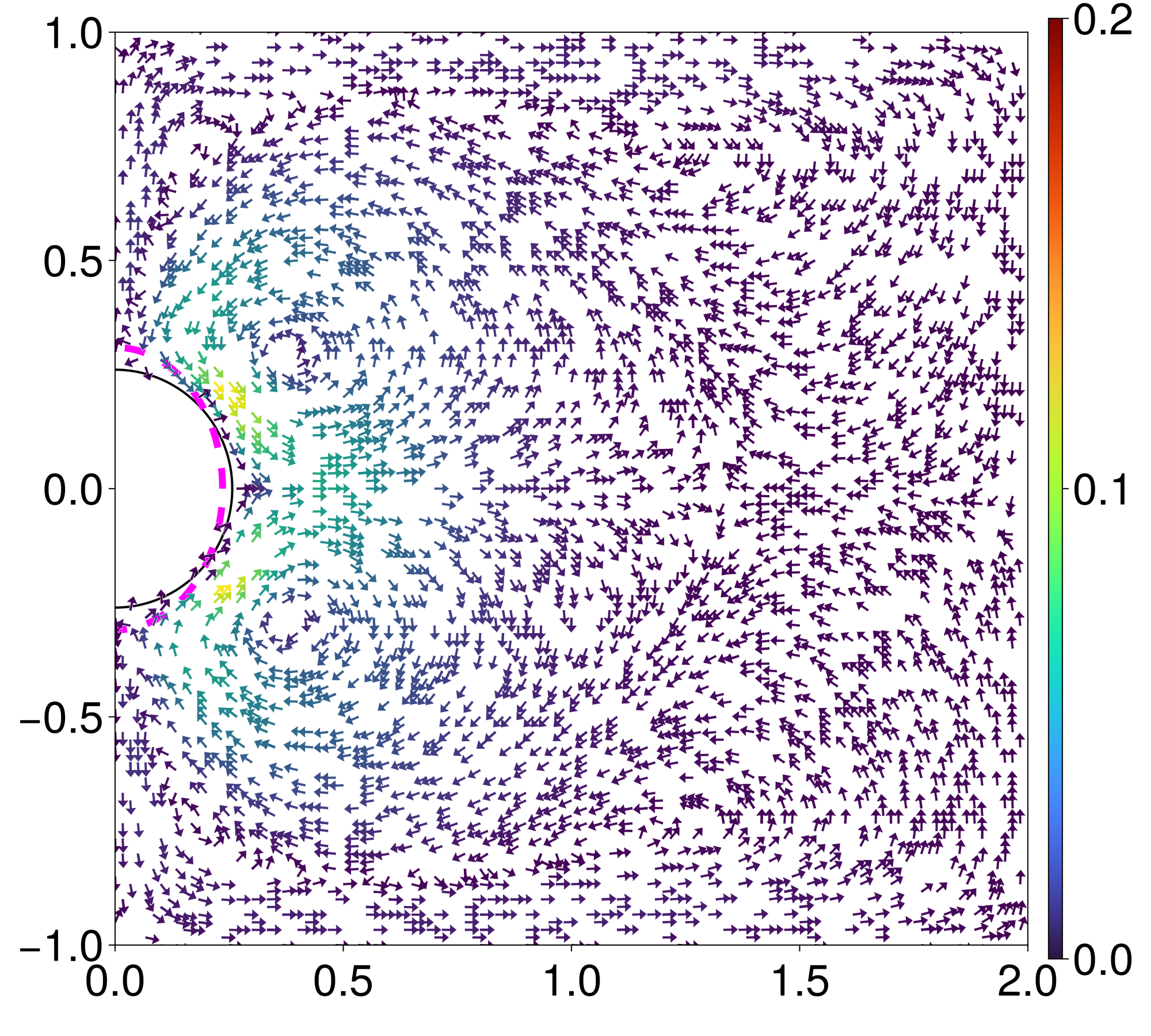}
       \caption{$\nu \cdot t = 0.25$}
   \end{subfigure}
    \hfill 
   \begin{subfigure}{0.245\textwidth}
       \includegraphics[width=\textwidth]{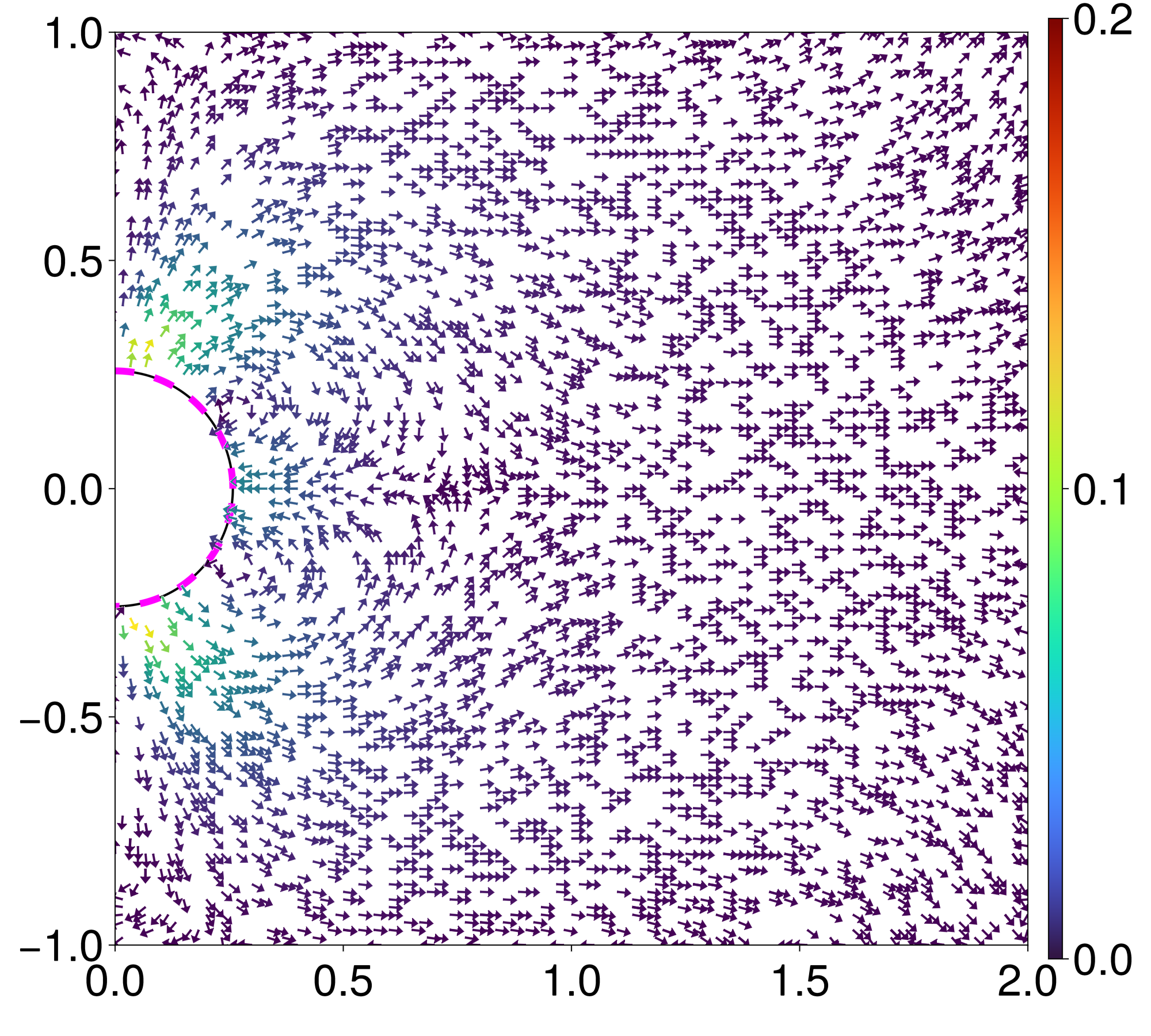}
        \caption{$\nu \cdot t = 0.5$}
   \end{subfigure}
   \hfill
   \begin{subfigure}{0.245\textwidth}
       \includegraphics[width=\textwidth]{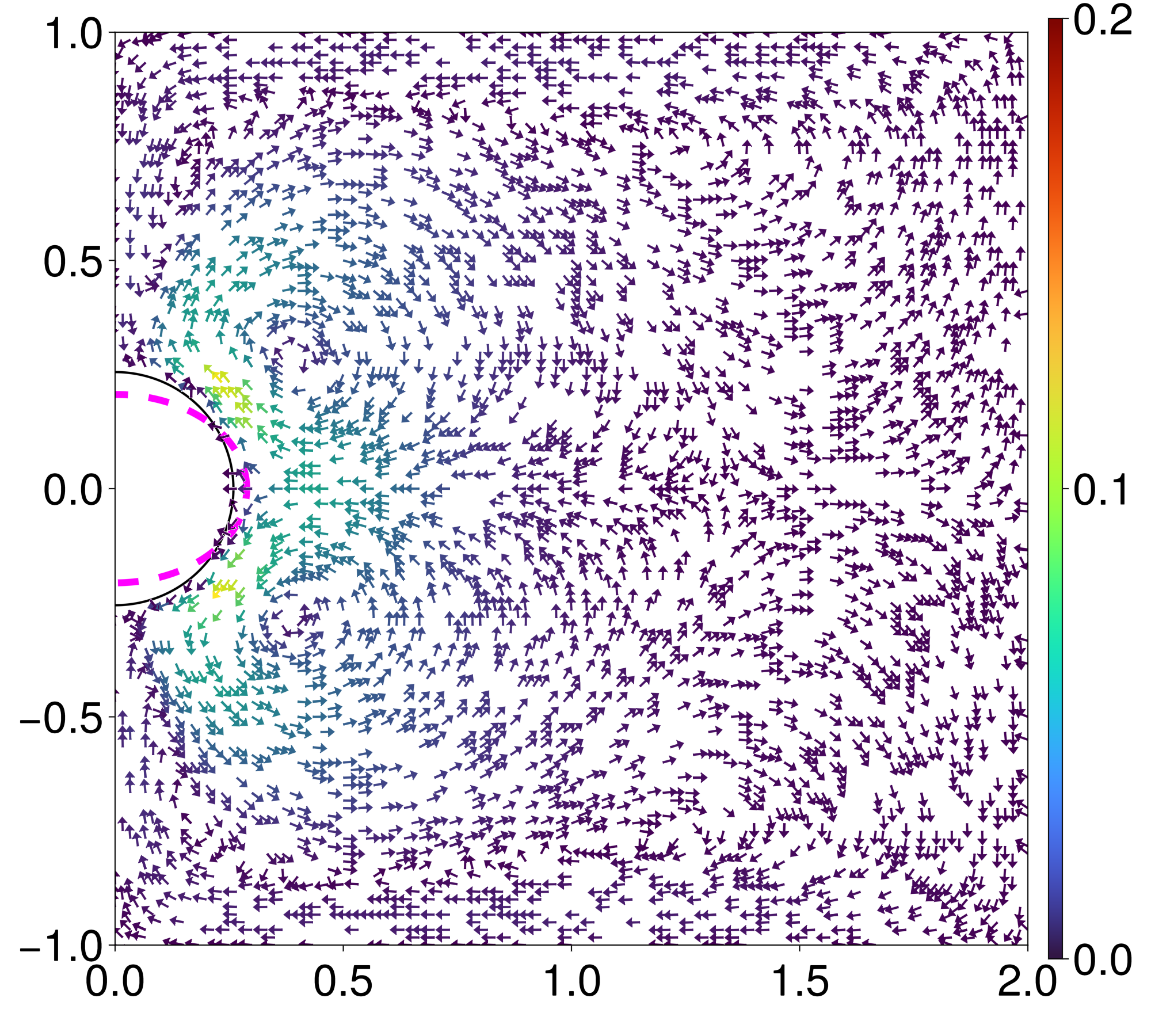}
        \caption{$\nu \cdot t = 0.75$}
   \end{subfigure}
   \hfill
   \begin{subfigure}{0.245\textwidth}
       \includegraphics[width=\textwidth]{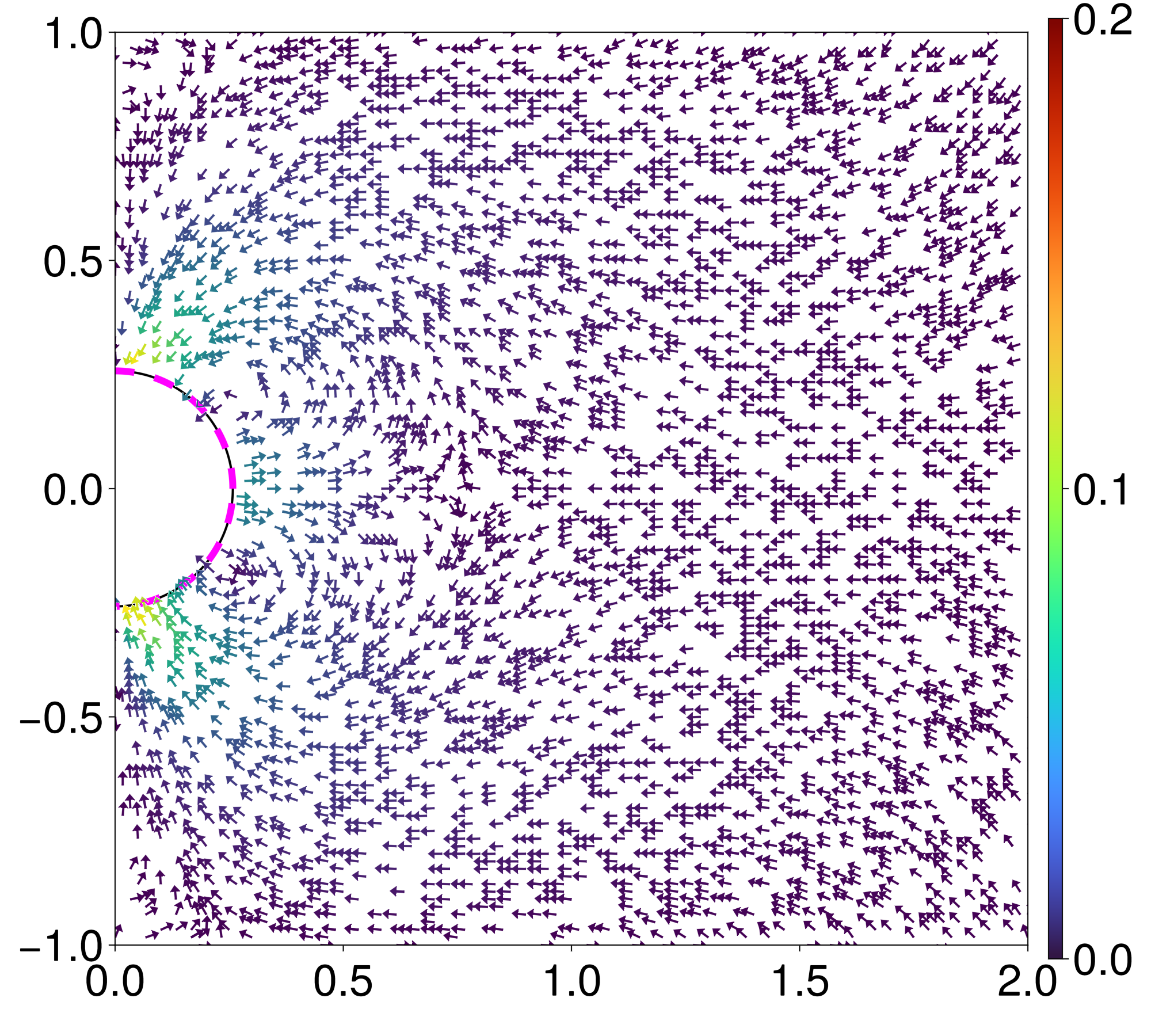}
        \caption{$\nu \cdot t = 1.0$}
   \end{subfigure}
    \caption{Velocity field for the flow past a bubble undergoing ellipsoidal deformation at $Re=0.1$ at selected fractions of the first oscillation period. The pink dashed line represents a fictitious depiction of the bubble, where the oscillation amplitude $A$ has been amplified by a factor of 20 for graphical clarity.}
    \label{fig:ellipsoidal_bubble}
\end{figure}

Motivated by the main application of this work, namely the study of the sorption kinetics of a surfactant around an oscillating bubble, we consider the flow generated by a bubble undergoing small-amplitude oscillations, resulting in a low Reynolds number flow. In the following test, taken from~\cite{ASTUTO2023111880}, we observe the fluid-flow past a transversely oscillating bubble, as well as a bubble undergoing ellipsoidal deformation at $Re = 0.1$.  The motion of the bubble is described using the parametric equations as:
\[\xi(\theta,t) = \xi_c(t) + \delta_{\xi}(t)\cos(\theta), \qquad \zeta(\theta,t)= \zeta_c(t) + \delta_{\zeta}(t)\sin(\theta),\]
where $(\xi_c(t),\zeta_c(t))$ denotes the center of the bubble, and $\delta_{\xi}(t),\delta_{\zeta}(t)$ control its deformation. The velocity $\boldsymbol{u}_B(\xi,\zeta)=(u_B(\xi,\zeta),v_B(\xi,\zeta))$ on the surface of the bubble ($\sqrt{\xi^2 + \zeta^2}=R_B$) is given by:
\[u_B(\xi,\zeta) = \xi_c'(t) + \delta_{\xi}'(t)\cos(\theta), \qquad v_B(\xi,\zeta) = \zeta_c'(t) + \delta_{\zeta}'(t)\sin(\theta),\quad \text{where ~} \theta = \arctan{\zeta/\xi}.\]
In the first test, we consider a vertical oscillation of the bubble:
\[ \xi(t) = 0, \qquad \zeta(t) = A \sin(2\pi \nu t), \qquad \delta_{\xi}(t)= \delta_{\zeta}(t) = R_B. \]
Next, we model the ellipsoidal deformation as:
\[ \xi(t) = \zeta(t) = 0, \qquad \delta_{\zeta}(t)= R_B (1+ A \sin(2\pi \nu t)), \qquad \delta_{\xi}(t) = \sqrt{\frac{R_B^3}{\delta_{\zeta}(t)}}. \]
We choose $A = 0.01,\, \nu = 10$, and the bubble radius $R_B = 0.258$. No-slip boundary conditions are enforced on the walls $(\boldsymbol{u} = 0 )$. The fluid velocity fields at selected fractions of the first oscillation period are shown in Fig.~\ref{fig:osc_bubble} and Fig.~\ref{fig:ellipsoidal_bubble} for the vertical oscillation and ellipsoidal deformation cases, respectively. In these figures, the pink dashed line represents a fictitious depiction of the bubble, where the oscillation amplitude $A$ has been amplified by a factor of 20 for graphical clarity.

\par For the transversely oscillating bubble, at $\nu \cdot t = 0.25$, a vortex is observed at $x\simeq0.5$. By $\nu \cdot t =0.5$, as the bubble returns to its initial configuration, a vortex forms slightly to the left of $x=0.5$. At $\nu \cdot t =0.5$, the vortex moves rightward, and at $\nu \cdot t =1$, when the bubble again returns to its initial configuration, a new vortex forms slightly to the left of $x=0.5$. This sequence repeats, demonstrating the periodic behavior of vortex formation induced by the vertical oscillations of the bubble. 
\par For the ellipsoidal deformation of the bubble, at $\nu \cdot t = 0.25$ and at $\nu \cdot t = 0.75$, two vortices are observed to the left of $x\simeq 0.5$, and two vortices appear between $x = 1$ and $x=1.5$. At $\nu \cdot t = 0.5,1$, the bubble returns to its initial configuration and two new vortices form near the bubble, while the vortices in the region between  $x = 1$ and $x=1.5$ vanish. This behavior demonstrates the periodic formation of vortices in the wake of the ellipsoidally deforming bubble, as observed in experimental results~\cite{2007JFM576191T}.



\section{Conclusions}
In this work we introduced a new numerical strategy for approximating the Navier–Stokes equations in arbitrary moving domains. The method combines a ghost finite element discretization in space with high-order IMEX time integrators. By tailoring Aslam’s extrapolation technique to the finite element framework, we successfully extended the scheme to handle domain motion in a robust and geometrically flexible manner.
We constructed and analyzed IMEX schemes of different orders, showing that the proposed two-stage and four-stage formulations achieve second- and third-order accuracy in time, respectively.
The approximation of the curved boundary by a polygon introduces a second-order geometrical error; higher order spatial discretizations are obtained via the shifted boundary method.
The numerical experiments confirm the expected convergence rates both in fixed and moving domains. In particular, the results demonstrate the effectiveness of Aslam  extrapolation in preserving accuracy when the computational domain evolves in time.

The aim of this research is the application described in the Introduction section, where an oscillating trap adsorbs the diffusing surfactant. At the current stage, we are able to capture only a one-way coupling, where the surrounding fluid is influenced by the motion of the obstacle. The next step is to incorporate the full two-way interaction, allowing the motion of the vibrating bubble to be affected by the fluid dynamics. This will enable us to simulate a scenario in which a domain contains two bubbles: initially, only one bubble oscillates, generating a fluid motion that subsequently induces oscillations in the second bubble through the surrounding flow.

\bibliographystyle{plain}
\bibliography{bibliography.bib}

\end{document}